\documentclass[12pt]{amsart}
\usepackage{amssymb,amscd,amsxtra}

\newcommand{\C}{\mathbb{C}} 
\newcommand{\R}{\mathbb{R}} 
\newcommand{\Z}{\mathbb{Z}} 
\newcommand{\Q}{\mathbb{Q}} 
\newcommand{\F}{\mathbb{F}} 
\newcommand{\cp}{\C P} 
\newcommand{\ch}{\C H} 
\newcommand{\E}{\mathcal{E}} 
\newcommand{\V}{\mathcal{V}} 
\newcommand{\T}{\mathcal{T}} 
\newcommand{\w}{\omega} 
\newcommand{\wbar}{{\bar\w}} 
\newcommand{\thetabar}{{\bar\theta}} 
\newcommand{\Hprim}{H_{0}} 
\newcommand{\hprim}{h_{0}} 
\newcommand{\forms}{\mathcal{C}}
\newcommand{\framed}{\mathcal{F}}
\newcommand{\moduli}{\mathcal{M}}
\renewcommand{\H}{\mathcal{H}}
\newcommand{\per}{g} 
\newcommand{\PGamma}{P\Gamma} 
\newcommand{\D}{\Delta} 
\newcommand{\wt}{\mathop{\rm wt}\nolimits}

\newcommand{\Arg}{\mathop{\rm Arg}\nolimits}
\newcommand{\dimension}{\dim}

\newcommand{\rhosublambda}{\rho_{\lambda^{\vphantom{-1}}}^{\vphantom{\prime}}}
				   
\renewcommand{\a}{\alpha}
\renewcommand{\b}{\beta}
\renewcommand{\c}{\gamma}
\renewcommand{\d}{\delta}
\newcommand{\dlowest}{\d_{\hbox{\rm\scriptsize lowest}}}
\newcommand{\e}{\varepsilon}
\newcommand{\A}{\mathcal{A}}
\newcommand{\B}{\mathcal{B}}
\newcommand{\sset}{\subseteq}
\newcommand{\isomorphism}{\cong}
\newcommand{\tensor}{\otimes}
\newcommand{\aut}{\mathop{\rm Aut}\nolimits}

\newcommand{\Hom}{\mathop{\rm Hom}\nolimits}
\newcommand{\re}{\mathop{\rm Re}\nolimits}
\newcommand{\GL}{{\rm GL}}
\newcommand{\SL}{{\rm SL}}
\newcommand{\PGL}{{\rm PGL}}
\newcommand{\Sp}{{\rm Sp}}
\newcommand{\U}{{\rm U}}
\newcommand{\semidirect}{\rtimes}
\newcommand{\Bhat}{\widehat{B}}
\newcommand{\Dhat}{\widehat{\Delta}}
\DeclareMathOperator{\isom}{\hbox{\rm Isom}}
\DeclareMathOperator{\diag}{\hbox{\rm diag}}

\newcommand{\udot}{{\dot{u}}}
\newcommand{\vdot}{{\dot{v}}}
\newcommand{\wdot}{{\dot{w}}}
\newcommand{\zdot}{{\dot{z}}}
\newcommand{\sdot}{{\dot{s}}}
\newcommand{\I}{{\mathcal{I}}}
\newcommand{\J}{{\mathcal{J}}}
\newcommand{\K}{{\mathcal{K}}}
\newcommand{\What}{\widehat{W}}
\newcommand{\phantominverse}{{\vphantom{-1}}}

\theoremstyle{plain}
\newtheorem{theorem}{Theorem}
\newtheorem{lemma}[theorem]{Lemma}

\numberwithin{theorem}{section}
\numberwithin{equation}{section}
\numberwithin{table}{section}
\numberwithin{figure}{section}

\theoremstyle{remark}
\newtheorem*{remark}{Remark}
\newtheorem*{remarks}{Remarks}

\newcommand{\mapright}[1]{\ \smash{ 
   \mathop{\longrightarrow}\limits^{#1}}\ } 
   
\def\mathllap#1{\mathchoice
{\llap{$\displaystyle #1$}}%
{\llap{$\textstyle #1$}}%
{\llap{$\scriptstyle #1$}}%
{\llap{$\scriptscriptstyle #1$}}}

\def\set#1#2{\left\{\,#1\mathllap{\phantom{#2}}\mathrel{\relax}\right|\left.#2\mathllap{\phantom{#1}}\,\right\}}
\def\ip#1#2{\left\langle#1\mathllap{\phantom{#2}}\right|\left.\!#2\mathllap{\phantom{#1}}\right\rangle}

\begin{document}

\title[Moduli of Cubic Threefolds]{The  Moduli Space of
  Cubic Threefolds as a Ball Quotient}
\author{Daniel Allcock}
\address{Department of Mathematics\\University of Texas at Austin\\Austin, TX 78712}
\email{allcock@math.utexas.edu}
\urladdr{http://www.math.utexas.edu/\textasciitilde allcock}
\author{James A. Carlson}
\address{Department of Mathematics\\
University of Utah\\
Salt Lake City, UT 84112;\\
\lower\smallskipamount\hbox{Clay Mathematics Institute}\\
One Bow Street\\
Cambridge, Massachusetts 02138}
\email{carlson@math.utah.edu; carlson@claymath.org}
\urladdr{http://www.math.utah.edu/\textasciitilde carlson}
\author{Domingo Toledo}
\address{Department of Mathematics\\University of Utah\\Salt
Lake City, UT 84112}
\email{toledo@math.utah.edu}
\urladdr{http://www.math.utah.edu/\textasciitilde toledo}
\date{August 9, 2006}
\thanks{First author partly supported by NSF grants DMS~0070930,
DMS-0231585 and DMS-0245120.
Second and third authors partly supported by NSF grants DMS~9900543
and DMS-0200877.}
\begin{abstract}
The moduli space of cubic threefolds in $\cp^4$, with some minor 
birational modifications, is the Baily-Borel compactification of the
quotient of the complex 10-ball by a discrete group.  We describe both
the birational modifications and the discrete group explicitly.
\end{abstract}
\maketitle
\tableofcontents

\section{Introduction}
\label{sec-intro}

One of the most basic facts in algebraic geometry is that the moduli
space of elliptic curves, which can be realized as plane cubic curves,
is isomorphic to the upper half plane modulo the action of linear
fractional transformations with integer coefficients.  In \cite{ACT},
we showed that there is an analogous result for cubic surfaces; the
analogy is clearest when
we view the upper half plane as complex hyperbolic 1-space, that is,
as the unit disk.  The result is that the moduli space of stable cubic
surfaces is isomorphic to a quotient of complex hyperbolic 4-space by
the action of a specific discrete group.  This is the group of matrices
which preserve the Hermitian form $\diag[-1,1,1,1,1]$ and which have
entries in the ring of Eisenstein integers: the ring obtained by
adjoining a primitive cube root of unity to the integers.  The idea of
the proof is not to use the Hodge structure of the cubic surface,
which has no moduli, but rather that of the cubic threefold obtained
as a triple cover of $\cp^3$ branched along the cubic surface.  The
resulting Hodge structures have a symmetry of order three, and the
moduli space of such structures is isomorphic to complex hyperbolic
4-space $\ch^4$. This is the starting point of the proof, which relies
crucially on the Clemens-Griffiths Torelli Theorem for cubic
threefolds \cite{Clemens-Griffiths}.

The purpose of this article is to extend the analogy to cubic
threefolds.  The idea is to use the period map for the cubic fourfolds
obtained as triple covers of $\cp^4$ branched along the threefolds,
using Voisin's Torelli theorem \cite{voisin} in place of that of Clemens and
Griffiths.  In this case, however, a new phenomenon occurs.  There is
one distinguished point in the moduli space of cubic threefolds which
is a point of indeterminacy for the period map.  This point is the one
represented by what we call a chordal cubic, meaning the
secant variety of a rational normal quartic curve in $\cp^4$.  The reason for
the indeterminacy is that the limit Hodge structure depends on the
direction of approach to the chordal cubic locus.  In fact, the limit
depends \emph{only} on this direction, and so the period map extends to
the blowup of the moduli space.

The natural period map for smooth cubic threefolds
\cite{Clemens-Griffiths} embeds the moduli space in a period domain for Hodge
structures of weight three, namely, a quotient of the Siegel upper
half space of genus five.   For this embedding, however, the target
space has dimension greater than that of the source.  For the
construction of this article, the dimensions of the source and target
are the same. 

To formulate the main result, let $\moduli_{ss}$ be the GIT moduli
space of cubic threefolds, and let $\widehat\moduli_{ss}$ be its blowup at
the point corresponding to the chordal cubics.  Let
$\moduli_s\sset\moduli_{ss}$ be the moduli space of stable cubic
threefolds, and let $\widehat\moduli_s$ be $\widehat\moduli_{ss}$
minus the proper transform of $\moduli_{ss}-\moduli_s$.  
Let $\moduli_0$ be the moduli space of smooth cubic threefolds.
Then we have
the following, contained in the statement of the main result,
theorem~\ref{thm-MAIN-THEOREM}:
\begin{theorem}
There is an arithmetic group $P\Gamma$ acting on complex
hyperbolic 10-space, such that the period map
\[
   \widehat\moduli_{s} \longrightarrow P\Gamma\backslash\ch^{10}
\]
is a isomorphism.  This map identifies $\moduli_s$ with the image in
$P\Gamma\backslash\ch^{10}$ of the complement of a hyperplane
arrangement $\H_c$.  It also identifies the discriminant in
$\widehat\moduli_{s}$ with the image of another hyperplane
arrangement $\H_{\D}\sset\ch^{10}$.  In particular, it identifies
$\moduli_0$ with
$P\Gamma\backslash\bigl(\ch^{10}-(\H_c\cup\H_{\D})\bigr)$.  Finally, the
period map extends to a morphism from $\widehat\moduli_{ss}$ to the
Baily-Borel compactification $\overline{\PGamma\backslash\ch^{10}}$.
\end{theorem}

\noindent
We also provide much more detailed information about all the objects
in the theorem, such as explicit descriptions of $P\Gamma$, $\H_c$ and
$\H_{\D}$, and an analysis of the part of $\widehat\moduli_{ss}$ lying
over the boundary points of $\overline{P\Gamma\backslash\ch^{10}}$.

Now we will say what the group $P\Gamma$ is.  Let $V$ be a cyclic
cubic fourfold, meaning a triple cover of $P^4$ branched over a cubic
threefold.  The primitive cohomology $\Lambda(V)$ of $V$ is naturally
a module for the Eisenstein integers, where a primitive cube root of
unity acts on cohomology as does the corresponding deck
transformation.  When $V$ is smooth, This Eisenstein module carries a natural %
nondegenerate Hermitian form of signature $(10,1)$, and $P\Gamma$ is
the projective isometry group of $\Lambda(V)$.

The architecture of our proof of theorem~\ref{thm-MAIN-THEOREM} dates back to before
\cite{allcock-threefolds}, and follows the pattern laid out in \cite{ACT} for cubic
surfaces.  But it is considerably more technical in its details.  We
therefore focus on the points where there are major differences or
where substantially more work must be done as compared with the case
of cubic surfaces.

In section~\ref{sec-smoothmoduli} we establish basic facts about $\Lambda(V)$ as an
Eisenstein module endowed with a complex Hodge structure, give an
inner product matrix for $\Lambda(V)$, and show that 
$\moduli_0\to\PGamma\backslash\ch^{10}$ is an isomorphism onto its
image.  The argument here follows that of \cite{ACT}, except that in
place of the Clemens-Griffiths theorem, we use Claire Voisin's Torelli
theorem for cubic fourfolds \cite{voisin}.  In this section we also
establish facts about the discriminant locus for cubic threefolds near
stable singularities.  These facts are used in section~\ref{sec-extension} for
extending the period map.

In section~\ref{sec-discr-near-chordal-cubic} we blow up the space $P^{34}$ of cubic threefolds
along the chordal cubic locus, and describe the proper transform of
the discriminant locus.  This is one of the most technical sections in
the paper, but it is required for the extension of the period map in
section~\ref{sec-extension}.  To give an idea of the main issue, consider a
one-parameter family of smooth cubic threefolds degenerating to a
chordal cubic.  We may write it as $F + tG=0$, where $F = 0$ defines
the chordal cubic.  The polynomial $G = 0$ cuts out on the rational
normal curve of $F = 0$ a set of twelve points.  Thus to a tangent
vector at a point of the chordal locus one associates a 12-tuple on
the projective line.  We show that the discriminant locus in the
blown-up $P^{34}$ has a local product structure, e.g., it is
homeomorphic to the product of the discriminant locus for 12-tuples in
$P^1$, times a disk representing the transverse direction, times a
twenty-one dimensional ball corresponding to the action of the
projective group $\PGL(5,\C)$.  
There is also a technical variation on this result which gives an
analytic model of the discriminant, not just a topological one.

In section~\ref{sec-extension} we extend the period map to the
semistable locus of the blown-up $P^{34}$.  This requires some
geometric invariant theory to say what the semistable locus is, and
here the work of Reichstein \cite{reichstein} is essential.  Then we study the
local monodromy groups at points in this semistable locus.  The
essential result for the extension of the period map is that these
groups are all finite or virtually unipotent.

In section~\ref{sec-hodge-theory-chordal} we show that the extended
 period map sends the chordal cubic locus to a divisor in
 $\PGamma\backslash\ch^{10}$.  The main point here is to identify the
 limit Hodge structure and in so doing show that the derivative of the
 extended map along the blowup of the chordal cubic locus is of rank
 nine.  We do this by establishing an isogeny between the limit Hodge
 structure and the sum of a 1-dimensional Hodge structure and a Hodge
 structure associated to a six-fold cover of the projective line
 branched at twelve points.  These Hodge structures first arose in the
 work of Deligne and Mostow \cite{Deligne-Mostow}.  Our analysis shows
 that the image of the period map is the quotient of a totally
 geodesic $\ch^9$ by a suitable subgroup of $\Gamma$.

Section~\ref{sec-hodge-theory-nodal} deals with the same issues for
 the divisor of nodal cubic threefolds.  Here the analysis is easier.
 We show that the Hodge structure is isomorphic to that of a special
 K3 surface, plus a 1-dimensional summand.  This K3 surface and its
 Hodge structure were treated in Kond\=o's work \cite{kondo:genus-4}
 on moduli of genus~4 curves.

In section~\ref{sec-main-theorem} we assemble the various pieces to
prove the main theorem, and in section~\ref{sec-group-and-arrangement}
we give some supplemental results on the monodromy group and the
hyperplane arrangements.

Another proof of the main theorem has been obtained by Looijenga and
Swierstra \cite{looijenga-period} .  Both proofs proceed by extending
the period map from the moduli space of smooth threefolds to a larger
space, but the extension process is quite different in the two proofs.
We use a detailed analysis of the discriminant in the space obtained
by blowing up the chordal cubic locus to extend the period map.
Looijenga and Swierstra use a general machinery developed earlier by
them \cite{looijenga-open} to handle extensions of period mappings.
We are grateful to Looijenga for sending us an early version of their
argument.
 
It is a pleasure to thank the many people whose conversations have
played a helpful role in the long gestation of this paper, including
Herb Clemens, Alessio Corti,  Johan de Jong, J\'anos Koll\'ar, 
Eduard Looijenga and Madhav Nori.  We are grateful to the Clay
Mathematics Institute for its hospitality.


\section{Moduli of smooth cubic threefolds}
\label{sec-smoothmoduli}

This section contains a number of foundational results, and its main
theorem is of interest in its own right.  We consider cyclic cubic
fourfolds, i.e., triple covers of $\cp^4$ branched along cubic
threefolds.  (1) The cohomology of such a fourfold is a module over
$\E=\Z[\w{=}{\root3\of1}]$, equipped with a Hermitian form.  (2) The
monodromy on this lattice as the threefolds acquire a node is a
complex reflection of order three; see Lemma
\ref{lem-meridians-act-by-complex-reflections}.  (3) To analyze the local
monodromy for more complicated singularities, we give a structure
theorem for the discriminant locus of the space of cubic threefolds
near a cubic threefold with singularities of type $A_n$ and
$D_4$.  See Lemma \ref{lem-discriminant-away-from-E}. Using this
result, we give an inner product matrix for the Hermitian form; see
Lemma \ref{thm-inner-product-matrix}.  (4) With the previous results
in hand, we define a framing of the Hodge structure of a cyclic cubic
fourfold and define the period map.  Finally, the main theorem of the
section is that the period
map for smooth cubic threefolds is an isomorphism onto its image;  see
Theorem \ref{thm-main-theorem-smooth-case}.

Let $\forms$ be the space of all nonzero cubic forms in variables
$x_0,\dots,x_4$.  For such a form $F$ let $T$ be the cubic threefold
in $\cp^4$ it defines, and let $V$ be the cubic fourfold in
$\cp^5$ defined by $F(x_0,\dots,x_4)+x_5^3=0$.  
Whenever we consider
$F\in\forms$, $T$ and $V$ will have these meanings.  $V$ is the triple
cover of $\cp^4$ branched along $T$.  We write $\forms_0$ for the set of
$F\in\forms$ for which $T$ is smooth (as a scheme) and $\D$ for the
discriminant $\forms-\forms_0$.  We will sometimes also write $\D$ for
its image in $P\forms$; context will make our meaning clear.
We write $\forms_s$ for the set of 
$F\in\forms$ for which $T$ is stable in the sense of geometric
invariant theory.  By \cite{allcock-threefolds} or \cite{yokoyama}, 
this holds if and only if $T$ has no
singularities of types other than $A_1,\dots,A_4$.  $\forms_s$ will
play a major role in sections~\ref{sec-extension}--\ref{sec-main-theorem}; in this section all we will use
is the fact that  $\forms_0$ lies within it.

Because we will vary our threefolds, we will need the universal
family $\T\sset\forms\times\cp^4$ of cubic threefolds,
$$
\T=\set{(F,[x_0{:}\dots{:}x_4])\in\forms\times\cp^4}{F(x_0,\dots,x_4)=0}\;,
$$
and the family of covers of $\cp^4$ branched over them,
$$
\V=\set{(F,[x_0{:}\dots{:}x_5])\in\forms\times\cp^5}{F(x_0,\dots,x_4)+x_5^3=0}\;.
$$
We will write $\pi_\T$ and $\pi_\V$ for the projections $\T\to\forms$
and $\V\to\forms$.
The total spaces of $\T$ and $\V$ are smooth.  We write
$\T_0$ and $\V_0$ for the topologically locally trivial fibrations
which are the restrictions of $\T$ and $\V$ to $\forms_0$.
The transformation
$\sigma:\C^6\to\C^6$ given by
\begin{equation}
\label{eq-def-of-sigma}
\sigma(x_0,\dots,x_5)=(x_0,\dots,x_4,\w x_5)\;, 
\end{equation}
where $\w$ is a fixed
primitive cube root of unity, plays an essential role in all that
follows.  We regard it as a symmetry of $\V$ and of the individual
$V$'s.  

Our period map $\forms_0\to\ch^{10}$ will be defined using the Hodge
structure of the fourfolds and its interaction with $\sigma$, so we need
to discuss $H^4(V)$ for $F\in\forms_0$.  To compute this it suffices
by the local triviality of $\V_0$ to consider the single
fourfold $x_0^3+\dots+x_5^3=0$; by thinking of it as an iterated
branched cover, one finds that its Euler characteristic is~$27$.  The
Lefschetz hyperplane theorem and Poincar\'e duality imply that
$H^i(V;\Z)$ is the same as $H^i(\cp^5;\Z)$ for $i\neq4$, so
$H^4(V;\Z)\isomorphism\Z^{23}$.   The class of a 3-plane
in $\cp^5$ pulls back to a class $\eta(V)\in H^4(V;\Z)$ of norm~$3$.
The primitive cohomology $\Hprim^4(V;\Z)$ is the orthogonal complement
of $\eta(V)$ in $H^4(V;\Z)$.  Since $H^4(V;\Z)$ is a unimodular
lattice, $\Hprim^4(V;\Z)$ is a $22$-dimensional lattice with
determinant equal to that of $\langle\eta(V)\rangle$, up to a sign, so
$\det\Hprim^4(V;\Z)=\pm3$. 

$\Hprim^4(V;\Z)$ is a module not only over $\Z$ but over the
Eisenstein integers $\E=\Z[\w]$ as well.  
To see this, observe that
the isomorphism $V/\langle\sigma\rangle\isomorphism \cp^4$ implies that
$H^4(\cp^4;\C)$ is the $\sigma$-invariant part of $H^4(V;\C)$.
Therefore $\sigma$ fixes no element of $H^4_0(V;\C)$ except $0$, hence
no element of $H^4_0(V;\Z)$ except $0$.  We define
$\Lambda(V)$ to be the $\E$-module whose underlying additive group is
$\Hprim^4(V;\Z)$, with the action of $\w\in\E$ defined as $\sigma^*$.
$\E$ is a unique factorization domain, so $\Lambda(V)$ is free of
rank~$11$.  We define a Hermitian form on $\Lambda(V)$ by the formula
\begin{equation}
\label{eq-defn-of-hermitian-form}
\ip{\a}{\b}=\frac{1}{2}\Bigl[
3\a\cdot\b-\theta\a\cdot({\sigma^*}^{-1}\b-\sigma^*\b)
\Bigr]\;,
\end{equation}
where the dot denotes the usual pairing $\a\cdot\b=\int_V\a\wedge\b$ and
$\theta=\w-\wbar=\sqrt{-3}$. 
The scale factor $1/2$ is chosen so that $\ip{\a}{\b}$ takes values in
$\E$; it is the smallest scale for which this is true.

\begin{lemma}
\label{lem-basic-properties-of-hermitian-pairing}
$\ip{{\cdot}}{{\cdot}}$ is an $\E$-valued Hermitian form, linear in
its first argument and antilinear in its second.  Furthermore,
$\ip{\a}{\b}\in\theta\E$ for all $\a,\b\in \Lambda(V)$.
\end{lemma}

\begin{proof}
$\Z$-bilinearity is obvious.  $\E$-linearity in its first argument
holds by the following calculation.  (Throughout the proof we write
$\sigma$ for $\sigma^*$.)
\begin{equation*}
\begin{split}
\ip{\theta\a}{\b}
&=\ip{\sigma\a-\sigma^{-1}\a}{\b}\\
&=\frac{1}{2}\Bigl[
3(\sigma\a-\sigma^{-1}\a)\cdot\b-\theta(\sigma\a-\sigma^{-1}\a)\cdot(\sigma^{-1}\b-\sigma\b)\Bigr]\\
&=\frac{\theta}{2}\Bigl[
\thetabar\sigma\a\cdot\b-\thetabar\sigma^{-1}\a\cdot\b\\
&\mathrel{\phantom{=}}\phantom{\frac{\theta}{2}\Bigl[}
-(\sigma\a\cdot\sigma^{-1}\b-\sigma\a\cdot\sigma\b
-\sigma^{-1}\a\cdot\sigma^{-1}\b+\sigma^{-1}\a\cdot\sigma\b)\Bigr]\\
&=\frac{\theta}{2}\Bigl[
\thetabar\a\cdot\sigma^{-1}\b-\thetabar\a\cdot\sigma\b\\
&\mathrel{\phantom{=}}\phantom{\frac{\theta}{2}\Bigl[}
-(\a\cdot\sigma^{-2}\b-\a\cdot\b-\a\cdot\b+\a\cdot\sigma^2\b)\Bigr]\\
&=\frac{\theta}{2}\Bigl[
2\a\cdot\b-\a\cdot(\sigma^{-2}\b+\sigma^2\b)
-\theta(\a\cdot\sigma^{-1}\b-\a\cdot\sigma\b)\Bigr]\\
&=\frac{\theta}{2}\Bigl[
3\a\cdot\b-\theta\a\cdot(\sigma^{-1}\b-\sigma\b)\Bigr]\\
&=\theta\ip{\a}{\b}\;.
\end{split}
\end{equation*}
In the second-to-last step we used the relation
$\sigma^{-2}+\sigma^2=\sigma+\sigma^{-1}=-1$.  That
$\ip{{\cdot}}{{\cdot}}$ is a $\C$-valued Hermitian form now follows
from $\ip{\a}{\b}=\overline{\ip{\b}{\a}}$; to prove this it suffices
to check that the imaginary part of \eqref{eq-defn-of-hermitian-form} changes sign when $\a$
and $\b$ are exchanged, i.e.,
\begin{equation*}
\begin{split}
\a\cdot(\sigma^{-1}\b-\sigma\b)
&=\a\cdot\sigma^{-1}\b-\a\cdot\sigma\b\\
&=\sigma\a\cdot\b-\sigma^{-1}\a\cdot\b\\
&=-(\sigma^{-1}\a-\sigma\a)\cdot\b\\
&=-\b\cdot(\sigma^{-1}\a-\sigma\a)\;.
\end{split}
\end{equation*}
Next we check that $\ip{{\cdot}}{{\cdot}}$ is $\E$-valued.  It is
obvious that its real part takes values in $\frac{1}{2}\Z$ and that
its imaginary part takes values in $\frac{\theta}{2}\Z$.  Since
$$
\E=\set{a/2+b\theta/2}{\hbox{$a,b\in\Z$ and $a\equiv b$ mod 2}}\;,
$$
it suffices to prove that 
$$
3\a\cdot\b\equiv -\a\cdot(\sigma^{-1}\b-\sigma\b)\qquad\hbox{(mod 2)},
$$
that is, that $2$ divides $\a\cdot(\b-\sigma\b+\sigma^{-1}\b)$.  This
follows from the relation $1-\w+\wbar=-2\w$ in $\E$.
Furthermore, \eqref{eq-defn-of-hermitian-form} shows that $\ip{\a}{\b}$ has real part in
$\frac{3}{2}\Z$; since every element of $\E$ with real part in
$\frac{3}{2}\Z$ is divisible by $\theta$, the proof is complete.
\end{proof}

We will describe $\Lambda(V)$ more precisely later, after considering
the Hodge structure of $V$.  

By the Griffiths residue calculus, $\Hprim^{p,q}(V;\C)$ is spanned for
$p+q=4$ by the residues of rational differential forms
\begin{equation}
\label{eq-Omega(A)}
\Omega(A)=\frac{A(x_0,\dots,x_5)\Omega}{\bigl(F(x_0,\dots,x_4)+x_5^3\bigr)^{q+1}}
\end{equation}
where 
$$
\Omega=\sum_{j=0}^5(-1)^jX_j\,
dX_0\wedge\dots\wedge\widehat{dX_j}\wedge\dots\wedge dX_5
$$ has degree~6 and $A$ is any homogeneous polynomial of degree
$3q-3$, so that the total degree of $\Omega(A)$ is~0.  Two such
polynomials define the same cohomology class, modulo
$\oplus_{p'>p}H^{p',4-p'}_0$, if and only if $A-A'$ lies in the
Jacobian ideal of $F+x_5^3$.  This gives Hodge numbers
$\hprim^{4,0}=\hprim^{0,4}=0$, $\hprim^{3,1}=\hprim^{1,3}=1$ and
$\hprim^{2,2}=20$.  Since $\sigma$ is holomorphic, its eigenspace
decomposition refines the Hodge decomposition.  In particular, the
generator $\Omega(1)=\Omega/(F+x_5^3)^2$ of $H^{3,1}$ has eigenvalue
$\w$.  Therefore $H^4_\w(V;\C)=H^{3,1}_\w\oplus H^{2,2}_\w$, the
summands being one- and ten-dimensional.

On $H^4(V;\C)$ there is the Hodge-theoretic Hermitian pairing
\begin{equation}
\label{eq-defn-of-Hodge-theoretic-pairing}
(\a,\b)=3\int_V \a\wedge\bar\b\;.
\end{equation}
The Hodge-Riemann bilinear relations \cite[p.~123]{griffiths-harris}
show that  $({\cdot},{\cdot})$
is positive-definite on $H_0^{2,2}$ and negative-definite on $H^{3,1}$.
It follows that $H^4_\w(V;\C)$ has signature $(10,1)$.

If $W$ is a complex vector space of dimension $n+1$ with a Hermitian
form of signature $(n,1)$ then we write $\ch(W)$ for the space of
lines in $W$ on which the given form is negative definite.  We call
this the complex hyperbolic space of $W$; it is an open subset of $PW$
and is biholomorphic to the unit ball in $\C^n$.  The previous two
paragraphs may be summarized by saying that the Hodge structure of $V$
defines a point in $\ch\bigl(H^4_\w(V;\C)\bigr)$.

We chose the factor $3$ in \eqref{eq-defn-of-Hodge-theoretic-pairing} so that $({\cdot},{\cdot})$ and
$\ip{{\cdot}}{{\cdot}}$ would agree in the sense of the next lemma.
To relate
the two Hermitian forms we consider the $\R$-linear map
$\Hprim^4(V;\R)\to H^4_\w(V;\C)$ which is the 
inclusion 
$$
H^4_0(V;\R)\to H^4(V;\R)\to H^4(V;\C)
$$ 
followed by projection to $\sigma$'s $\w$-eigenspace.  This is an
isomorphism of real vector spaces.  Since the $\Z$-lattice underlying
$\Lambda(V)$ is $\Hprim^4(V;\Z)\sset\Hprim^4(V;\R)$, we get a map
$Z:\Lambda(V)\tensor_\Z\R\to H^4_\w(V;\C)$.    Since
the complex structure on $\Lambda(V)\tensor_\Z\R$ is defined by taking
$\w$ to act as $\sigma^*$, and since $H^4_\w$ is defined as a space on
which $\sigma^*$ acts by multiplication by $\w$, $Z$ is
complex-linear.  Since
$\Lambda(V)\tensor_\Z\R=\Lambda(V)\tensor_\E\C$, we may regard $Z$ as
an  isomorphism $\Lambda(V)\tensor_\E\C\to
H^4_\w(V;\C)$ of complex vector spaces.

\begin{lemma}
\label{lem-Hermitian-forms-agree}
For all $\a,\b\in \Lambda(V)$, $\bigl(Z(\a),Z(\b)\bigr)=\ip{\a}{\b}$.
\end{lemma}

\begin{proof}
Since both $({\cdot},{\cdot})$ and $\ip{{\cdot}}{{\cdot}}$ 
are Hermitian forms, it suffices to check that
$(Z\a,Z\a)=\ip{\a}{\a}$ for all $\a$.  We write $\a_\w$
and $\a_\wbar$ for the projections of $\a\in\Hprim^4(V;\R)$ to
$H^4_\w(V;\C)$ and $H^4_\wbar(V;\C)$. By definition,
$$
(Z\a,Z\a)=3\int_V \a_\w\wedge\overline{\a_\w}
=3\int_V\a_\w\wedge\a_\wbar\;.
$$
We will write $\sigma$ for $\sigma^*$ throughout the proof.
Since $\sigma\a=\w\a_\w+\wbar\a_\wbar$ and
$\sigma^{-1}\a=\wbar\a_\w+\w\a_\wbar$, we deduce
$\a_\w=-\frac{1}{\theta}(\w\sigma\a-\wbar\sigma^{-1}\a)$ and
$\a_\wbar=-\frac{1}{\theta}(-\wbar\sigma\a+\w\sigma^{-1}\a)$.
Therefore
\begin{equation*}
\begin{split}
(Z\a,Z\a)
&=3\int \frac{1}{\theta^2}\Bigl[
\w\sigma\a\wedge(-\wbar\sigma\a)+\w\sigma\a\wedge\w\sigma^{-1}\a\\
&\mathrel{\phantom{=}}\phantom{3\int \frac{1}{\theta^2}\Bigl[}
-\wbar\sigma^{-1}\a\wedge(-\wbar\sigma\a)
-\wbar\sigma^{-1}\a\wedge\w\sigma^{-1}\a\Bigr]\\
&=-\int\Bigl[
-\sigma\a\wedge\sigma\a+\wbar\sigma\a\wedge\sigma^{-1}\a\\
&\mathrel{\phantom{=}}\phantom{-\int\Bigl[}
+\w\sigma^{-1}\a\wedge\sigma\a-\sigma^{-1}\a\wedge\sigma^{-1}\a\Bigr]\\
&=-\int\Bigl[-2\a\wedge\a-\sigma\a\wedge\sigma^{-1}\a\Bigr]\\
&=2\a\cdot\a+\int\sigma\a\wedge\sigma^{-1}\a\;.
\end{split}
\end{equation*}
To evaluate the second term, we use
\begin{equation}
\label{eq-foo}
\int \sigma\a\wedge\a+\int\sigma\a\wedge\sigma\a+\int\sigma\a\wedge\sigma^{-1}\a=0\;,
\end{equation}
which follows from 
$1+\sigma+\sigma^{-1}=0$.
Since $\sigma$ is an isomorphism and $\wedge$ is symmetric, the first
and last terms of \eqref{eq-foo} are equal, so we get 
$$
\int\sigma\a\wedge\sigma^{-1}\a
=-\frac{1}{2}\int\sigma\a\wedge\sigma\a
=-\frac{1}{2}\a\cdot\a\;.
$$
This yields
$$
(Z\a,Z\a)=\frac{3}{2}\a\cdot\a=\ip{\a}{\a}
$$
as desired.
\end{proof}

 From the lemma it follows that $\Lambda(V)$ has signature $(10,1)$ and
that $\ch(H^4_\w(V;\C))$ is naturally identified with
$\ch(\Lambda(V)\tensor_\E\C)$.  We write $\ch(V)$ for either of these
complex hyperbolic spaces.

In order to study the variation of the Hodge structure of $V$ we must
realize our constructions in local systems over $\forms_0$.  To do
this we use the fact that the family $\V_0$ over $\forms_0$ gives rise
to a sheaf $R^4\pi_*(\Z)$ over $\forms_0$.  Recall that this is the
sheaf associated to the presheaf $U\mapsto H^4(\pi_\V^{-1}(U);\Z)$.
Since $\V_0$ is topologically locally trivial, $R^4\pi_*(\Z)$ is a
local system of 23-dimensional $\Z$-lattices isomorphic to
$H^4(V;\Z)$.  The map $\eta:F\mapsto\eta(V)$, for $F\in\forms_0$, is a
section over $\forms_0$, and the subsheaf
$(R^4\pi_*(\Z))_0$, whose local sections are the local
sections of $R^4\pi_*(\Z)$ orthogonal to $\eta$, is a local system of
22-dimensional $\Z$-lattices isomorphic to $\Hprim^4(V;\Z)$.  Since
$\sigma$ acts  on $\V_0$, it acts on $R^4\pi_*(\Z)$; since it
preserves $\eta$ it acts on $(R^4\pi_*(\Z))_0$, giving the sheaf
the structure of a local system of $\E$-modules isomorphic to
$\Lambda(V)$.  We call this local system $\Lambda(\V_0)$.  The formula
\eqref{eq-defn-of-hermitian-form} endows $\Lambda(\V_0)$ with the
structure of a local system of Hermitian $\E$-modules.  We write
$\ch(\V_0)$ for the corresponding local system of hyperbolic spaces.

We can also consider the sheaf $R^4\pi_*(\C)$ over $\forms_0$; it is a
local system because it is the complexification of $R^4\pi_*(\Z)$.
Now, $\sigma$ acts on $R^4\pi_*(\C)$ and we consider its
$\w$-eigensheaf $(R^4\pi_*(\C))_\w$, which is a local system of
Hermitian vector spaces isometric to $H^4_\w(V;\C)$, hence of
signature $(10,1)$, with a corresponding local system of complex
hyperbolic spaces.  The map $Z$ identifies the local systems
$\Lambda(\V_0)\tensor_\E\C$ and $(R^4\pi_*(\C))_\w$, and therefore
identifies the two local systems of complex hyperbolic spaces.
Therefore we may regard the inclusion $H^{3,1}(V)\to H^4_\w(V)$ as a
defining a section
\begin{equation}
\label{eq-per-map-as-section}
\per:\forms_0\to\ch(\V_0)\;.
\end{equation}
It is holomorphic since the Hodge filtration varies holomorphically.
This is the period map; all our results refer to various formulations
of it.

Next we obtain a concrete description of the $\E$-lattice
$\Lambda(V)$, by
investigating the monodromy of $\Lambda(\V_0)$.
Fix a basepoint $F\in\forms_0$, let $\Gamma(V)$ be the isometry group
of $\Lambda(V)$, and let $\rho$ be the monodromy
representation $\pi_1(\forms_0,F)\to\Gamma(V)$.
By a meridian around a divisor, such as $\D$,
we mean the boundary circle of a small disk transverse
to the divisor at a generic point of it, traversed once positively.

If $W$ is a complex vector space then a complex reflection of $W$ is a
linear transformation that fixes a hyperplane pointwise and has
finite order${}>1$.  If this order is $2$, $3$ or $6$ then it is
called a biflection, triflection or hexaflection.  If $W$ has a
Hermitian form $\ip{\cdot}{\cdot}$, $r\in W$ has nonzero norm and
$\zeta$ is a primitive $n$th root of unity, then the transformation
\begin{equation}
\label{eq-defn-of-complex-reflection}
x\mapsto x-(1-\zeta)\frac{\ip{x}{r}}{\ip{r}{r}}r
\end{equation}
is a complex reflection of order $n$ and preserves
$\ip{{\cdot}}{{\cdot}}$.  It fixes $r^\perp$ pointwise and sends $r$
to $\zeta r$; we call
it the $\zeta$-reflection in $r$.  If $r$ has norm~3 and lies in an
$\E$-lattice in which $\theta$ divides all inner products, such as
$\Lambda(V)$, then \eqref{eq-defn-of-complex-reflection} shows that $\w$-reflection in $r$ also
preserves the lattice.

\begin{lemma}
\label{lem-meridians-act-by-complex-reflections}
The image of a meridian under the monodromy representation
$\rho:\pi_1(\forms_0,F)\to\Gamma(V)=\aut \Lambda(V)$ is the
$\w$-reflection in an element of $\Lambda(V)$ of norm~$3$.
\end{lemma}

\begin{proof}
The argument is much the same as for lemma~5.4 of \cite{ACT}.  Let $D$ be
a small disk in $\forms$, meeting $\D$ only at its center, and
transversally there.  We write $F_0$ for the form at the center of $D$.
Suppose without loss of generality that the basepoint $F$ of
$\forms_0$ is on $\partial D$, and let $\c$ be the element of
$\pi_1(\forms_0,F)$ that traverses $\partial D$ once positively.  The
essential ingredients of the proof are the following.  First, $T_0$
has an $A_1$ singularity, so $V_0$ has an $A_2$ singularity; this
means that in
suitable local analytic coordinates it is given by
$x_1^2+\dots+x_4^2+x_5^3=0$.  Second, the vanishing cohomology for
this singularity, i.e., the Poincar\'e dual of the kernel of 
$$
H_4(V;\Z)\to H_4(\V|_D;\Z)\isomorphism H_4(V_0;\Z)\;,
$$
is a (positive-definite) copy of the $A_2$ root lattice.  (An $A_2$
surface singularity has vanishing cohomology a negative-definite copy
of this lattice, and the sign changes when the dimension increases
by~2.)  Third, following Sebastiani-Thom \cite{sebastiani-thom}, this
lattice may be described as 
$$
V(2)\tensor V(2)\tensor V(2)\tensor V(2)\tensor V(3)\;,
$$
where $V(k)$ is the $\Z$-module spanned by the differences of the
$k$th roots of unity,  $\c$ acts by
$$
-1\tensor-1\tensor-1\tensor-1\tensor\w
$$
and $\sigma$ acts by
$$
1\tensor1\tensor1\tensor1\tensor\w\;.
$$
Here, $\pm1$ (resp.\ $\w$) indicates the action on $V(2)$ (resp.\ $V(3)$)
given by sending each square root (resp. cube root) of unity to itself times
$\pm1$ (resp.\ $\w$).  This shows that the vanishing cohomology is a
1-dimensional $\E$-lattice; we write $r$ for a generator.  It also shows
that $\c$ acts on $\langle r\rangle$ in the same way that $\sigma^*$
does.  Since $\w$'s action on $\Lambda(V)$ is defined to be $\sigma^*$,
$\c$ acts on $\langle r\rangle$ by $\w$.  Fourth, $\c$ acts trivially
on the orthogonal complement of the vanishing cohomology in
$H^4(V;\Z)$; this implies that $\c$ is the $\w$-reflection in $r$.
Finally, since the roots of the $A_2$ lattice have norm~2, we see by
\eqref{eq-defn-of-hermitian-form} that $\ip{r}{r}=3$.
\end{proof}

The following two lemmas play only a small role in this section, at
one point in the proof of theorem~\ref{thm-inner-product-matrix}, to
which the reader could skip right away.  However, they will be very
important in section~\ref{sec-extension}, where we extend the domain
of the period map.  Their content is that the discriminant has nice
local models, which make many homomorphisms from braid groups into
$\pi_1(\forms_0)$ visible.  
We also show that distinct braid generators
have distinct monodromy actions.  

We recall that the fundamental
group of the discriminant complement of an $A_n$ singularity is the
braid group $B_{n+1}$, also known as the Artin group $\A(A_n)$ of type
$A_n$.  More generally, the Artin group of an $A_n$, $D_n$ or $E_n$
Dynkin diagram has one generator for each node, with two of the
generators braiding ($aba=bab$) or commuting, corresponding to whether
the corresponding nodes are joined or not.  It is the fundamental
group of the discriminant complement of that corresponding
singularity  \cite{brieskorn-pi1s}. 
Only $\A(A_n)$ and
$\A(D_4)$ will be relevant to this paper.

\begin{lemma}
\label{lem-discriminant-away-from-E}
Suppose $F\in\forms$ defines a cubic threefold with singularities
 $s_1,\dots,$ $s_m$, each having one of the types $A_n$ or $D_4$, and
 no other singularities.  Let $K_{i=1,\dots,m}$ be the base of a miniversal
 deformation of a singularity having the type of $s_i$, with
 discriminant locus $\D_i\sset K_i$.  Then there is a neighborhood $U$
 of $F$ in $\forms$ diffeomorphic to $K_1\times\dots\times K_m\times
 B^N$, where $N=35-\sum\dimension K_i$, such that $U-\D$ corresponds
 to
$$
(K_1-\D_1)\times\dots\times(K_m-\D_m)\times B^N\;.
$$
In particular, $\pi_1(U-\D)$ is the direct product of $m$ Artin
groups, the $i$th factor having the type of the singularity $s_i$. 
\end{lemma}

\begin{proof}
This is essentially the assertion that $\forms$ contains a
simultaneous versal deformation of all the singularities of $T$.  By
theorem~1.1 of \cite{du-plessis-wall}, it suffices to show that the sum of
the Tjurina numbers of $s_1,\dots,s_m$ is less than~16.  Because the
singularities of $T$ are quasihomogeneous, their Tjurina numbers
coincide with their Milnor numbers.

We will write $\mu_i$ for the Milnor number of $T$'s singularity at
$s_i$, and $Z_i\sset H_0^4(V)$ for the vanishing cohomology of the
corresponding singularity of $V$.  If $T$ has a $D_4$ singularity at
$s_i$, then $\mu_i=4$, and $V$ has an $\tilde E_6$ singularity there, with
$\dimension Z_i=8$ and $\dimension(Z_i\cap Z_i^\perp)=2$.  When $T$
has an $A_n$ singularity at $s_i$, we have $\mu_i=n$, and $V$ has a
singularity locally modeled on $x_0^2+x_1^2+x_2^2+y^3+z^n=0$.  By
\cite[p.~77]{arnold-et-al}, $Z_i$ has a basis
$a_1,\dots,a_n,b_1,\dots,b_n$ with $a_i^2=b_i^2=2$, $a_i\cdot
a_{i\pm1}=-1$, $b_i\cdot b_{i\pm1}=-1$, $a_i\cdot b_i=-1$, $a_i\cdot
b_{i-1}=1$, and all other inner products zero.  For $n>11$ this
quadratic form has a negative-definite subspace of dimension${}\geq4$,
so it cannot lie in $H_0^4(V)$, which has signature $(20,2)$.
Therefore cubic threefolds cannot have $A_{n>11}$ singularities.  For
$n<12$, $Z_i$ is nondegenerate except for $n=5$ and~$11$, when
$\dimension (Z_i\cap Z_i^\perp)=2$.  We made these calculations using
PARI/GP \cite{pari}.

In every case, we have
$\mu_i\leq\frac{2}{3}\dimension\bigl(Z_i/(Z_i\cap Z_i^\perp)\bigr)$.
Since $i\neq j$ implies $Z_i\perp Z_j$, we have 
$$
\sum\dimension \bigl(Z_i/(Z_i\cap Z_i^\perp)\bigr)
\leq \dimension H_0^4(V)=22\;.
$$
Putting these inequalities together yields
$\sum\mu_i\leq\frac{2}{3}\cdot22<16$, so that \cite{du-plessis-wall} applies.
This gives the claimed description of $\D$ near $F$, and the
description of $\pi_1(U-\D)$ follows immediately. 
\end{proof}

\begin{lemma}
\label{lem-local-monodromy-noncyclic}
In the situation of the previous lemma, suppose $g$ and $g'$ are two
of the standard generators for $\pi_1(U-\D,F')$, where $F'$ is a
basepoint in $U-\D$.
Then $\rho(g),\rho(g')\in\Gamma(V')$ are the $\w$-reflections in
linearly independent roots
$r,r'\in \Lambda(V')$.
\end{lemma}

\begin{proof}
If $T$ has two $A_1$ singularities, then near $F$, $\D$ has two
components, and $g$ and $g'$ are meridians around them. As we saw in
the proof of lemma~\ref{lem-meridians-act-by-complex-reflections}, $V$ has two $A_2$ singularities, and the
vanishing cohomology of each of them is a $\sigma$-invariant
sublattice of $H_0^4(V')$.  In fact, they are the $\E$-spans
of $r$ and $r'$.  With respect to the cup product, vanishing cocycles
for distinct singularities are orthogonal, and it follows from
\eqref{eq-defn-of-hermitian-form} that they are also orthogonal under $\ip{}{}$.  Therefore
$r\perp r'$.  Since $r$ and $r'$ have nonzero norm, they must be
linearly independent.

Now allow $T$ to have $A_n$ and/or $D_4$ singularities, and suppose
$g$ and $g'$ commute.  Then there exists $F_0\in U$ and a neighborhood
$U_0$ of $F_0$ in $U$ with the following properties.  $T_0$ has just
two singularities, both nodes, $U_0$ contains $F'$,
$\pi_1(U_0-\D,F')\isomorphism\Z^2$, and $g,g'\in\pi_1(U-\D,F')$ are
represented by loops in $U_0$ that are meridians around the two components
of $\D$ at $F_0$.  These properties follow from Brieskorn's
description \cite{brieskorn-versals}
of the versal deformations of simple singularities in
terms of the corresponding Coxeter groups.  By the previous paragraph,
$r$ and $r'$ are linearly independent.

The same ideas show that if $g$ and $g'$ braid but do not commute,
then there exists $F_0\in U$ and a neighborhood $U_0$ of $F_0$ in $U$
with the following properties.  $T_0$ has just one singularity, of
type $A_2$, $F'\in U_0$, $\pi_1(U_0-\D,F')\isomorphism B_3$, and
$g,g'\in\pi_1(U-\D,F')$ are represented by the standard generators of
$B_3$.  It is easy to see that the set of cubic threefolds having an
$A_2$ singularity and no other singularities is irreducible.
Therefore: if there is a single counterexample to the lemma, then
$\rho(g)=\rho(g')$ for every pair of generators of $\pi_1(U-\D,F')$,
for every $F$ as in the lemma.  So it suffices to treat the case that
$T$ has a single singularity, of type $A_3$, and $g$ and $g'$ are the
first two standard generators of $\pi_1(U-\D,F')\isomorphism B_4$.
Adjoining the relation $g=g'$ reduces $B_4$ to $\Z$.  If
$\rho(g)=\rho(g')$, then $\rho|_H$ factors through $\Z$, which implies
that all three generators of $B_4$ have the same $\rho$-image.  This
is impossible by the previous paragraph, because the first and last
generators commute.
\end{proof}

\begin{theorem}
\label{thm-inner-product-matrix}
For $F\in\forms_0$, $\Lambda(V)$ is isometric to the $\E$-lattice with
inner product matrix
\begin{equation}
\label{eq-inner-product-matrix}
\Lambda:= \begin{pmatrix}3\end{pmatrix}
\oplus
\begin{pmatrix}
3&\theta&0&0\\
\thetabar&3&\theta&0\\
0&\thetabar&3&\theta\\
0&0&\thetabar&3
\end{pmatrix}
\oplus
\begin{pmatrix}
3&\theta&0&0\\
\thetabar&3&\theta&0\\
0&\thetabar&3&\theta\\
0&0&\thetabar&3
\end{pmatrix}
\oplus
\begin{pmatrix}
0&\theta\\
\thetabar&0
\end{pmatrix}
\;.
\end{equation}
\end{theorem}

\begin{remarks}
Regarding \eqref{eq-inner-product-matrix} as an $11\times11$ matrix $(\lambda_{ij})$, this
means that $\Lambda=\E^{11}$, with 
$$
\ip{\strut(x_1,\dots,x_{11})}{(y_1,\dots,y_{11})}=\sum_{i,j}\lambda_{ij}x_i\bar
y_j\;.
$$
The four-dimensional lattice appearing twice among the
summands is called  $E_8^\E$, because its
underlying $\Z$-lattice is a scaled copy of the $E_8$ root lattice.
\end{remarks}

\begin{proof}
By section~2 of \cite{allcock-threefolds}, the cubic threefold $T_0$ defined by 
$$
F_0=x_2^3+x_0x_3^2+x_1^2x_4-x_0x_2x_4-2x_1x_2x_3+x_4^3
$$ has an $A_{11}$ singularity at $[1{:}0{:}\dots{:}0]\in\cp^4$ and no other
singularities.  By lemma~\ref{lem-discriminant-away-from-E}, we may
choose a neighborhood $U$ of $F_0$ and $F\in U-\D$ such that
$\pi_1(U-\D,F)\isomorphism B_{12}$.  By
lemma~\ref{lem-meridians-act-by-complex-reflections}, the standard
generators of $B_{12}$ act on $\Lambda(V)$ by the $\w$-reflections
$R_1,\dots,R_{11}$ in pairwise linearly independent roots
$r_1,\dots,r_{11}\in \Lambda(V)$.  The commutation relations imply
$r_i\perp r_j=0$ if $j\neq i\pm1$.  By the argument of
\cite[\S5]{allcock-inventiones}, the relation
$R_iR_{i+1}R_i=R_{i+1}R_iR_{i+1}$ implies that $\bigl|\langle
r_i|r_{i+1}\rangle\bigr|=\sqrt3$, so after multiplying some of the
$r_i$ by scalars, we may take $\ip{r_i}{r_{i+1}}=\theta$.

The rank of the inner product matrix of the $r_i$ is $10$.  Therefore,
if they were linearly independent then they would span $\Lambda$ up to
finite index, and the Hermitian form on $\Lambda$ would be
degenerate.  It is not, so the span must be only $10$-dimensional.
By the argument of \cite[\S5]{allcock-inventiones}, the $r_i$ span a
sublattice of $\Lambda(V)$ isometric to the direct sum (call it
$\Lambda_{10}$) of the last three summands of
\eqref{eq-inner-product-matrix}.  In \cite{allcock-inventiones} we
used a form of signature $(1,9)$ rather than $(9,1)$; this difference
is unimportant.  One can check directly that $\theta
\Lambda_{10}^*=\Lambda_{10}$; the 
underlying reason is that the real forms of $E_8^\E$ and
$\bigl(\begin{smallmatrix}0&\theta\\\thetabar&0\end{smallmatrix}\bigr)$
are scaled copies of even unimodular $\Z$-lattices.  Since
$\Lambda(V)\sset\theta \Lambda(V)^*$, $\Lambda_{10}$ is a direct summand of
$\Lambda(V)$, so $\Lambda(V)\isomorphism (n)\oplus \Lambda_{10}$ for some
$n\in\Z$.  We have $n>0$ because $\Lambda(V)$ has signature $(10,1)$.

For an $\E$-lattice $M$ we define $M^\Z$ to be the $\Z$-module
underlying $M$, equipped with the $\Z$-bilinear pairing
$\a\cdot\b=\frac{2}{3}\re\ip{\a}{\b}$.  Computation shows that
$(n)^\Z$ has inner product matrix
$\bigl(\begin{smallmatrix}2n/3&-n/3\\-n/3&2n/3\end{smallmatrix}\bigr)$,
$(E_8^\E)^\Z$ is the even unimodular $\Z$-lattice $E_8$, and
$\bigl(\begin{smallmatrix}0&\theta\\\thetabar&0\end{smallmatrix}\bigr)\strut^\Z$ is the
even unimodular $\Z$-lattice
$II_{2,2}=\bigl(\begin{smallmatrix}0&1\\1&0\end{smallmatrix}\bigr)\oplus\bigl(\begin{smallmatrix}0&1\\1&0\end{smallmatrix}\bigr)$.
Since $(\Lambda(V))^\Z=H^4_0(V;\Z)$ has determinant~$\pm3$, we must
have $n=3$.
\end{proof}

We define a framing of a form $F\in\forms_0$ to be an equivalence
class $[\phi]$ of isometries $\phi:\Lambda(V)\to \Lambda$, two
isometries being
equivalent if they differ by multiplication by a scalar.  Sometimes we
write $\phi$ rather than $[\phi]$ and leave it to the reader to check
that the construction at hand depends only on $[\phi]$.  We define
$\framed_0$ to be the set of all framings of all smooth cubic forms.
Since the stalk of $\Lambda(\V_0)$ at $F\in\forms_0$ is canonically
isomorphic to $\Lambda(V)$, the set $\framed_0$ is in natural
bijection with the subsheaf of $P\!\Hom\bigl(\Lambda(\V_0),\forms_0\times
\Lambda\bigr)$ consisting of projective equivalence classes of
homomorphisms which are isometries on each stalk.  This bijection
gives $\framed_0$ the structure of a complex manifold.  We refer to an
element $(F,[\phi])$ of $\framed_0$ as a framed smooth cubic form.

We write $\Gamma$ for $\aut \Lambda$ and $P\Gamma$ for $P\!\aut
\Lambda$.
On $\framed_0$ are defined commuting actions of $\PGamma$ and
$G=\GL_5\C/D$, where $D$ is the group $\{I,\w I,\wbar I\}$.  An element
$\gamma$ of $\PGamma$ acts on the left by
$$
\gamma.(F,[\phi])=(F,[\gamma\circ\phi])\;.
$$
This action realizes $\PGamma$ as the group of deck
transformations of the covering space $\framed_0\to\forms_0$.  An
element $g$ of $\GL_5\C$ acts on the right by
\begin{equation}
\label{eq-GL5C-action-on-framed}
\bigl(F,[\phi]\bigr).g=\bigl(F\circ g,[\phi\circ g^{-1}{}^*]\bigr)\;.
\end{equation}
Here, $\GL_5\C$ acts on $\C^5$ on the left, hence acts on $\forms$ on
the right by $(F.g)(x)=F(g.x)$, i.e., $F.g=F\circ g$.  We extend
$\GL_5\C$'s action on $\C^5$ to $\C^6=\C^5\oplus\C$ by the trivial
action on the $\C$ summand.  This induces a right action on $\V$ by
$(F,x).g=(F.g,g^{-1}x)$.  That is, $g$ carries the zero-locus of
$(F+x_5^3).g$ to the zero-locus of $F+x_5^3$.  The $g^{-1}{}^*$
appearing in \eqref{eq-GL5C-action-on-framed} is the inverse of the
induced map on cohomology, which respects the $\E$-module structure
since $g$ commutes with $\sigma$.  The subgroup $D\sset \GL_5\C$ acts
trivially on $\framed_0$ because it acts trivially on $\forms_0$ and
by scalars on every $\Lambda(V)$.

We now introduce the moduli spaces $\moduli_0$ and $\moduli_0^f$ of
smooth and framed smooth cubic threefolds.  Since
$\forms_0\sset\forms_s$, $G$ acts properly on $\forms_0$, with the
quotient $\moduli_0=\forms_0/G$ a complex analytic orbifold and a
quasiprojective variety.  The properness on $\forms_0$ implies
properness on $\framed_0$, so $\moduli_0^f=\framed_0/G$ is a complex
analytic orbifold and an analytic space.  We will see in lemma~\ref{lem-G-acts-freely-and-per-map-has-full-rank}
that $\moduli_0^f$ is a complex manifold.  Since the $G$-stabilizer of
a point $(F,[\phi])$ of $\framed_0$ is a subgroup of the
$G$-stabilizer of $F\in\forms_0$, the covering map
$\framed_0\to\PGamma\backslash\framed_0=\forms_0$ descends to an
orbifold covering map
$\moduli_0^f=\framed_0/G\to\forms_0/G=\moduli_0$.  Counting dimensions
shows that $\moduli_0$ and $\moduli_0^f$ are 10-dimensional.

We write $\ch^{10}$ for $\ch(\Lambda\tensor_\E\C)$.  Recall that
$H^{3,1}(V;\C)$ is a negative line in the Hermitian vector space
$H^4_\w(V;\C)$ and that $Z$ is an isometry $\Lambda(V)\tensor_\E\C\to
H^4_\w(V;\C)$.  We reformulate the period map
\eqref{eq-per-map-as-section} as the holomorphic map
$\per:\framed_0\to\ch^{10}$ given by
\begin{equation}
\label{eq-per-map-formulated-with-framings}
g(F,[\phi])=\phi\bigl(Z^{-1}\bigl(H^{3,1}(V;\C)\bigr)\bigr)\;.
\end{equation}
On the right we have written just $\phi$ for $\phi$'s $\C$-linear
extension $\Lambda(V)\tensor_\E\C\to \Lambda\tensor_\E\C$.  Since
$\ch^{10}$ is a 10-ball and bounded holomorphic functions on $\GL_5\C$
are constant, $\per$ is constant along $\GL_5\C$-orbits, so it
descends to a holomorphic map $g:\moduli_0^f\to\ch^{10}$, also called
the period map.  This map is equivariant with respect to the action of
$\PGamma$, so it in turn descends to a map
\begin{equation}
\label{eq-per-map-as-map-to-ball-quotient}
\per:\moduli_0
=\PGamma\backslash\moduli_0^f
=\PGamma\backslash\framed_0/\GL_5\C
\to\PGamma\backslash\ch^{10}\;, 
\end{equation}
again called the period map.

\begin{lemma}
\label{lem-G-acts-freely-and-per-map-has-full-rank}
$G$ acts freely on $\framed_0$, so that $\moduli_0^f$ is a complex
manifold, not just an orbifold.  The period map
$g:\moduli_0^f\to\ch^{10}$ has rank~$10$ at every point of
$\moduli_0^f$. 
\end{lemma}

\begin{proof}
We prove the second assertion first.  Let $F\in\forms_0$, let
$F'\in\forms$ be different from $F$, and let $\e>0$ be small enough
that the disk $D=\set{F+tG}{\hbox{$t\in\C$ and $|t|\leq\e$}}$ lies in
$\forms_0$.  Writing $F_t$ for $F+tG$, we know from the discussion
surrounding \eqref{eq-Omega(A)} that
$H^{3,1}(V_t)$ is spanned by the residue of $\Omega/(F_t+x_5^3)^2$.
Since $\V$ trivializes over $D$, we may unambiguously translate this
class into $H^4(V;\C)$;  this gives a map $h:D\to H^4(V;\C)$.  For
sufficiently small $t$, $h(F_t)$ is the element of
$\Hom\bigl(H_4(V;\Z),\C\bigr)$ given by 
$$
\hbox{(an integral 4-cycle $C$)}\mapsto
\int_{\partial N}\frac{\Omega}{(F_t+x_5^3)^2}
$$
where $N$ the the part of the boundary of a tubular neighborhood of
$V$ in $\cp^5$ that lies over a submanifold of $V$ representing $C$.
Therefore we may differentiate with respect to $t$ under
the integral sign, so the derivative of $h$ at the center of $D$ is
the element of $\Hom\bigl(H_4(V;\Z),\C\bigr)$ given by 
$$
C\mapsto\int_{\partial N}\frac{\Omega}{(F_t+x_5^3)^3}
\cdot(-2)\frac{\partial}{\partial t}(F_t+x_5^3)\Biggr|_{t=0}
=
-2\int_{\partial N}\frac{\Omega F'(x_0,\dots,x_4)}{(F+x_5^3)^3}\;.
$$ This lies in $H^4_\w(V;\C)$, and it lies in $H^{3,1}$ if and only
if $F'$ lies in the Jacobian ideal of $F+x_5^3$, i.e., if and only if
$F'$ lies in the Jacobian ideal of $F$, i.e., if and only if the pencil
$\langle F,F'\rangle$ in $\forms$ is tangent to the $G$-orbit of $F$.

Upon choosing a framing $\phi$ for $F$ and lifting $D$ to a disk
$\tilde D=\{(F_t,[\phi_t])\}$ in $\framed_0$ passing through
$(F,[\phi])$, it follows that the derivative of
$g:\framed_0\to\ch^{10}$ along $\tilde D$ at $(F,[\phi])$ is zero if
and only if $\tilde D$ is tangent to the $G$-orbit of $(F,[\phi])$.
Since the orbit has codimension~10, $g$ has rank~$10$.

To prove the first assertion, recall that an orbifold chart about the
image of $(F,[\phi])$ in $\moduli_0^f$ is $U\to
U/H\sset\moduli_0^f$, where $H$ is the $G$-stabilizer of $(F,[\phi])$
and $U$ is a small $H$-invariant transversal to the the $G$-orbit of
$(F,[\phi])$.  We have just seen that the composition $U\to
U/H\sset\moduli_0^f\to\ch^{10}$ has rank~10 and is hence a
diffeomorphism onto its image.  It follows that $H=\{1\}$.
\end{proof}

\begin{theorem}
\label{thm-main-theorem-smooth-case}
The period map $\per:\moduli_0\to\PGamma\backslash\ch^{10}$ is
an isomorphism onto its image.
\end{theorem}

\begin{proof}
We begin by proving that if $F$ and $F'$ are generic elements of
$\forms_0$ with the same image under $\per$ then they are
$G$-equivalent, i.e., $T$ and $T'$ are projectively equivalent.  By
hypothesis there exists an isometry $b:\Lambda(T)\to \Lambda(T')$
which carries $Z^{-1}(H^{3,1}(V))\in \Lambda(V)\tensor_\E\C$ to
$Z^{-1}(H^{3,1}(V'))\in \Lambda(V')\tensor_\E\C$.  Passing to the
underlying integer lattices, $b$ is an isometry
$\Hprim^4(V;\Z)\to\Hprim^4(V';\Z)$ carrying $H^{3,1}(V;\C)$ to
$H^{3,1}(V';\C)$.  By complex conjugation it also identifies
$H^{1,3}(V;\C)$ with $H^{1,3}(V';\C)$, and by considering the
orthogonal complement of $H^{3,1}\oplus H^{1,3}$ we see that it
identifies $\Hprim^{2,2}(V;\C)$ with $\Hprim^{2,2}(V';\C)$.  That is,
it induces an isomorphism of Hodge structures.

Next: one of $\pm b$ extends to an isometry $H^4(V;\Z)\to
H^4(V';\Z)$ carrying $\eta(V)$ to $\eta(V')$.  This follows from some
lattice-theoretic considerations: if $L$ is a nondegenerate 
primitive sublattice of a
unimodular lattice $M$, that is, $L=(L\tensor\Q)\cap M$,
then the projections of $M$ into $L\tensor\Q$ and $L^\perp\tensor\Q$
define an isomorphism of $L^*/L$ with $(L^\perp)^*/L^\perp$.  Here
the asterisk denotes the dual lattice.  It is easy to check that an
isometry of $L$ and an isometry of $L^\perp$ together give an isometry
of $M$ if and only if their actions on $L^*/L$ and
$(L^\perp)^*/L^\perp$ coincide under this identification.  Since
$\bigl\langle\eta(V)\bigr\rangle^*/\bigl\langle\eta(V)\bigr\rangle\isomorphism
\Z/3$, it 
follows that exactly one of the isometries
$$
\bigl\langle\eta(V)\bigr\rangle\oplus\Hprim^4(V;\Z)\to 
\bigl\langle\eta(V')\bigr\rangle\oplus\Hprim^4(V';\Z)\;,
$$
given on the first summand by $\eta(V)\mapsto\eta(V')$, and on the
second by $\pm b$, extends to an isometry $H^4(V;\Z)\to H^4(V';\Z)$.

 From Claire Voisin's theorem \cite{voisin} we deduce that there is a
projective transformation $\b$ carrying $V$ to $V'$.  
One can check that the variety $S$ of smooth cubic fourfolds admitting
a triflection is irreducible, so that one can speak of a generic such
fourfold.  Furthermore, a generic such fourfold admits
only one triflection (and its inverse).  Since $V$ and $V'$ admit the
triflections $\sigma^{\pm1}$ and are generic points of $S$, $\b$ carries the fixed-point set $T$ of
$\sigma$ in $V$ to the fixed-point set $T'$ of $\sigma$ in $V'$.
That is, $T$ and $T'$ are projectively equivalent.

We have proven that the period map from $\moduli_0$ to
$\PGamma\backslash\ch^{10}$ is generically injective, and 
the previous lemma shows that it is a local isomorphism.  It follows that it
is an isomorphism onto its image.
\end{proof}

\section{The discriminant near a chordal cubic}
\label{sec-discr-near-chordal-cubic}

In the next section we will enlarge the domain of the period map
$\forms_0\to\PGamma\backslash\ch^{10}$, in order to obtain a map from
a compactification of $\moduli_0$ to the Baily-Borel compactification
$\overline{\PGamma\backslash\ch^{10}}$.  In order to do this we will
need to understand the local structure of the discriminant
$\D\sset\forms_0$, at least near the threefolds to which we will
extend $g$.  In \cite{allcock-threefolds} (see also \cite{yokoyama}),
the GIT-stability of cubic threefolds is completely worked out.  There
is one distinguished type of threefold, which we call a chordal cubic,
which is the secant variety of the rational normal quartic curve.
Except for the chordal cubics and those cubics that are GIT-equivalent
to them, a cubic threefold is semistable if and only if it has
singularities only of types $A_1,\dots,A_5$ and $D_4$.  At such a
threefold the local structure of $\D$ is given by
lemma~\ref{lem-discriminant-away-from-E}.

The rest of this section addresses the nature of $\D$ near the chordal
cubic locus.  It turns out (see the remark following
theorem~\ref{thm-chordal-maps-onto-divisor}) that the period map
$P\forms_0\to\PGamma\backslash\ch^{10}$ does not extend to a regular
map $P\forms_{ss}\to\overline{\PGamma\backslash\ch^{10}}$.  The
problem is that it does not extend to the chordal cubic locus.
Therefore it is natural to try to enlarge the domain of the period map
not to $P\forms_{ss}$ but rather to $(\widehat{P\forms})_{ss}$, where
$\widehat{P\forms}$ is the blowup of $P\forms$ along the closure of
the chordal cubic locus.  The details concerning the GIT analysis and
the extension of the period map appear in section~\ref{sec-extension};
at this point we are only motivating the study of the local structure
of the proper transform $\Dhat$ of $\D$ along the exceptional divisor
$E$.  Recall that we defined $\D$ as a subset of $\forms$, but will
also write $\D$ for its image in $P\forms$.

If $T\in P\forms$ is a chordal cubic then we write $E_T$ for
$\pi^{-1}(T)\sset E$, where $\pi$ is the natural projection $\widehat
{P\forms}\to P\forms$.  (There are a number of projection maps in this
paper, such as $\pi_\T$ and $\pi_\V$ in
section~\ref{sec-smoothmoduli}, and some others introduced later.  To
keep them straight, we will use a subscript to indicate the domain for
all of them except this one.)  $E_T$ may be
described as the set of unordered 12-tuples in the rational normal
curve $R_T$
which is the singular locus of $T$.  To see this, one
counts dimensions to find that the chordal cubic locus has
codimension~13 in $P\forms$, so $E_T$ is a copy of $P^{12}$.  To
identify $E_T$ with the set of unordered 12-tuples in $R_T$, consider
a pencil of cubic threefolds degenerating to $T$.  The 12-tuple may be
obtained as the intersection of $R_T$ with a generic member of the
pencil; since $R_T$ has degree~4, this intersection consists of 12
points, counted with multiplicity.  (If every member of the pencil
vanishes identically on $R_T$ then the pencil is not transverse to the
chordal cubic locus.)  We will indicate an element of $E_T$ by a pair
$(T,\tau)$, where $\tau$ is an unordered 12-tuple in $R_T$.

We will describe $\Dhat\sset\widehat{P\forms}$ by proving the
following two theorems, which are similar to but weaker than
lemma~\ref{lem-discriminant-away-from-E}.  The first is weaker because
it asserts a homeomorphism with a standard model of the discriminant,
rather than a diffeomorphism.  The second gives a complex-analytic
isomorphism, but refers to a finite cover of (an open set in)
$\widehat{P\forms}$, branched over $E$.  But we don't know any reason
that the homeomorphism is
theorem~\ref{thm-discr-along-E-topological-model} couldn't be promoted
to a diffeomorphism.

\begin{theorem}
\label{thm-discr-along-E-topological-model}
Suppose $T$ is a chordal cubic and $\tau$ is a $12$-tuple in
$R_T$, with $m$ singularities, of types $A_{n_1},\dots,A_{n_m}$, where
an $A_n$ singularity means a point of multiplicity $n+1$.  Let
$K_{i=1,\dots,m}$ be the base of a miniversal deformation of an $A_{n_i}$ singularity,
with discriminant locus $\D_i\sset K_i$.  Then there is a neighborhood
$U$ of $(T,\tau)$ in $\widehat{P\forms}$ homeomorphic to $B^1\times
K_1\times\dots\times K_m\times B^N$, where $N=33-\sum\dimension K_i$, such that $E$ corresponds to
$$
\{0\}\times K_1\times\dots\times K_m\times B^N
$$
and $U-\widehat\D$  to
$$
B^1\times(K_1-\D_1)\times\dots\times(K_m-\D_m)\times B^N\;.
$$
In particular, 
$$
\pi_1\bigl(U-(\Dhat\cup E)\bigr)
\isomorphism
\Z\times B_{n_1+1}\times\dots\times B_{n_m+1}\;,
$$
where the $\Z$ factor is generated by a meridian of $E$ and the
standard generators of the braid group factors are meridians of
$\Dhat$. 
\end{theorem}

\begin{theorem}
\label{thm-discr-along-E-analytic-model}
Suppose $(T,\tau)$, $m$, $n_1,\dots,n_m$, $K_1,\dots,K_m$,
$\D_1,\dots,$ $\D_m$ and $N$ are as in
theorem~\ref{thm-discr-along-E-topological-model}.  Then there exists a
neighborhood $U$ of $(T,\tau)$ in  $\widehat{P\forms}$ diffeomorphic to
$B^1\times B^{33}$, with $U\cap E$ corresponding to $\{0\}\times
B^{33}$, such that the following holds.  We write 
$\pi_{\tilde U}:\tilde U\to U$ 
for the $6$-fold cover of $U$ branched over $U\cap E$, and
$(T,\tau)\sptilde$ for the point $\pi_{\tilde U}^{-1}(T,\tau)$.  There is a
neighborhood $V$ of $(T,\tau)\sptilde$ in $\tilde U$ diffeomorphic to
$B^1\times K_1\times \dots\times K_m\times B^N$, such that 
\begin{equation}
\label{eq-local-model-description-of-E}
V\cap\pi_{\tilde U}^{-1}(E)=\{0\}\times K_1\times\dots\times K_m\times B^N
\end{equation}
and
\begin{equation}
\label{eq-local-model-of-discriminant-complement-along-E}
V-\pi_{\tilde U}^{-1}(\widehat\D)=
B^1\times(K_1-\D_1)\times\dots\times(K_m-\D_m)\times B^N\;.
\end{equation}
\end{theorem}


The rest of this section is devoted to proving
theorems~\ref{thm-discr-along-E-topological-model}
and~\ref{thm-discr-along-E-analytic-model}.  It is rather
technical, especially lemma~\ref{lem-rigidity-of-discriminant} and beyond; although these theorems
are analogues of lemma~\ref{lem-discriminant-away-from-E}, the
proofs are much more complicated.  

\begin{lemma}
\label{lem-Deltahat-intersect-E}
Suppose $T$ is a chordal cubic and $(T,\tau)\in E_T$.
\renewcommand\theenumi{\roman{enumi}}
\begin{enumerate}
\item
\label{item-Deltahat-intersect-E}
$(T,\tau)$ lies in $\Dhat$ if and only if $\tau$ has a multiple point.
\item
\label{item-limits-of-Ans}
If $\tau$ has a point of multiplicity $n+1$, then $(T,\tau)$ is a
limit of points of $\widehat{P\forms}$ representing cubic threefolds
with $A_n$ singularities.
\end{enumerate}
\end{lemma}

Nowhere else in the paper do we refer to any result or notation
introduced from here to the end of this section.  

When we refer to the ``standard chordal cubic'', we
mean the one defined by 
\begin{equation}
\label{eq-standard-chordal-cubic}
F(x_0,\dots,x_4)=\det
\begin{pmatrix}
x_0&x_1&x_2\\
x_1&x_2&x_3\\
x_2&x_3&x_4
\end{pmatrix}
=0\;,
\end{equation}
which is the secant variety of the rational normal curve parameterized
by 
\begin{equation}
\label{eq-parameterization-of-RNC}
s\mapsto[1,s,s^2,s^3,s^4]\qquad(s\in P^1).
\end{equation}
We write $P$ for the
point $[1,0,0,0,0]$.

\begin{proof}[Proof of lemma~\ref{lem-Deltahat-intersect-E}]
We prove \eqref{item-limits-of-Ans} first.
We take $T$ to be the standard chordal cubic, and place the multiple
point of $\tau$ at $P\in R_T$. 
We suppose without loss of generality that
$[0,0,0,0,1]\in R_T$ is not one of the points of $\tau$.  Observe that 
\begin{equation}
\label{eq-family-of-Gs-for-A11}
\begin{split}
G_{u_1,\dots,u_{12}}=F+x_4^3
&{}
+u_1x_4^2x_3
+u_2x_4x_3^2
+u_3x_3^3
+u_4x_3^2x_2
\\&{}
+u_5x_3x_2^2
+u_6x_2^3
+u_7x_2^2x_1
+u_8x_2x_1^2
\\&{}
+u_9x_1^3
+u_{10}x_1^2x_0
+u_{11}x_1x_0^2
+u_{12}x_0^3
\end{split}
\end{equation}
restricts to $R_T$ as the polynomial
\begin{equation}
\label{eq-restriction-of-Gs-to-rational-normal-curve}
s^{12}+u_1s^{11}+u_2s^{10}+\dots+u_{11}s+u_{12}\;,
\end{equation}
where $R_T$ is parameterized as in \eqref{eq-parameterization-of-RNC}.  Since $\tau$ has no
point at $s=\infty$, there is a choice of $u_1,\dots,u_{12}$ such that
$\tau$ is the limiting direction of the pencil $\langle
F,G_{u_1,\dots,u_{12}}\rangle$.  Since $\tau$ has a point of
multiplicity $n+1$ at $P$, we have $u_{12}=\dots=u_{12-n}=0$ and
$u_{11-(n+1)}\neq0$.  Then singularity analysis as in
\cite[sec.~2]{allcock-threefolds} shows that a generic member of the
pencil has an $A_n$ singularity at $P$.  This proves
\eqref{item-limits-of-Ans}.

Now we prove \eqref{item-Deltahat-intersect-E}.  If $\tau$ has a
multiple point then \eqref{item-limits-of-Ans} shows that
$(T,\tau)\in\widehat\D$.  If $\tau$ has no multiple points, and
$\langle F,F'\rangle$ is a pencil in $P\forms$ with limiting direction
$\tau$, then $dF'$ vanishes at no point of $R_T$, because otherwise
$\tau$ would have a multiple point.  Since $dF'$ is nonvanishing along
the locus $dF=0$, there exists $\e>0$ and a neighborhood $W$ of $F'$
in $\forms$ such that $dF+\eta\,dF''$ is nowhere-vanishing, for all
$0<|\eta|<\e$ and all $F''\in W$, so $F+\eta F''$ defines a smooth
threefold.  Therefore $(T,\tau)$ has a neighborhood in
$\widehat{P\forms}$ disjoint from $\Dhat-E$, so $(T,\tau)\notin\Dhat$.
\end{proof}

Because of the action of $PG$, proving theorems~\ref{thm-discr-along-E-topological-model} and~\ref{thm-discr-along-E-analytic-model}
reduces to a similar but lower-dimensional problem.  Let $T$ be the
standard chordal cubic.  In our arguments, the 12-tuple consisting of
12 points all concentrated at $P$ will play a special role; we 
call it $\tau_0$.  Let $A$ be the affine $11$-space in $P\forms$
consisting of the $G_{0,u_2,\dots,u_{12}}$ of
\eqref{eq-family-of-Gs-for-A11}.  Let $B$ be 
the projective space spanned by $A$ and $T$, and $\Bhat$ be its proper
transform.  $\Bhat$ is where most of our work will take place.  Near
$(T,\tau_0)$, $\widehat{P\forms}$ is a product $\Bhat\times B^{22}$,
in a sense made precise by the following lemma.  In order to state it,
we observe that the
stabilizer of $(T,\tau_0)$ in $PG$ is $\C\semidirect\C^*$, of
codimension~$22$.  Therefore there exists a small $B^{22}\sset PG$
transverse to $\C\semidirect\C^*$ at $1\in PG$.

\begin{lemma}
\label{lem-Bhat-transverse-to-orbit}
The map $\Bhat\times B^{22}\to\widehat{P\forms}$ given by
$(b,g)\mapsto b.g$ is a local diffeomorphism at
$\bigl((T,\tau_0),1\bigr)$. 
\end{lemma}

\begin{proof}
We write $Y$ for the $PG$-orbit of $(T,\tau_0)$.  The lemma amounts
to the transversality of $Y$ and $\Bhat$ in
$\widehat{P\forms}$ at $(T,\tau_0)$.  Since $B$ is transverse to the
chordal cubic 
locus at $T$, it suffices to prove that $Y\cap E_T$ and $\Bhat\cap
E_T$ are transverse in $E_T$ at $(T,\tau_0)$.  We use $u_1,\dots,u_{12}$
as coordinates around $(T,\tau_0)$ in $E_T$ as in the proof of
lemma~\ref{lem-Deltahat-intersect-E}.  Then $\Bhat\cap E_T$ has equation $u_1=0$.  
$Y\cap E_T$ is the curve consisting of binary 12-tuples
$$
(s-\lambda)^{12}=s^{12}+12\lambda s^{11}+\dots+\lambda^{12}\;,
$$
which passes through $(T,\tau_0)$ when $\lambda=0$ and is transverse to
$\Bhat$ there because the $s^{11}$ coefficient (i.e., the
$u_1$-coordinate) is linear in $\lambda$. 
\end{proof}

\begin{remark}
It doesn't matter for us, but we note that 
$\Bhat$ and $Y$ are not transverse everywhere.  Since $Y\cap E_T$ is a
rational normal curve of degree~12 and $\Bhat\cap E_T$ is a hyperplane in
$E_T$, they intersect in 12 points, counted with multiplicity.
Besides $(T,\tau_0)$, the only place they intersect is at
$\lambda=\infty$, so they  make 11th-order contact there.
\end{remark}

The analogues of theorems~\ref{thm-discr-along-E-topological-model} and~\ref{thm-discr-along-E-analytic-model} in this lower-dimensional
setting are the the following.  It turns out (see the proof of
theorem~\ref{thm-discr-along-E-topological-model}) that restricting attention to $\tau_0$, rather than
treating general $\tau$, is sufficient.

\begin{theorem}
\label{thm-homeomorphism-model-lower-dimensional}
There exists a neighborhood $U'$ of $(T,\tau_0)$ in $\Bhat$ which is
homeomorphic to $B^1\times(U'\cap E)$, such that $U'\cap E$
corresponds to $\{0\}\times (U'\cap E)$ and $U'\cap\Dhat$ to
$B^1\times(U'\cap E\cap\Dhat)$.
\end{theorem}

\begin{theorem}
\label{thm-analytic-model-lower-dimensional}
There exists a neighborhood $U'$ of $(T,\tau_0)$ in $\Bhat$ which is
diffeomorphic to $B^1\times B^{11}$, with $U'\cap E$ corresponding to
$\{0\}\times B^{11}$, such that the following holds.  We write
$\pi_{\tilde U'}:\tilde U'\to U'$ for the $6$-fold cover of $U'$,
branched over $U'\cap E$, and $(T,\tau_0)\sptilde$ for the preimage
therein of $(T,\tau_0)$.  There is a neighborhood 
$\tilde V'$ of $(T,\tau_0)\sptilde$ in $\tilde U'$,
and a neighborhood $W'$ of $(T,\tau_0)$ in $\Bhat\cap E$, such that
$\tilde V'$ is diffeomorphic to $B^1\times W'$, such that
$$
\tilde V'\cap\pi_{\tilde U'}^{-1}(E)=\{0\}\times W'
$$
and
$$
\tilde V'\cap\pi_{\tilde U'}^{-1}(\Dhat)=B^1\times\bigl(W'\cap\Dhat\bigr)\;.
$$
\end{theorem}

These theorems describe $\Dhat$ in a neighborhood of $(T,\tau_0)$ in $\Bhat\cap E$
in terms of its intersection with $E$.  Therefore we need to
understand $\Bhat\cap E\cap\Dhat$:

\begin{lemma}
\label{lem-versal-deformations-inside-E-cap-Bhat}
Suppose $(T,\tau)\in\Bhat\cap E$, and that none of the points of
$\tau$ is $[0,0,0,0,1]$.  Let $m$, $n_1,\dots,n_m$, $K_1,\dots,K_m$
and $\D_1,\dots,\D_m$ be as in
theorem~\ref{thm-discr-along-E-topological-model}.  Let
$N'=11-\sum\dimension K_i$.  Then there is a neighborhood $Z$ of
$(T,\tau)$ in $\Bhat\cap E$ diffeomorphic to $K_1\times\dots\times
K_m\times B^{N'}$, such that $Z-\Dhat$ corresponds to
$$
(K_1-\D_1)\times\dots\times(K_m-\D_m)\times B^{N'}\;.
$$
\end{lemma}

\begin{proof}
Using coordinates $u_2,\dots,u_{12}$ around $(T,\tau_0)$ as in the
proof of lemma~\ref{lem-Deltahat-intersect-E}, and
parameterizing $R_T$ by $s$ as in \eqref{eq-parameterization-of-RNC}, the $\tau$'s treated in this
lemma are parameterized by the functions
$$
f(s)=s^{12}+u_2s^{10}+\dots+u_{11}s+u_{12}\;,
$$
i.e., as the monic polynomials with root sum equal to zero.  The lemma
amounts to the assertion that any singular function in this family
admits a simultaneous versal deformation of all its singularities,
within the family.  The family of functions
\begin{equation}
\label{eq-foo-3}
(s-s_0)^n+c_2(s-s_0)^{n-2}+\dots+c_{n-1}(s-s_0)+c_{n}
\end{equation}
provides a versal deformation of $(s-s_0)^n$, with every member of the
family having the same root sum, namely $ns_0$.  Given $f$ as above,
we may take a product of terms like \eqref{eq-foo-3}, one for each
singularity of $f$.  This obviously provides a simultaneous versal
deformation, and the root sum of any member of the family is that of
$f$, namely $0$.  Therefore the family lies in $\Bhat\cap E$.
\end{proof}

\begin{proof}[Proof of
    theorem~\ref{thm-discr-along-E-topological-model}, given
    theorem~\ref{thm-homeomorphism-model-lower-dimensional};] 
We first claim that $(T,\tau_0)$ has a neighborhood $Z$ in
$\Bhat\cap E$ and a neighborhood $U$ in $\widehat{P\forms}$, such that
$U$ is homeomorphic to $B^1\times Z\times B^{22}$, with $U\cap E$
corresponding to $\{0\}\times Z\times B^{22}$ and $U\cap\Dhat$ to
$B^1\times(Z\cap\Dhat)\times B^{22}$.  To get this, apply
theorem~\ref{thm-homeomorphism-model-lower-dimensional} to obtain $U'\sset\Bhat$ with the properties stated there, and
set $Z$ equal to $U'\cap E$.  By shrinking $U'$ and $B^{22}\sset
PG$ if necessary, we may suppose by lemma~\ref{lem-Bhat-transverse-to-orbit} that $U'\times
B^{22}\to\widehat{P\forms}$ is a diffeomorphism onto a neighborhood of
$(T,\tau_0)$ in $\widehat{P\forms}$, which we take to be $U$.  Then
$$
U
\isomorphism 
U'\times B^{22}
\isomorphism
B^1\times(U'\cap E)\times B^{22}
=
B^1\times Z\times B^{22}\;.
$$
Here, the first `$\isomorphism$' is a diffeomorphism and the second
is a homeomorphism.  Also, $U\cap E$ corresponds to $(U'\cap E)\times
B^{22}=\{0\}\times Z\times B^{22}$, and 
$U\cap\Dhat$ to $(U'\cap\Dhat)\times
B^{22}=B^1\times(Z\cap\Dhat)\times B^{22}$.

Now we observe that by the nature of the claim, the same conclusions
apply when $\tau_0$ is replaced by any $\tau\in Z$.  Now we use
lemma~\ref{lem-versal-deformations-inside-E-cap-Bhat}, which describes
$Z\cap\Dhat$ (possibly after shrinking $Z$ to a smaller neighborhood
of $(T,\tau)$, which shrinks $U$).  The result
is that $U$ 
is homeomorphic to
$$
B^1\times K_1\times\dots\times K_m\times B^{N'}\times B^{22}\;,
$$
such that 
$$
U\cap E=\{0\}\times K_1\times\dots\times K_m\times B^{N'}\times B^{22}
$$
and
$$
U-\Dhat=B^1\times(K_1-\D_1)\times\dots\times(K_m-\D_m)\times
B^{N'}\times B^{22}\;.
$$
Therefore theorem~\ref{thm-discr-along-E-topological-model} holds for
$(T,\tau)$.   Obviously, it also 
holds for and $(T,\tau')\in E_T$ that is equivalent to  some
$(T,\tau)\in Z$ under the stabilizer $\PGL(2,\C)$ of $T$ in $PG$.  It
is easy to see that 
this accounts for every $\tau'$, so the proof of
theorem~\ref{thm-discr-along-E-topological-model} is 
complete. 
\end{proof}

The proof of theorem~\ref{thm-discr-along-E-analytic-model}, given theorem~\ref{thm-analytic-model-lower-dimensional}, is essentially the
same.  Therefore it remains only to prove theorems~\ref{thm-homeomorphism-model-lower-dimensional} and~\ref{thm-analytic-model-lower-dimensional}.
We will treat theorem~\ref{thm-analytic-model-lower-dimensional} first.

\begin{lemma}
\label{lem-Bhat-transform-equals-transform-Bhat}
In a neighborhood of $(T,\tau_0)$, $\widehat{B\cap\D}=\Bhat\cap\Dhat$.
\end{lemma}

\begin{proof}
By lemma~\ref{lem-Bhat-transverse-to-orbit}, we may choose a neighborhood $U'$ of $(T,\tau_0)$ in
$\Bhat$, and shrink $B^{22}\sset PG$ if necessary, so that $U'\times
B^{22}\to\widehat{P\forms}$ is a diffeomorphism onto a neighborhood of
$(T,\tau_0)$ in $\widehat{P\forms}$.  Under this identification, $\Dhat-E$
corresponds to $\bigl((\Dhat-E)\cap U'\bigr)\times B^{22}$.  To get
$\Bhat\cap\Dhat$, we take the closure and then intersect with
$U'\times\{\hbox{point}\}$, and to get $\widehat{B\cap\D}$, we
intersect with $U'\times\{\hbox{point}\}$ and then take the closure.
Clearly, both give the same result.
\end{proof}

The following technical lemma is impossible to motivate without seeing
its use in the proof of theorem~\ref{thm-analytic-model-lower-dimensional}; the reader should skip it and
refer back when needed.

\begin{lemma}
\label{lem-rigidity-of-discriminant}
Suppose $\d(u_2,\dots,u_{12})$ is a quasihomogeneous polynomial of
weight $132$, where $\wt(u_i)=i$.  Suppose also that the $u_{12}^{11}$
and $u_{12}^{11-i}u_{11}^iu_i$ $(i=11,\dots,2)$ terms of $\d$ have
nonzero coefficients. Suppose $g:(\C^{11},0)\to(\C^{11},0)$ is the
germ of a diffeomorphism such that $\d\circ g$ has no terms of
weight~$<132$.  Then $g$ preserves the weight filtration, in the sense
that
\begin{equation}
\label{eq-vs-in-term-of-us}
u_i\circ g=c_iu_i+p_i(u_2,\dots,u_{12})+q_i(u_2,\dots,u_{12})
\end{equation}
for each $i$, where $c_i$ is a nonzero constant, $p_i$ is 
quasihomogeneous of weight~$i$ with no linear terms, and $q_i$
is an analytic function whose power series expansion has only terms of
weight~$>i$. 
\end{lemma}

\begin{proof}
We write $v_i$ for $u_i\circ g$, regarded as a function of
$u_2,\dots,u_{12}$.  One obtains $(\d\circ g)(u_2,\dots,u_{12})$ by
beginning with $\d(u_2,\dots,u_{12})$ and replacing each $u_i$ by
$v_i(u_2,\dots,u_{12})$.
This leads to a big mess, with the coefficients of $\d\circ g$
depending on those of $\d$ and the $v_i$ in a complicated way.
Nevertheless, there are some coefficients of $\d\circ g$ to which only
one term of $g$ can contribute, and this will allow us to deduce that
various coefficients of the $v_i$ vanish.  To be able to speak
precisely, we make the following definitions.  When we refer to a term
or monomial of $\d$
(resp. $v_i$), we mean a monomial whose 
coefficient in $\d$ (resp. $v_i$) is nonzero.  If $m$ is a monomial
$u_{i_1}\cdots u_{i_n}$, and $\mu_j(u_2,\dots,u_{12})$ is a monomial
of $v_{i_j}$ for each $j$, then we say that $m$  produces the
monomial 
$$
\mu=\mu_1(u_2,\dots,u_{12})\cdots\mu_n(u_2,\dots,u_{12})\;.
$$
For example, if $v_{12}=u_{12}+u_2^2$ and $v_{11}=u_{11}+u_2$, then
$m=u_{12}^2u_{11}$ produces the monomials $u_{12}^2u_{11}$,
$u_{12}u_2^2u_{11}$, $u_2^4u_{11}$, $u_{12}^2u_2$, $u_{12}u_2^3$ and
$u_2^5$.  When we wish to be more specific, we say that $m$ produces
$\mu$ by the substitution $\mu_j$ for each $u_{i_j}$.  Continuing the
example, we would say that $u_{12}^2u_{11}$ produces $u_{12}u_2^3$ by
substituting $u_{12}$ for one factor $u_{12}$ of $m$, $u_2^2$ for the
other, and $u_2$ for the factor $u_{11}$.
We note that even if a monomial $m$ of $\d$ produces a monomial $\mu$, the
coefficient of $\mu$ in $\d\circ g$ (or even in $m\circ g$) may
still be zero, because of possible cancellation.

We will prove the lemma by proving the following assertions $\A_d$ by
increasing induction on $d$, and we will prove each $\A_d$ by proving
the assertions $\B_{d,i}$ by decreasing induction on $i$.  
$\A_0$ is vacuously true.  The first two steps in the induction,
$\B_{1,12}$ and $\B_{1,11}$, require special treatment.

Assertion $\A_d$: No $v_i$ has any term of degree${}\leq d$ and
weight${}<i$.

Assertion $\B_{d,i}$: $v_i$ has no term of degree $d$ and
weight${}<i$.

Proof of $\B_{1,12}$:  if $v_{12}$ has a term $u_{j<12}$, then the
monomial $u_{12}^{11}$ of $\d$ produces $\mu=u_j^{11}$, of
weight${}<132$.  No other monomial of $\d$ can produce $\mu$, because
the only monomial of weight $132$ and degree${}<12$ is
$u_{12}^{11}$.  Since $\d\circ g$ has no term of weight${}<132$,
$v_{12}$ cannot have a term $u_{j<12}$, proving $\B_{1,12}$.

Proof of $\B_{1,11}$:  if $v_{11}$ has a term $u_{j<11}$, then the
monomial $u_{11}^{12}$ of $\d$ produces $\mu=u_j^{12}$, of
weight${}<132$.  We claim that no other monomial $m$ of $\d$ can
produce $\mu$.  If $\deg m<12$ then $m=u_{12}^{11}$, and $\B_{1,12}$
implies that the only degree~12 monomials that $m$ can produce have the
form
$$
(\hbox{a quadratic monomial})\cdot u_{12}^{10}\neq\mu\;.
$$
If $\deg m>12$ then $\deg m>\deg\mu$, so $m$ cannot produce $\mu$.  If
$m$ has degree~12, then $m$ could only produce $\mu$ by a linear
substitution for each factor.  If $m$ has
a factor $u_{12}$, then any monomial that $m$ produces by such a
substitution is divisible by $u_{12}$, by $\B_{1,12}$.  We have proven
that a monomial $m$ of $\d$ that produces $\mu$ has degree~12 and no
$u_{12}$ factors.  Since the average weight of the factors
is $132/12=11$, every factor is $u_{11}$, so $m=u_{11}^{12}$,
proving our claim.  Since $\d\circ g$ has no term $u_{j<11}^{12}$,
$v_{11}$ has no term $u_{j<11}$.  This proves $\B_{1,11}$.

Having proven $\B_{1,12}$ and $\B_{1,11}$, we observe that $v_{12}$
has a linear term because $g$ is a diffeomorphism, and since $v_{12}$
has no term $u_{11},\dots,u_2$, it does have a term $u_{12}$.
Similarly, using $\B_{1,11}$ and the fact that $v_{11}$ and $v_{12}$
have linearly independent linear parts, we see that $v_{11}$ has a
$u_{11}$ term.  We will use these facts in the rest of the proof.

Proof that $\A_d$ and $\B_{d,12},\dots,\B_{d,i+1}$ imply $\B_{d,i}$
(for $i=11,\dots,2$, except for $\B_{1,11}$, treated above): If $v_i$
had a term $t$ of degree $d$ and weight${}<i$, then the monomial
$u_{12}^{11-i}u_{11}^iu_i$ of $\d$ would produce the term
$\mu=u_{12}^{11-i}u_{11}^it$, of degree $11+d$ and weight${}<132$.
Here we are using that $v_{12}$ has a $u_{12}$ term and $v_{11}$ has a
$u_{11}$ term.  We claim that no other monomial $m$ of $\d$ can
produce $\mu$; the argument is similar to the proof of $\B_{1,11}$,
just more complicated.  First suppose that $m$ has degree${}>12$.  The
only way that $m$ can produce a term of degree $11+d$ is by
substituting a monomial of degree${}<d$ for each factor.  By
$\A_{d-1}$, the resulting monomial will have weight at least that of
$m$ and therefore cannot be $\mu$.  Next, suppose that $m$ has
degree${}<12$, so that $m=u_{12}^{11}$.  If we substitute a term of
degree${}\leq d$ for each factor $u_{12}$, then by $\A_{d-1}$ and
$\B_{d,12}$, the resulting monomial has weight at least that of $m$,
hence is not $\mu$.  On the other hand, if we substitute a term of
degree${}>d$ for some factor $u_{12}$, then either that term has
degree exactly $d+1$ and $u_{12}$ is substituted for each of the other
$u_{12}$'s, or else the degree of the resulting monomial is more than
$11+d$.  In neither of these cases can the result be $\mu$,
because $\mu$ has fewer than $10$ factors $u_{12}$ and has degree
$11+d$.  Finally, suppose $m$ has degree~12.  If we substitute a term
of degree${}<d$ for each factor of $m$, then by $\A_{d-1}$, the
resulting monomial has weight at least that of $m$, so is not $\mu$.
So the only way $m$ could produce $\mu$ is by a degree~$d$
substitution for one of the 12 factors and linear substitutions for
the others.  If the exceptional factor is one of
$u_{12},\dots,u_{i+1}$, then by $\B_{d,12},\dots,\B_{d,i+1}$, the
resulting monomial again has weight at least that of $m$, so is not
$\mu$.  We have shown that if $m$ produces $\mu$ then $m$ has degree
$12$, has a factor $u_{j\leq i}$, and that all the other factors are
replaced by linear terms.  The latter condition implies that the
exponent of $u_{12}$ in $m$ is at most that of $\mu$, since (by
$\B_{1,12}$) $u_{12}\to u_{12}$ is the only possible linear
substitution for $u_{12}$.  Therefore
$$
m=u_{12}^{p\leq11-i}u_{j\leq i}\cdot(\hbox{$11-p$ factors, none of
  which is $u_{12}$}).
$$
The average weight of the $11-p$ other factors is
$$
\frac{132-12p-j}{11-p}\geq11\;,
$$
with equality only if $p=11-i$ and $j=i$.  
This forces both of these equalities to hold, and each of
the remaining factors to be $u_{11}$.
That is, $m=u_{12}^{11-i}u_{11}^iu_i$.  This proves our
claim that this is the only monomial of $\d$ that produces $\mu$.
Since $\mu$ is not a term of $\d\circ g$, $v_i$ has no term $t$,
proving $\B_{d,i}$.

Proof that $\A_{d-1}$ and $\B_{d,12},\dots,\B_{d,2}$ imply $\A_d$:
trivial.

Proof that $\A_{d-1}$ implies $\B_{d,12}$ (except for $\B_{1,12}$,
treated above):  if $v_{12}$ has a term $t$ of degree~$d$ and
weight${}<12$, then $u_{12}^{11}$ produces $\mu=u_{12}^{10}\,t$, of
degree $10+d$ and weight${}<132$.  We claim that no other monomial $m$
of $\d$ can produce $\mu$.  If $m$ has degree${}\geq12$ then the only
way $m$ can produce a monomial of degree $10+d$ is by replacing each
factor by a term of degree${}<d$.  By $\A_{d-1}$, the result of such a
substitution has weight at least that of $m$, so is not $\mu$.  Since
$u_{12}^{11}$ is the only monomial of weight $132$ and
degree${}<12$, we have proven our claim.  Since $\d\circ g$ has no
term $u_{12}^{10}\,t$, $u_{12}$ has no term $t$, proving $\B_{d,12}$. 

The induction proves $\A_1,\dots,\A_5$ successively, and $\A_5$
implies that the $v_i$ have no terms other than those in \eqref{eq-vs-in-term-of-us}.  To
see that $c_i\neq0$ for all $i$, apply the argument used to prove that
$v_{12}$ (resp. $v_{11}$) has a term $u_{12}$ (resp. $u_{11}$).
Namely, $\B_{1,12},\dots,\B_{1,i}$ plus the linear independence of the
linear parts of $v_{12},\dots,v_i$ imply that $v_i$ has a term $u_i$.
\end{proof}

\begin{proof}[Proof of theorem~\ref{thm-analytic-model-lower-dimensional}:]
We will need standard coordinates around $(T,\tau_0)$ in $\Bhat$ and in
the 6-fold cover.  Let $\gamma:\C\times A\to\Bhat$ be the map
lying over the map $\C\times A\to B$ given by
$(\lambda,a)\mapsto F+\lambda(a-F)$.  For $a\in A$, $\c(0,a)$ is the
point of $E_T$ corresponding to the limiting direction of the pencil
$\langle F,a\rangle$.  
In particular, $\c(0,G_{0,\dots,0})=(T,\tau_0)$.
We take $U'$ to be the image of $\c$.
If $\lambda\neq0$ then $F+\lambda(a-F)$ has an isolated singularity at
$P$, so it is not a chordal cubic.  
It follows that
$\c^{-1}(E)=\{0\}\times A$, verifying the claimed property of $U'$.
Then $\tilde U'=\C\times A$, with $\pi_{\tilde U'}:\tilde U'\to U'$
given by $(\lambda,a)\mapsto(\lambda^6,a)$.  For brevity, we will
write $\b$ for this map.

The first idea (of five) is to use a 1-parameter group to work out a
defining equation for $(\c\circ\b)^{-1}(\Dhat)$ in terms of a defining
equation for $\D\cap A$.  We write $\d'$ for the defining equation for
$\D\cap A$ with respect to the coordinates $u_2,\dots,u_{12}$.  Our
basic tool is the 1-parameter subgroup
\begin{equation}
\label{eq-1-parameter-group}
\sigma_\lambda:(x_0,\dots,x_4)\mapsto (\lambda^{-2}x_0,
\lambda^{-1}x_1,x_2,\lambda x_3, \lambda^2x_4)
\end{equation}
of $G$.  For $\lambda\neq0$, $(\c\circ\b)(\lambda,u_2,\dots,u_{12})$
lies in $\D$ if and only if its image under $\sigma_{\lambda}^{-1}$
does.  This image is
\begin{align*}
\bigl(\c\circ\b(\lambda,u_2,&\dots,u_{12})\bigr).\sigma_{\lambda}^{-1}\\
&{}=\bigl(F+\lambda^6x_4^3+\lambda^6u_2x_4x_3^2+\cdots+\lambda^6u_{12}x_0^3\bigr).\sigma_\lambda^{-1}\\
&{}=F+x_4^3+\lambda^2u_2x_4x_3^2+\dots+\lambda^{12}u_{12}x_0^3\;,
\end{align*}
which is the point of $A$ with coordinates
$(\lambda^2u_2,\dots,\lambda^{12}u_{12})$.  Therefore a defining
equation for $(\c\circ\b)^{-1}(\Dhat)$ in $\C^*\times A$ is
\begin{equation}
\label{eq-foo2}
\d'(\lambda^2u_2,\dots,\lambda^{12}u_{12})=0\;.
\end{equation}

The second idea is to determine the lowest-weight terms of $\d'$, with
respect to the weights $\wt(u_i)=i$, by using our knowledge of
$\Dhat\cap E$.  The closure in
$\C\times A$ of the variety \eqref{eq-foo2} meets $\{0\}\times A$ in the variety defined by
$\dlowest'(u_2,\dots,u_{12})$, where $\dlowest'$ consists of the
lowest-weight terms of $\d'$.  By lemma~\ref{lem-Bhat-transform-equals-transform-Bhat}, we have
$\widehat{B\cap\D}=\Bhat\cap\Dhat$ near $(T,\tau_0)$.  Therefore $\dlowest'$
defines $(\c\circ\b)^{-1}(E\cap\Dhat)\sset\{0\}\times A$.  By
lemma~\ref{lem-Deltahat-intersect-E}\eqref{item-Deltahat-intersect-E}, we know that $E_T\cap\Dhat$ consists of those
$(T,\tau)$ where $\tau$ has a multiple point.
This forces $\dlowest'$ to be a power of the
standard $A_{11}$ discriminant $\d$, which is defined by the property that
$\d(u_2,\dots,u_{12})=0$ if and only if
$$
s^{12}+u_2s^{10}+u_3s^9+\dots+u_{12}
$$
has a multiple root.  Since $\d$ is quasihomogeneous of weight~132,
we have
\begin{equation}
\label{eq-lowest-order-terms}
\d'(u_2,\dots,u_{12})=\d^{\,p}(u_2,\dots,u_{12})
+(\hbox{terms of weight${}>132p$})
\end{equation}
for some $p\geq1$.  Although it is not essential, we will soon see
that $p=1$.

The third idea is to use singularity theory to describe $\D\cap A$ in
a neighborhood of $G_{0,\dots,0}$, in terms of a different set of
local coordinates.  One can show that the threefold defined by $G_{0,\dots,0}$
has only one singularity, which lies at $P$ and has type $A_{11}$.
Furthermore, the family $A$ of threefolds provides a versal
deformation of this singularity.  Therefore there is a neighborhood
$W$ of $G_{0,\dots,0}$ in $A$ such that $\D\cap W$ is a copy of the standard
$A_{11}$ discriminant.  That is, there are analytic coordinates
$v_2,\dots,v_{12}$ on $W$, centered at $G_{0,\dots,0}$, such that $\D\cap W$ is
defined by $\d(v_2,\dots,v_{12})=0$.  
This tells us
that $p=1$, because $\d'(u_2,\dots,u_{12})=0$ and
$\d(v_2,\dots,v_{12})=0$ define the same variety, and therefore
vanish to the same order at the origin.
We suppose without loss
of generality that $W$ is the unit polydisk with respect to the
$v_i$ coordinate system.   

The fourth idea is to use a sort of rigidity of the $A_{11}$
discriminant.  Briefly: since the variety defined by
$\d(u_2,\dots,u_{12})=0$ is close to that defined by  $\d(v_2,\dots,v_{12})=0$,
the $u_i$ must be close the $v_i$.  This relationship between
coordinate systems will be crucial later in the proof.  To make
this idea precise, consider the diffeomorphism-germ $g$ of $W$ at
$G_{0,\dots,0}$ 
given by $u_i\circ g=v_i$.  Since the image of the locus
$\d(v_2,\dots,v_{12})=0$ is the locus $\d(u_2,\dots,u_{12})=0$,
and since the first of these is also the locus $\d'(u_2,\dots,u_{12})=0$,
we have $\d\circ g=\d'$.  Now, a computer calculation using Maple
\cite{maple} shows that the 11 terms of $\d$ specified by the hypothesis
of lemma~\ref{lem-rigidity-of-discriminant} are nonzero, and \eqref{eq-lowest-order-terms} implies
that $\d'$ has no terms of weight${}<132$.  The lemma then implies
\begin{equation}
\label{eq-v-in-terms-of-u}
v_i=c_iu_i+p_i(u_2,\dots,u_{12})+q_i(u_2,\dots,u_{12})\;,
\end{equation}
where the $c_i$, $p_i$ and $q_i$ have the properties stated there.  We will
actually need not this but rather its inverse, giving the $u_i$ in
terms of the $v_i$.  To compute this it suffices to work in the formal
power series ring; one  writes \eqref{eq-v-in-terms-of-u} as
$$
u_i=\frac{1}{c_i}v_i-\frac{1}{c_i}p_i(u_2,\dots,u_{12})
+\frac{1}{c_i}q_i(u_2,\dots,u_{12})\;,
$$
substitutes these expressions into themselves, and repeats this process
infinitely many times.  The result is
\begin{equation}
\label{eq-the-us-in-terms-of-the-vs}
u_i=c_i'v_i+p_i'(v_2,\dots,v_{12})+q_i'(v_2,\dots,v_{12})\;,
\end{equation}
where the $c_i'$, $p_i'$ and $q_i'$ satisfy the same conditions as the
$c_i$, $p_i$ and $q_i$, with respect to the weights $\wt(v_i)=i$.

The fifth idea is to combine the quasihomogeneous scalings of the
$u_i$ and $v_i$ to define a map $\a$ whose image will be the set
$\tilde V'$ whose existence is claimed by the lemma.  For every
$0<|\lambda|<1$ we define $\rhosublambda:W\to W$ to be the
quasihomogeneous scaling with
respect to the $v$-coordinates:
$$
\rhosublambda(v_2,\dots,v_{12})=(\lambda^2v_2,\dots,\lambda^{12}v_{12})\;.
$$
For every $\lambda\in\C^*$, we define $\eta_\lambda:A\to A$ to be the
quasihomogeneous scaling with respect the the $u$-coordinates:
$$
\eta_\lambda(u_2,\dots,u_{12})=(\lambda^2u_2,\dots,\lambda^{12}u_{12})\;.
$$
The $\eta_\lambda$ are related to the 1-parameter group \eqref{eq-1-parameter-group} by
\begin{equation}
\label{eq-relation-to-1-par-subgroup}
\c\circ\b\bigl(\lambda,\eta_{\lambda}(a)\bigr)
=\sigma_\lambda^{-1}(a)\;,
\end{equation}
which can be  verified by expressing both sides in terms of the
$u$-coordinates and expanding.
We define 
$$
\a:\bigl(B^1-\{0\}\bigr)\times W\to\bigl(B^1-\{0\}\bigr)\times A
$$
by $\a(\lambda,w)=(\lambda,\eta_{\lambda^{-1}}\rhosublambda(a))$.
It is easy to see that $\a$ is injective.

The first property of $\a$ is that the preimage of the discriminant
under it has a very simple form, namely
\begin{equation}
\label{eq-main-property-of-alpha}
(\c\circ\b\circ\a)^{-1}(\Dhat)=\bigl(B^1-\{0\}\bigr)\times(\D\cap W)\;.
\end{equation}
To see this, observe that 
$$
\c\circ\b\circ\a(\lambda,w)
=\c\circ\b\bigl(\lambda,\eta_{\lambda^{-1}}\rhosublambda(w)\bigr)
=\sigma_\lambda(\rhosublambda(w))\;,
$$
which lies in $\D$ if and only if $\rhosublambda(w)$ does, hence if
and only if $w$ does.  

The second property of $\a$ is that it extends to a
holomorphic map
$B^1\times W\to\C\times A$.  To prove this, it suffices by Riemann
extension to show that it has a continuous extension.  The key step is
to compute
$\lim_{\lambda\to0}\a(\lambda,w)$.
If $w\in W$, then its $v$-coordinates are
$v_2(w),\dots,v_{12}(w)$, and the $v$-coordinates of $\rhosublambda(w)$
are $\lambda^2v_2(w),\dots,$ $\lambda^{12}v_{12}(w)$.  Using \eqref{eq-the-us-in-terms-of-the-vs}, the
$u$-coordinates of $\rhosublambda(w)$ are
\begin{align*}
u_i\bigl(\rhosublambda(w)\bigr) ={}&
c_i'\lambda^iv_i(w)
+p_i'\bigl(\lambda^2v_2(w),\dots,\lambda^{12}v_{12}(w)\bigr)\\
&\qquad
+q_i'\bigl(\lambda^2v_2(w),\dots,\lambda^{12}v_{12}(w)\bigr)\\
{}={}&
\lambda^i\cdot\Bigl(c_i'v_i(w)+p_i'\bigl(v_2(w),\dots,v_{12}(w)\bigr)\Bigr)\\
&\qquad
+\bigl(\hbox{terms of degree${}>i$ in $\lambda$}\bigr)\;.
\end{align*}
Therefore the $u$-coordinates of $\eta_{\lambda^{-1}}\rhosublambda(w)$
are
\begin{align*}
u_i\bigl(\eta_{\lambda^{-1}}\rhosublambda(w)\bigr)
=\lambda^{-i}u_i\bigl(\rhosublambda(w)\bigr)
={}&
c_iv_i(w)+p_i'\bigl(v_2(w),\dots,v_{12}(w)\bigr)\\
&\qquad
+\hbox{terms involving $\lambda$}\;.
\end{align*}
The limit as $\lambda\to0$ obviously exists, and provides the desired
extension.  

The third property of $\a$ is that after shrinking $W$,
we may suppose that $\a:B^1\times W\to\C\times A$ is a diffeomorphism
onto a neighborhood of $(0,G_{0,\dots,0})$.  To see this we use the
fact that
$$
u_i\bigl(\a(0,w)\bigr)=c_i'v_i(w)+p_i'\bigl(v_2(w),\dots,v_{12}(w)\bigr)\;;
$$
since the $c_i'$ are nonzero and the $p_i'$ have no linear terms, we
may shrink $W$ so that $w\mapsto\a(0,w)$ is injective.  Then
$\a:B^1\times W\to\C\times A$ gives  a diffeomorphism from
$\{0\}\times W$ onto a neighborhood
of $(0,G_{0,\dots,0})$ in $\{0\}\times A\sset \tilde U'$.  We will
call this diffeomorphism $\a_0$.

Because $\a:B^1\times W\to\tilde U'$ is injective and $\tilde U'$ is
normal, $\a$ is a diffeomorphism onto a
neighborhood of $(T,\tau_0)\sptilde$ in $\tilde U'$.  We define
$\tilde V'$ to
be the image of $\a$.  Unwinding the definitions gives
\begin{equation}
\label{eq-foobar-1}
\bigl(\pi_{\tilde U'}\circ\a\bigr)^{-1}(E)=\{0\}\times W
\end{equation}
and
\begin{equation}
\label{eq-foobar-2}
\bigl(\pi_{\tilde U'}\circ\a\bigr)^{-1}(\Dhat)=B^1\times(W\cap\Dhat)\;.
\end{equation}
This is exactly what we want, except that $W$ is a subset of $A$, not
a subset of $\Bhat\cap E$.  However, $\a_0$ identifies $W$ with a
subset of $\{0\}\times A$, which is identified under $\c\circ\b$ with
a subset of $\Bhat\cap E$, which we take to be $W'$.  The
identification $W\isomorphism W'$ identifies $W\cap\Dhat$ with
$W'\cap\Dhat$, so we may replace $W$ by $W'$ in \eqref{eq-foobar-1} and \eqref{eq-foobar-2}.
This completes the proof.
\end{proof}

\begin{proof}[Proof of theorem~\ref{thm-homeomorphism-model-lower-dimensional}]

Theorem~\ref{thm-analytic-model-lower-dimensional} gives us an
analytic description of the 6-fold branched cover $\tilde V'$ of a
neighborhood of $(T,\tau_0)$ in $\Bhat$.  The idea is to take the
quotient by the deck group $\Z/6$ and see what we get.  Recall that
$\tilde U'\isomorphism\C\times A$, where $A\isomorphism\C^{11}$, and
the deck group is generated by $(\lambda,a)\mapsto(\lambda\zeta,a)$,
where $\zeta=e^{\pi i/3}$.  We will write $\xi$ for this map.
$(T,\tau_0)\sptilde$ is the point $(0,G_{0,\dots,0})\in\tilde U'$.
Also, $\a:B^1\times W\to\tilde U'$ is an embedding onto a neighborhood
$\tilde V'$ of $(T,\tau_0)\sptilde$, where $W$ is a polydisk around
$G_{0,\dots,0}$ in $A$.  We can't regard $\xi$ as a self-map of
$B^1\times W$, because $\tilde V'$ may not be a $\xi$-invariant subset
of $\tilde U'$.  However, we can take the intersection of the finitely
many translates of $\tilde V'$, and let $Z\sset B^1\times W$ be the
$\a$-preimage of this intersection.  The action of $\xi$ on $Z$ can be
worked out by using the definition of $\a$.  The result is
$\xi(0,w)=(0,w)$, and
$$
\xi(\lambda,w)=\bigl(\lambda\zeta,\rho_\lambda^{-1}\circ\rho_\zeta^{-1}\circ\eta_\zeta\circ\rho_\lambda(w)\bigr)
$$
for $\lambda\in B^1-\{0\}$.
Now, $\pi_{\tilde U'}\circ\a$ carries $Z/\langle\xi\rangle$
homeomorphically to a neighborhood of $(T,\tau_0)$ in $\Bhat$.
$(T,\tau_0)$ corresponds to the image of $(0,G_{0,\dots,0})$, and
$\Dhat\cap\Bhat$ to the image of
$Z\cap\bigl(B^1\times(W\cap\Dhat)\bigr)$.  Therefore our goal is to
describe $Z/\langle\xi\rangle$ and the image therein of
$Z\cap\bigl(B^1\times(W\cap\Dhat)\bigr)$.

To take the quotient $Z/\langle\xi\rangle$, we first observe that the
`slice of pie'
$$
\Sigma=\set{\lambda\in B^1}{\hbox{$\lambda=0$ or $\Arg \lambda\in[0,\pi/3]$}}
$$
is a fundamental domain for $\lambda\mapsto \lambda\zeta$ acting on
$B^1$, and the quotient $B^1/(\Z/6)$ is got by gluing one edge
$$
E_0=\set{\lambda\in B^1}{\hbox{$\lambda=0$ or $\Arg \lambda=0$}}
$$
to the other
$$
E_{\pi/3}=\set{\lambda\in B^1}{\hbox{$\lambda=0$ or $\Arg \lambda=\pi/3$}}
$$
in the obvious way.  Similarly, $Z/\langle\xi\rangle$ is homeomorphic
to $\bigl(Z\cap(\Sigma\times W)\bigr)/{\sim}$, where $\sim$ is the
equivalence relation that each $(s,w)\in Z\cap(E_0\times W)$ is identified
with $\xi(s,w)\in Z\cap(E_{\pi/3}\times W)$.  On the other hand,
consider $(\Sigma\times W)/{\approx}$, where $\approx$ is the
equivalence relation that each $(s,w)\in E_0\times W$ is identified
with $(s\zeta,w)$.  The key step in the proof is the definition of the
following map
$$
\Phi:\bigl(Z_0\cap(\Sigma\times W)\bigr)/{\sim}
\ \to\ 
(\Sigma\times W)/{\approx}\;,
$$
where $Z_0$ is a suitable neighborhood of $(0,G_{0,\dots,0})$ in $Z$,
defined below.   The map is $\Phi(\lambda,w)=(\lambda,w)$ for
$\lambda\in E_0$, and
$$
\Phi(\lambda,w)=
\Bigl(\lambda,\ \rho_{(\Arg \lambda)/(\pi/3)}^{-1}\circ
(\rho_\lambda^{-1}\rho_\zeta^{-1}\eta^\phantominverse_\zeta\rho^\phantominverse_\lambda)
\circ\rho_{(\Arg \lambda)/(\pi/3)}
(w)\Bigr)
$$ 
otherwise.  The definition of $Z_0$ is essentially the subset of
$Z$ on which the formula makes sense.  That is, if $(\lambda,w)\in
Z_0$ with $\lambda\in\Sigma- E_0$, then
$$
\eta^\phantominverse_\zeta\rho^\phantominverse_\lambda\rho^\phantominverse_{(\Arg \lambda)/(\pi/3)}(w)
\in
\rho^\phantominverse_\lambda\rho^\phantominverse_{(\Arg \lambda)/(\pi/3)}(W)\;.
$$

Assuming for a moment the existence of this set $Z_0$, it is easy to
complete the proof.  The key point is that the restriction of $\Phi$
to $Z_0\cap(E_0\times W)$ is the obvious inclusion into $E_0\times W$,
while the $\Phi$-image of $(\lambda,w)\in(E_{\pi/3}\times W)$ is
$\xi(\lambda\zeta^{-1},w)$.  Furthermore, for any $(\lambda,w)\in
Z_0\cap(\Sigma\times W)$, $\Phi(\lambda,w)$ lies in
$\Sigma\times(W\cap\Dhat)$ if and only if $(\lambda,w)$ itself does.
(This uses the fact that the $\rho_\lambda$ are quasihomogeneous
scalings of $W$ that preserve $W\cap\Dhat$.)  It follows that
$Z/\langle\xi\rangle$ is homeomorphic to a neighborhood of
$(0,G_{0,\dots,0})$ in $(\Sigma\times W)/{\approx}$, i.e., in
$B^1\times W$, with its intersection with $\Dhat$ being
$\bigl(\Sigma\times(\Dhat\times W)\bigr)/{\approx}$, i.e.,
$B^1\times(\Dhat\times W)$.  This implies the theorem.

All that remains is to construct $Z_0$.  First observe that $Z$
contains $B^1\times\{G_{0,\dots,0}\}\sset B^1\times W$.  Therefore
there exists a neighborhood $W_0$ of $G_{0,\dots,0}$ in $W$ such that
$B^1\times W_0\sset Z$.  Since $\xi$ is defined on $B^1\times W_0$, we
have $\eta_\zeta\rho_{\lambda'}(w)\in\rho_{\lambda'}(W)$ for all
$(\lambda',w)\in B^1\times W_0$.  Setting $\lambda'=\lambda\cdot(\Arg
\lambda)/(\pi/3)$, we see that $\Phi$ is defined on $B^1\times W_0$.
Now we set $Z_0$ equal to the intersection of the $\xi$-translates of
$B^1\times W_0$.
\end{proof} 

\section{Extension of the period map}
\label{sec-extension}

In this section we extend the period map
$g:P\forms_0\to\PGamma\backslash\ch^{10}$ defined in
section~\ref{sec-smoothmoduli} to a larger domain.  It turns out that
$g$ does extend to all of $P\forms_s$ but not to $P\forms_{ss}$.
Replacing $\PGamma\backslash\ch^{10}$ by its Baily-Borel
compactification $\overline{\PGamma\backslash\ch^{10}}$ allows us to
extend $g$ to most of $P\forms_{ss}$ but not all.  The problem is that
it does not extend to the chordal cubic locus.  We will see a proof of
this in section~\ref{sec-hodge-theory-chordal}, but for now we just
refer to this to motivate the blowing-up of the chordal cubic locus in
$P\forms$ to obtain $\widehat{P\forms}$, and extending $g$ to a
regular map
$(\widehat{P\forms})_{ss}\to\overline{\PGamma\backslash\ch^{10}}$.  In
this section we will construct the extension; we rely heavily on
lemma~\ref{lem-discriminant-away-from-E} and
theorems~\ref{thm-discr-along-E-topological-model}
and~\ref{thm-discr-along-E-analytic-model}, which describe the
local structure of the discriminant.  Recall that we write $\pi$ for
the projection $\widehat{P\forms}\to P\forms$, $E$ for the exceptional
divisor, $E_T\isomorphism P^{12}$ for the points of $E$ lying over a
chordal cubic $T$, and $(T,\tau)$ for a point in $E_T$, where $\tau$
is an unordered 12-tuple in the rational normal curve $R_T$ of which $T$ is
the secant variety.

The first thing we need to do is describe $(\widehat{P\forms})_s$ and
$(\widehat{P\forms})_{ss}$, using \cite{allcock-threefolds} and \cite{reichstein}.
To lighten the notation we will write just $\widehat{P\forms}_s$ and
$\widehat{P\forms}_{ss}$.  In order to discuss
GIT-stability we need to choose a line bundle on
$\widehat{P\forms}$.  The following lemma shows that this choice
doesn't matter very much; it follows directly from the considerations
of \cite[sec.~2]{reichstein}.

\begin{lemma}
\label{lem-choice-of-line-bundle-irrelevant}
For large enough $d$, the stable and semistable loci of
$\widehat{P\forms}$, with respect to the standard $\SL(5,\C)$-action on 
\begin{equation}
\label{eq-choice-of-line-bundle}
\mathcal{O}(-E)\tensor\pi^*\bigl(\mathcal{O}(d)\bigr)\;,
\end{equation}
are independent of $d$.
\qed
\end{lemma}


Our notation $\widehat{P\forms}_s$ and $\widehat{P\forms}_{ss}$ refers
to the linearization \eqref{eq-choice-of-line-bundle} for large enough
$d$.  Reichstein's work allows us to describe these sets explicitly:

\begin{theorem}
\label{thm-GIT-analysis}
Suppose $T$ is a cubic threefold not in the closure of the chordal
cubic locus, regarded as an element of $\widehat{P\forms}$.  Then
\renewcommand\theenumi{\roman{enumi}}
\begin{enumerate}
\item
\label{item-stable-threefolds}
$T$ is stable if and only if each singularity of $T$ has type
$A_1$, $A_2$, $A_3$ or $A_4$;
\item
\label{item-semistable-threefolds}
$T$ is semistable if and only if each singularity of $T$ has type $A_1,\dots,A_5$
  or $D_4$;
\item
\label{item-closed-orbits-of-semistable-threefolds}
$T$ is strictly semistable with closed orbit in
  $\widehat{P\forms}_{ss}$ if and only if $T$ is projectively equivalent to one
  of the threefolds defined by
$$
x_0x_1x_2+x_3^3+x_4^3
$$
or
$$
F_{A,B}=Ax_2^3+x_0x_3^2+x_1^2x_4-x_0x_2x_4+Bx_1x_2x_3\;,
$$
with $A,B\in\C$ and $4A\neq B^2$.
\end{enumerate}
Now suppose instead that $T$ is in the closure of the chordal cubic locus
but is not a 
chordal cubic.  Then every element of
$\pi^{-1}(T)\sset\widehat{P\forms}$ is unstable.  Finally, suppose
that $T$ is a
chordal cubic, and $\tau$ is an unordered $12$-tuple in the rational
normal curve $R_T$, so that $(T,\tau)\in E_T\sset\widehat{P\forms}$.  Then
\begin{enumerate}
\setcounter{enumi}{3}
\item
\label{item-stable-12-tuples}
$(T,\tau)$ is stable if and only if $\tau$ has no points of
  multiplicity${}\geq6$;
\item
\label{item-semistable-12-tuples}
$(T,\tau)$ is semistable if and only if $\tau$ has no points of
  multiplicity greater than~$6$;
\item
\label{item-closed-orbits-of-semistable-12-tuples}
$(T,\tau)$ is strictly semistable with closed orbit in
  $\widehat{P\forms}_{ss}$ if and only if $\tau$
  consists of two distinct points of multiplicity~$6$.
\end{enumerate}
Finally, the points \eqref{item-closed-orbits-of-semistable-12-tuples} of
$\widehat{P\forms}$ lie in the closure of the union of the orbits of the $T_{A,B}$ from
\eqref{item-closed-orbits-of-semistable-threefolds}.
\end{theorem}

\begin{remarks}
The first threefold described in
\eqref{item-closed-orbits-of-semistable-threefolds} has three $D_4$
singularities, and is the unique such cubic threefold.  The
2-parameter family $T_{A,B}$ 
really describes
only a 1-parameter set of orbits, because the projective equivalence
class is determined by the ratio $4A/B^2\in\cp^1-\{1\}$.  These
threefolds have exactly two singularities, both of type $A_5$, except
when $A=0$, when there is also an $A_1$ singularity.  Every cubic
threefold with two $A_5$ singularities is projectively equivalent to
one of these.  If $4A$ and $B^2$ were allowed to be equal and nonzero,
then $F_{A,B}$ would define a chordal cubic.  All of these assertions are
proven in section~5 of \cite{allcock-threefolds}.
\end{remarks}

\begin{proof}
Throughout the proof, we will write $L\sset P\forms$ for the closure
of the chordal cubic locus.  By theorems~2.1 and~2.3 of
\cite{reichstein}, a point of $\widehat{P\forms}-E$ is unstable as an
element of $\widehat{P\forms}$ if and only if either (a) it is
unstable as an element of $P\forms$, or (b) it is GIT-equivalent in
$P\forms$ to an element of $L$.  Referring to the stability of cubic
threefolds, given by theorems~1.3 and~1.4 of
\cite{allcock-threefolds}, this says that $\widehat{P\forms}_{ss}-E$
is the set of $T$'s having no singularities of types other than
$A_1,\dots,A_5$ and $D_4$.  This justifies
\eqref{item-semistable-threefolds}.  The same theorems of
\cite{reichstein} say that a point of $\widehat{P\forms}-E$ is stable
as an element of $\widehat{P\forms}$ if and only if it is stable as an
element of $P\forms$.  Referring again to \cite{allcock-threefolds},
this says that $\widehat{P\forms}_s-E$ is the set of $T$ having no
singularities of types other than $A_1,\dots,A_4$, justifying
\eqref{item-stable-threefolds}.  Now we prove
\eqref{item-closed-orbits-of-semistable-threefolds}.  If
$T\in\widehat{P\forms}_{ss}-E$ is not in $\widehat{P\forms}_s$, then
it has a singularity of type $A_5$ or $D_4$.  Then theorem~1.3(i,ii) of \cite{allcock-threefolds} implies that $T$ is GIT-equivalent
in $P\forms$, hence in $\widehat{P\forms}$, to one of the threefolds
given in \eqref{item-closed-orbits-of-semistable-threefolds}.
Theorem~1.2 of \cite{allcock-threefolds} implies that the threefolds
given explicitly in \eqref{item-closed-orbits-of-semistable-threefolds}
have closed orbits in $P\forms_{ss}$; since the orbits miss $L$, they are
also closed in $\widehat{P\forms}_{ss}$.  It follows that these orbits are
the only orbits in $\widehat{P\forms}_{ss}-E$ that are strictly semistable and
closed in $\widehat{P\forms}_{ss}$.  This justifies
\eqref{item-closed-orbits-of-semistable-threefolds}.

Now suppose $T\in L$.  If $T$ is not a chordal cubic then it is
unstable by theorem~1.4(i) of \cite{allcock-threefolds}, so every
point of $\widehat{P\forms}$ lying over $T$ is unstable by theorem~2.1
of \cite{reichstein}.  It remains only to discuss stability of pairs
$(T,\tau)$ with $T$ a chordal cubic.  Our key tool is theorem~2.4 of
\cite{reichstein}.  This says that $(T,\tau)$ is unstable if and only
if it lies in the proper transform of the set of cubic threefolds that
are GIT-equivalent to chordal cubics.  So our job is to determine this
proper transform.  If $\tau$ has a point of multiplicity${}>6$, then
by lemma~\ref{lem-Deltahat-intersect-E}\eqref{item-limits-of-Ans} it
is a limit of threefolds having $A_{n>5}$ singularities.  Since
$T$ lies in $P\forms_{ss}$ and $P\forms_{ss}$ is open in $P\forms$,
$(T,\tau)$ is a limit of semistable threefolds having $A_{n>5}$
singularities.  By theorem~1.3 of \cite{allcock-threefolds}, such
threefolds are GIT-equivalent to chordal cubics.  Then Reichstein's
theorem~2.4 shows that $(T,\tau)$ is unstable.  Reichstein's theorem
also asserts that $(T,\tau)\in\widehat{P\forms}_{ss}$ is non-stable if and only if it lies in
the proper transform of $P\forms_{ss}-P\forms_s$.  If $\tau$ has a
point of multiplicity~$6$, then 
lemma~\ref{lem-Deltahat-intersect-E}\eqref{item-limits-of-Ans} shows
that 
$(T,\tau)$ is a limit of semistable threefolds having $A_5$
singularities, so it is
not stable.  This justifies the `if' parts of
\eqref{item-stable-12-tuples} and \eqref{item-semistable-12-tuples}.

Now, suppose $\tau$ has no point of multiplicity${}>6$.  Since
$P\forms_{ss}$ is open, $T$ has a neighborhood $U\sset P\forms$ with
every member of $U-L$ having only $A_n$ and $D_4$ singularities.  (In
fact, $D_4$ singularities can be excluded, but this doesn't matter
here.)  By lemma~\ref{lem-discriminant-away-from-E}, every member of $U-L$ admits in $P\forms$  a
simultaneous versal deformation of all its singularities.  If some
member of $U-L$ had an $A_{n\geq6}$ singularity, then at some point 
of $U-L$, $\D$ would be locally modeled on the $A_n$ discriminant
(times a ball of the appropriate dimension).  On
the other hand, it follows from theorem~\ref{thm-discr-along-E-analytic-model} that after shrinking $U$
we may suppose that at every point of $U-L$, $\D$ is locally modeled
on 
\begin{equation}
\label{eq-foo3}
\bigcup_{i=1}^m K_1\times\dots\times K_{i-1}\times\D_i\times
K_{i+1}\times\dots\times K_m\times B^N\;,
\end{equation}
where the notation is as in lemma~\ref{lem-discriminant-away-from-E}.
In particular, the $\D_i\sset K_i$ are copies of the $A_k$ discriminants for
various $k$'s that are at most $5$.  Since \eqref{eq-foo3} is not a copy of an
$A_{n\geq6}$ discriminant, $U$ contains no points with an $A_{n\geq6}$
singularity, hence no points GIT-equivalent to chordal cubics.  By
theorem~2.4 of \cite{allcock-threefolds}, $(T,\tau)$ is not unstable,
which is to say that it is semistable.
This proves the `only if' part of \eqref{item-semistable-12-tuples}.  The same argument, using
the fact that members of $\widehat{P\forms}_{ss}-E$ have $A_5$
or $D_4$ singularities, proves the `only if' part of \eqref{item-stable-12-tuples}.

Finally, if $\tau$ has a point of multiplicity~$6$, then $(T,\tau)$'s
orbit closure in $E_T\cap\widehat{P\forms}_{ss}$ contains $(T,\tau')$,
where $\tau'$ has two points of multiplicity~$6$.  This is a classical
fact about point-sets in $P^1$.  This proves
\eqref{item-closed-orbits-of-semistable-12-tuples}.  To prove the last
claim of the theorem, just observe that the restrictions of the
$F_{A,B}$ in \eqref{item-closed-orbits-of-semistable-threefolds} to
the singular locus of the standard chordal cubic  (defined by
$F_{1,-2}$) consists of $[1,0,0,0,0]$ and $[0,0,0,0,1]$, each with
multiplicity~$6$.  Let $A\to 1$ and $B\to -2$.
\end{proof}

Now that we know how much to enlarge the domain of $g$, we will
construct the extension.  This relies on an analysis of the local
monodromy group at a point of $\widehat{P\forms}$, by which we mean the
following.  In section~\ref{sec-smoothmoduli} we considered the local
system $\Lambda(\V_0)$ over $\forms_0$ and its associated local system
$\ch(\V_0)$ of complex hyperbolic spaces.  Now, $\Lambda(\V_0)$ does
not descend to a local system on $P\forms_0$, but $\ch(\V_0)$ does,
because the scalars $\{I,\w I,\wbar I\}\sset\GL(5,\C)$ act on each
$\Lambda(V)$ by scalar multiplication.  After fixing a basepoint
$F\in\forms_0$, we defined
$$
\rho:\pi_1(\forms_0,F)\to\Gamma(V):=\aut \Lambda(V)
$$
to be the monodromy of $\Lambda(\V_0)$.  
Analogously, we define, for $T\in P\forms_0$, 
\begin{equation}
\label{eq-defn-of-projective-monodromy-rep}
P\rho:\pi_1(P\forms_0,T)\to\PGamma(V)\sset\isom\bigl(\ch(V)\bigr)\;.
\end{equation}

Henceforth, all references to monodromy refer to $P\rho$ unless
otherwise stated.
In the arguments below, we will compute the monodromy of various
elements of $\pi_1(P\forms_0)$.  For convenience we will perform
various monodromy calculations with roots of $\Lambda(V)$, but these
could all be rephrased in terms of elements of $P\Gamma(V)$.

Now suppose $T_0$ or $(T_0,\tau_0)$ is an element of
$\widehat{P\forms}$ and $U$ is a suitable small neighborhood of
it; for example, $U$ could be as in
lemma~\ref{lem-discriminant-away-from-E} or
theorem~\ref{thm-discr-along-E-topological-model}.  By the local
fundamental group we mean $\pi_1\bigl(U-(\widehat\D\cup E),T\bigr)$,
where $T$ is a basepoint.  By the local monodromy action we mean the
restriction of $P\rho$ to the local fundamental group, and by the
local monodromy group we mean the image of this homomorphism in
$\PGamma(V)$.  We will see that $\widehat{P\forms}_s$ is exactly
the subset of $\widehat{P\forms}_{ss}$ where the local monodromy group
is finite.  In order to establish this, we will need to know the
monodromy around a meridian of $E$:

\begin{lemma}
\label{lem-monodromy-around-E}
Suppose $\c$ is a meridian around $E$ in $\widehat{P\forms}$, $T$ is a
point of $\c$, and $P\rho(\c)$ is the monodromy action of $\c$ on
$\ch(V)$. 
Then there is a direct sum decomposition
$$
\Lambda(V)=\Lambda_1\oplus\Lambda_{10}\;,
$$ where $\Lambda_1$ is the span of a norm $3$ vector $s$, $P\rho(\c)$
acts on $\ch(V)$ as a hexaflection in $s$, and $\Lambda_{10}$ is
isometric to the sum of the last three summands in
\eqref{eq-inner-product-matrix}.
\end{lemma}

This lemma resembles
lemma~\ref{lem-meridians-act-by-complex-reflections}; each shows that
a certain monodromy action is a complex reflection in a norm~3 vector
of $\Lambda(V)$.  But there is an essential difference.  We have
already defined a root of $\Lambda(V)$ to be any norm~3 vector $r$; we
refine the language by calling $r$ nodal or chordal root according to
whether $\ip{r}{\Lambda(V)}$ it $\theta\E$ or $3\E$.  It is easy to
see that every root is either nodal or chordal.
Lemma~\ref{lem-monodromy-around-E} asserts that the monodromy of a
meridian around $E$ is a hexaflection in a chordal root.
Lemma~\ref{lem-meridians-act-by-complex-reflections} asserts that the
monodromy of a meridian around $\Dhat$ is a triflection in a root, and
a simple argument shows that this root must be nodal.  (Namely, by
considering a threefold with an $A_2$ singularity, one finds two
meridians of $\Dhat$, which by
lemma~\ref{lem-local-monodromy-noncyclic} act by the $\w$-reflections
in linearly independent roots $r$ and $r'$, and satisfy the braid
relation.  This relation forces $\bigl|\ip{r}{r'}\bigr|=\sqrt3$, so
$\ip{r}{r'}$ is a unit times $\theta$.)  The `nodal' and `chordal'
language reflects the fact that these monodromy transformations arise
by considering a degeneration to a nodal threefold or to a chordal
cubic.  We caution the reader that while it is true that every nodal
(resp. chordal) root of $\Lambda(V)$ comes from a nodal
(resp. chordal) degeneration, we have not yet proven it.  In
theorem~\ref{lem-transitivity-on-nodal/chordal-hyperplanes} we show
that $\Gamma$ is transitive on nodal and chordal roots of $\Lambda$.
(The proof of
theorem~\ref{lem-transitivity-on-nodal/chordal-hyperplanes} is
independent of the rest of the paper, so it could be read at this
point.)

\begin{proof}[Proof of lemma~\ref{lem-monodromy-around-E}:]
Let $T_0$ be the standard chordal cubic, and $\tau_0$ a 12-tuple in
$R_{T_0}$ concentrated at one point.  By
theorem~\ref{thm-discr-along-E-topological-model}, 
the local fundamental group at $(T_0,\tau_0)$ is $\Z\times
B_{12}$, where $\c$ is a generator of $\Z$ and we write
$a_1,\dots,a_{11}$ for standard generators for the braid
group.  By lemma~\ref{lem-meridians-act-by-complex-reflections}, the
$a_i$ act on $\ch(V)$ as triflections, and the 1-dimensional
eigenspaces of (lifts of the $a_i$ to) $\Lambda(V)$ are spanned by
vectors $r_i$ of norm~3.  We take $\Lambda_{10}$ to be the span of the
$r_i$.  Following the proof of theorem~\ref{thm-inner-product-matrix}
shows that $\Lambda_{10}$ is a copy of the direct sum of the last
three summands of \eqref{eq-inner-product-matrix}, that $\Lambda_{10}$
is a summand of $\Lambda(V)$, and that $\Lambda_{10}^\perp$ is spanned
by a vector of norm~$3$.  We write $s$ for such a vector and
$\Lambda_1$ for its span.  Since $\c$ commutes with the $a_i$, any
lift of $P\rho(\c)$ to $\Lambda(V)$ multiplies each $r_i$ by a scalar.
Since $r_i\cdot r_{i+1}\neq0$, it multiplies all the $r_i$ by the same
scalar, so that it acts on $\Lambda_{10}$ as that scalar.  Therefore
$P\rho(\c)$ acts on $\ch(V)$ as a complex reflection in $s$ of order
$2$, $3$ or $6$ (or acts trivially).

Now we show that $P\rho(\c)$ has order~6.  We may find a neighborhood
of $E-\widehat\D$ in $\widehat{P\forms}$ which is a disk bundle over
$E-\widehat\D$.  (We do not need all the fibers to ``have the same
radius''---they can shrink as one approaches $\widehat\D$.)  We choose
a ball $B$ around $T_0$ in the chordal cubic locus, and write $N$ for
the restriction of this disk bundle to $\pi^{-1}(B)\sset E$.  $N$ is a
neighborhood of $E_{T_0}-\Dhat$ in $\widehat{P\forms}-\Dhat$, and we may
suppose without loss of generality that $\c$ and $a_1,\dots,a_{11}$
lie in $N-E$.  Now, $N-E$ is a punctured-disk bundle over
$\pi^{-1}(B)$, which in turn is a punctured $P^{12}$-bundle over $B$.
Since $B$ is a ball, $\pi^{-1}(B)\to B$ trivializes topologically, so
up to homotopy, $N-E$ is a circle-bundle over $E_{T_0}-\widehat\D$.
Now, $\pi_1(E_{T_0}-\widehat\D)$ is the 12-strand {\it spherical\/}
braid group $B_{12}(S^2)$, so $\pi_1(N-E)$ is a central extension of
$B_{12}(S^2)$ by $\Z=\langle\c\rangle$.  Furthermore, the local
description of the discriminant shows that the generators
$a_1,\dots,a_{11}$ map to the corresponding standard generators for
$B_{12}(S^2)$.  Now, $w=a_1\cdots a_{10}a_{11}^2a_{10}\cdots a_1\in
B_{12}$ represents the braid in which the leftmost strand moves in a
large circle around all the other strands.  Since this is trivial in
$B_{12}(S^2)$, $w$ is homotopic in $N-E$ to a member of the central
$\Z$, i.e., to a power of $\c$.  One can write out the $r_i\in
\Lambda(V)$ explicitly, as in \cite[sec.~5]{allcock-inventiones}, and
then matrix multiplication shows that $w$ acts on $\ch(V)$ with
order~6.  Since $P\rho(\c)$ has order dividing 6, and some power of it
has order 6, $P\rho(\c)$ itself has order 6.
\end{proof}

Now we will extend the domain of $g$, in two steps.  
We will begin with the map $g:P\framed_0\to\ch^{10}$ obtained from
\eqref{eq-per-map-formulated-with-framings}, where  $P\framed_0$ is the quotient of $\framed_0$ by
the action of $\C^*\sset\GL(5,\C)$ given in
\eqref{eq-GL5C-action-on-framed}.  
We enlarge $P\framed_0$ to a space $P\framed_s$,
which is the branched cover of $\widehat{P\forms}_s$ associated to the
covering space $P\framed_0\to P\forms_0$.  Formally, we define
$p:P\framed_s\to \widehat{P\forms}_s$ to be the Fox completion of the
composition $P\framed_0\to P\forms_0\to\widehat{P\forms}_s$.  That is,
a point of $P\framed_s$ lying over a point $T$ of
$\widehat{P\forms}_s$ is a function $\a$ which assigns to each
neighborhood $W$ of $T$ a connected component $\a(W)$ of $p^{-1}(W\cap
P\forms_0)$, in such a way that if $W'\sset W$ then
$\a(W')\sset\a(W)$.  $P\framed_s$ has a natural topology; for details
see \cite{fox}.  By the naturality of the Fox completion, the actions
of $P\Gamma$ and $PG$ extend to $P\framed_s$.

Since $P\framed_s\to\widehat{P\forms}_s$ is branched over $\Dhat\cup
E$, 
it is clear that the local structure of
$\Dhat$ and $E$ plays a key role in the nature of $P\framed_s$; by
studying it we will show that $P\framed_s$ is a complex manifold.
The analysis follows (3.3)--(3.10) of \cite{ACT}, but is more
complicated.

We first need to assemble some known results about certain
complex reflection groups.  Coxeter
\cite{coxeter-braid-group-quotients} noticed that for $n=1,\dots,$
$4$, if one adjoins to the $(n+1)$-strand braid group the relations
that the $n$ standard generators have order~$3$, then one obtains a
finite complex reflection group.  We call this group $R_n$.  One can
describe the group concretely by choosing vectors
$r_1,\dots,r_n$ that span an $n$-dimensional Euclidean complex vector
space $V_n$, such that the $i$th generator acts as $\w$-reflection in
$r_i$.  One may scale the roots in any convenient manner; we take
$r_i^2=3$ and refer to them as roots.  Then the braid and commutation
relations imply that 
$\bigl|\langle r_i|r_{i\pm1}\rangle\bigr|=\sqrt3$ and all other inner
products vanish.  By multiplying $r_2,\dots,r_n$ in turn by scalars,
we can take $r_i\cdot r_{i+1}=\theta$ for all $i$.  The group
generated by the reflections in $r_1,\dots,r_n$ is what we call $R_n$.  In each case,
$r_1,\dots,r_n$ generate an $\E$-lattice, and it turns out that the
reflections in $R_n$ are exactly the triflections in the norm~$3$
vectors of this lattice.   We write $\H_n$ for the union of the
orthogonal complements of all these vectors.  

\begin{theorem}
\label{thm-finite-reflection-group-facts}
For any $n=1,\dots,4$, the pair $\bigl(V_n/R_n,\H_n/R_n\bigr)$ is
diffeomorphic to $(\C^n,\D_{A_n})$, where $\D_{A_n}$ is the standard
$A_n$ discriminant.  $R_n$ acts freely on $V_n-\H_n$, so
$V_n-\H_n\to\C^n-\D_{A_n}$ is a covering map.  The subgroup of
$B_{n+1}=\pi_1(\C^n-\D_{A_n})$ corresponding to this covering space is
the kernel of the homomorphism $B_{n+1}\to R_n$ described above.
Finally, $V_n\to\C^n$ is the Fox completion of the composition
$$
V_n-\H_n\to\C^n-\D_{A_n}\to\C^n\;.
$$
\end{theorem}

\begin{proof}
That $V_n/R_n\isomorphism\C^n$ is the same as the ring of
$R_n$-invariants on $V_n$ being a polynomial ring, which it is by work
of Shephard and Todd \cite{shephard-todd}.  That $\H_n/R_n$
corresponds to the $A_n$ discriminant is part of the main result of
Orlik and Solomon \cite[cor.~2.26]{orlik-solomon}.  It is known that
any finite complex reflection group acts freely on the complement of
the mirrors of its reflections.
The subgroup $H$ of $B_{n+1}$ corresponding to the
covering space contains the cubes of the meridians of $\D_{A_n}$,
since $R_n$ contains the triflections across the components of
$\H_n$.  Since modding out $B_{n+1}$ by the cubes of the standard
generators yields a copy of $R_n$, the cubes
of meridians generate $H$, and $B_{n+1}/H\isomorphism  R_n$ under the indicated
homomorphism.  

The claim about the Fox completion is a special case of the following:
suppose $G$ is a finite group acting linearly and faithfully on a
finite-dimensional real vector space $V$, and contains no (real)
reflections.  Then, writing $V_0$ for the open subset of $V$ on which
$G$ acts freely, $V\to V/G$ is the Fox completion of $V_0\to
V_0/G\to V/G$.  (One just verifies that $V\to V/G$ satisfies the
definition of a completion of $V_0\to V/G$.  The absence of real
reflections in $G$ is required for $V_0$ to be locally connected in
$V$, in Fox's terminology.)
\end{proof}

Now we can describe the Fox completion $P\framed_s$.  First we
describe it away from the chordal locus, and then at a point in the
chordal locus.

\begin{theorem}
\label{thm-fox-completion-away-from-E}
Suppose $T\in\widehat{P\forms}_s-E$ has $n_i$ singularities of type
$A_i$, for each $i=1,\dots,4$.  Suppose $\breve T\in P\framed_s$
lies over $T$.  Then near $\breve T$, $P\framed_s$ has a complex
manifold structure, indeed a unique one for which $P\framed_s\to
P\forms_s$ is holomorphic.  With respect to this structure, $\breve T$
has a neighborhood in $P\framed_s$ diffeomorphic to
$$
(B^1)^{n_1}\times(B^2)^{n_2}\times(B^3)^{n_3}\times(B^4)^{n_4}\times B^N\;,
$$
where $N=34-n_1-2n_2-3n_3-4n_4$, such that $P\framed_0$ corresponds to 
$$
(B^1-\H_1)^{n_1}\times(B^2-\H_2)^{n_2}\times(B^3-\H_3)^{n_3}\times(B^4-\H_4)^{n_4}\times B^N\;.
$$
The stabilizer of $\breve T$ in $P\Gamma$ is isomorphic to
$G_1^{n_1}\times\dots\times G_4^{n_4}$, acting in the obvious way, and
the map to $\widehat{P\forms}_s$ is the quotient by this group action.
\end{theorem}

\begin{proof}
By lemma~\ref{lem-discriminant-away-from-E}, $T$ has a neighborhood $U\sset P\forms$ diffeomorphic to
\begin{equation}
\label{eq-local-model-of-Fox-completion}
(B^1/G_1)^{n_1}\times\dots\times(B^4/G_4)^{n_4}\times B^N\;,
\end{equation}
such that $U\cap P\forms_0$ corresponds to (by theorem~\ref{thm-finite-reflection-group-facts})
\begin{equation}
\label{eq-local-model-of-hyperplane-complement-in-Fox-completion}
\bigl((B^1-\H_1)/G_1\bigr)^{n_1}\times\dots\times\bigl((B^4-\H_4)/G_4\bigr)^{n_4}\times B^N\;,
\end{equation}
and the local fundamental group is $B_2^{n_1}\times\dots\times
B_5^{n_4}$.  We write $T'$ for a basepoint in $U-\Dhat$, so that we can
refer to its associated fourfold $V'$.   By
lemma~\ref{lem-meridians-act-by-complex-reflections}, any standard generator of any of the braid group factors
acts on $\ch(V')$ as the $\w$-reflection in a root $r\in \Lambda(V')$.
We write $H\sset\Gamma(V')$ for the group generated by all these
reflections.  The local 
monodromy group is by definition the projectivization of $H$.

By lemma~\ref{lem-local-monodromy-noncyclic}, distinct generators of the local fundamental group
give linearly independent roots.  Therefore
the discussion before theorem~\ref{thm-finite-reflection-group-facts} shows
that $H$ is $R_1^{n_1}\times\dots\times R_4^{n_4}$.  Since the
$\E$-sublattice spanned by the roots is positive-definite, it
has lower dimension than $\Lambda(V')$, so $H$ contains no scalars.
Therefore $P\rho\bigl(\pi_1(U-\Dhat)\bigr)$ is a copy of $H$.  By
theorem~\ref{thm-finite-reflection-group-facts}, the covering space of $U-\Dhat$ associated to the
kernel of this monodromy is 
$$
(B^1-\H_1)^{n_1}\times\dots\times(G^4-\H_4)^{n_4}\times B^N\;,
$$ 
with the deck group being $H$, acting in the obvious way.
Furthermore, the Fox completion  over $U$ is
then
$$
(B^1)^{n_1}\times\dots\times(B^4)^{n_4}\times B^N\;,
$$
with $\breve T$ being the point at the center.  

Since \eqref{eq-local-model-of-Fox-completion} is a diffeomorphism, not just a homeomorphism,
$(B^1)^{n_1}\times\dots\times(B^4)^{n_4}\times B^N\to U$ is complex
analytic when the domain is equipped with the standard complex
manifold structure.  This gives the Fox completion a complex manifold
structure such that $P\framed_s\to P\forms_s$ is holomorphic.  A
standard argument using Riemann extension shows that this structure is
unique.
\end{proof}

\begin{theorem}
\label{thm-fox-completion-over-E}
Suppose $(T,\tau)\in E\cap\widehat{P\forms}_s$, where $\tau$ has $n_i$
points of multiplicity $i+1$, for each $i=1,\dots,4$.  Suppose
$(T,\tau)\spbreve\in P\framed_s$ lies over $(T,\tau)$.  Then near
$(T,\tau)\spbreve$, $P\framed_s$ has a complex manifold structure,
indeed a unique one for which $P\framed_s\to\widehat{P\forms}_s$ is
holomorphic.  With respect to this structure, $(T,\tau)\spbreve$ has a
neighborhood diffeomorphic to
$$
B^1\times(B^1)^{n_1}\times\dots\times(B^4)^{n_4}\times B^{N-1}\;,
$$ 
where $N$ is as in theorem~\ref{thm-fox-completion-away-from-E}, such that the preimage of $E$
corresponds to 
\begin{equation}
\label{eq-foo-6}
\{0\}\times(B^1)^{n_1}\times\dots\times(B^4)^{n_4}\times B^{N-1}
\end{equation}
and the preimage of $P\forms_0$ corresponds to
$$
\bigl(B^1-\{0\}\bigr)\times\bigl(B^1-\H_1\bigr)^{n_1}\times\dots\times\bigl(B^4-\H_4\bigr)^{n_4}\times B^{N-1}\;.
$$ The stabilizer of $(T,\tau)\spbreve$ in $P\Gamma$ is isomorphic to
$\Z/6\times G_1^{n_1}\times\dots\times G_4^{n_4}$, with the $G_i$'s
acting in the obvious way.  The $\Z/6$ acts freely away from
\eqref{eq-foo-6}.
\end{theorem} 

\begin{proof}
This is much the same as the previous proof.  The difference is
that we don't have a local analytic description of $\Dhat$ near $E$,
only weaker results, theorems~\ref{thm-discr-along-E-topological-model} and~\ref{thm-discr-along-E-analytic-model}.  We begin with the
local monodromy analysis.  Theorem~\ref{thm-discr-along-E-topological-model} provides a neighborhood $U$
of $(T,\tau)$ with 
$$
\pi_1\bigl(U-(E\cup\Dhat)\bigr)\isomorphism\Z\times
(B_2)^{n_1}\times\dots\times (B_5)^{n_4}\;,
$$ 
where a generator for the $\Z$ factor is a meridian $\c$ of $E$,
and the standard generators for the braid group factors are meridians
of $\Dhat$.  As in the previous proof, we write $T'$ for a basepoint
in $U-(E\cup\Dhat)$, so we can refer to the associated fourfold $V'$.
We write $H$ for the subgroup of $\Gamma(V')$ generated by the
reflections in the roots associated to the braid group factors.  By
lemma~\ref{lem-monodromy-around-E}, $P\rho(\c)$ is a hexaflection of $\ch(V')$, which is the
projectivization of a hexaflection $S$ of $\Lambda(V')$ in a chordal
root $s$ of $\Lambda(V')$.  We write $H'$ for $\langle H,S\rangle$.
The local monodromy group
$P\rho\bigl(\pi_1\bigl(U-(E\cup\Dhat)\bigr)\bigr)$ is the
projectivization of $H'$.  Following the previous proof shows that
$H\isomorphism G_1^{n_1}\times\dots\times G_4^{n_4}$.  We claim that $s$
is orthogonal to all the roots of the braid group factors.  To prove
this, we use the fact that $S$ commutes with $H$, so that for every
nodal root $r$ of a braid group factor, the triflection $R$ in $r$
carries $s$ to a multiple of itself.  Therefore, either $s$ is
orthogonal to all the $r$'s, or else it is proportional to one of
them.  The latter is impossible because then $s$ would be both nodal
and chordal, which is impossible.  Since $s$ is orthogonal to the
$r$'s, $H'=\Z/6\times H$.  Arguing as in the previous proof, $H'$
contains no scalars, so it maps isomorphically to its
projectivization.  Continuing as before proves the corollary, with
``diffeomorphic'' replaced by ``homeomorphic''.

To prove the existence of the complex manifold structure, we proceed
in two steps.  First, we take $\tilde U$ to be the 6-fold cover of
$U$, branched over $U\cap E$.  This clearly has a complex manifold
structure such that $\pi_{\tilde U}:\tilde U\to U$ is holomorphic.
Writing $(T,\tau)\sptilde$  for the
preimage of $(T,\tau)$, theorem~\ref{thm-discr-along-E-analytic-model} provides us with a
neighborhood $\tilde V$ 
of $(T,\tau)\sptilde$ with the properties stated
there.  The important property is that $\tilde V$ is {\it
  diffeomorphic\/} to
$B^1\times(B^1/R_1)^{n_1}\times\dots\times(B^4/R_4)^{n_4}\times B^{N-1}$, such
that $\tilde V-\tilde\D$ corresponds to
$$
B^1\times\bigl((B^1-\H_1)/R_1\bigr)\times\dots\times
\bigl((B^4-\H_4)/R_4\bigr)\times B^{N-1}\;.
$$
Taking the branched
cover of $\tilde V$ with deck group $H$ gives the claimed complex manifold
model of $P\framed_s$ near $(T,\tau)\spbreve$.  The uniqueness of the complex
manifold structure 
again follows from Riemann extension.
\end{proof}

Extending the period map to $P\framed_s$ is now easy.  If $r$ is a
root of $\Lambda$, then we will call $r^\perp\sset\ch^{10}$ a
discriminant hyperplane or chordal hyperplane according to whether $r$
is a nodal or chordal root.  By the chordal (resp. discriminant) locus
of $P\framed_s$, we mean the preimage of $E\sset\widehat{P\forms}_s$
(resp. $\Dhat$).

\begin{theorem}
\label{thm-period-map-extends-to-stable}
The period map $P\framed_0\to\ch^{10}$ extends to
a holomorphic  map
$g:P\framed_s\to\ch^{10}$, which is invariant under $PG$ and
equivariant under
$P\Gamma$.  
The chordal (resp. discriminant) locus of $P\framed_s$ maps into the
chordal (resp. discriminant) hyperplanes of $\ch^{10}$.
\end{theorem}

\begin{proof}
$P\framed_s$ is a complex manifold by
  theorems~\ref{thm-fox-completion-away-from-E}
  and~\ref{thm-fox-completion-over-E}; 
  since $g$ is a map to a bounded domain, the extension exists by
  Riemann extension.  The
preimages of $\widehat\D$ and $E$ map into hyperplanes as claimed
because of $P\Gamma$-equivariance.  Namely, a generic point of
$P\framed_s$ lying over $\Dhat$ has stabilizer $\Z/3$ in $P\Gamma$, so it
maps to the fixed-point set of $\Z/3$ in $\ch^{10}$, which is a discriminant
hyperplane.  The same idea applies when $\Dhat$ is replaced by $E$.
\end{proof}

Since the period map of theorem~\ref{thm-period-map-extends-to-stable} is $\PGamma$-equivariant, it
induces a map
$$
\widehat{P\forms}_s=P\Gamma\backslash P\framed_s\to
P\Gamma\backslash\ch^{10}\;,
$$
which we will now extend further.  As before, we continue to use the
notation $g$.  The argument relies on a monodromy analysis near a
threefold with an $A_5$ or $D_4$ singularity; we give the key point as
a lemma:

\begin{lemma}
\label{lem-monodromy-for-A5-and-D4}
Suppose $F\in\forms_{ss}$ defines a threefold $T$ with a singularity
of type $A_5$, and let $U$ be a neighborhood of $F$ as in
lemma~\ref{lem-discriminant-away-from-E}, with basepoint $F'$.  Let $a_1,\dots,a_5$ be standard generators for the
corresponding factor $B_6$ of $\pi_1(U-\D,F')$, and let $r_1,\dots,r_5$
be roots of $\Lambda(V')$, by whose $\w$-reflections the $a_i$ act.
Suppose $\ip{r_i}{r_{i\pm1}}=\pm\theta$ and all other inner products
are zero.  Then $\xi=r_1-\theta r_2-2r_3+\theta r_4+r_5$ is a nonzero
isotropic vector of $\Lambda(V')$, and $(a_1\dots a_5)^6$ acts on
$\Lambda(V')$ by the unitary transvection in $\xi$, namely
$$
x\mapsto x-\frac{\ip{x}{\xi}}{\theta}\xi\;.
$$
Now suppose the singularity has type $D_4$ rather than $A_5$, and that
$a_1$, $a_2$, $a_3$ and $b$ are standard generators for
$\A(D_4)\sset\pi_1(U-\D)$, with $b$ corresponding to the central node
of the $D_4$ diagram.  Suppose $r_1$, $r_2$, $r_3$ and $r'$ are roots
for the corresponding $\w$-reflections of $\Lambda(V')$, scaled so
that $\ip{r_i}{r'}=\theta$ for $i=1,2,3$.  Then
$\xi=r_1+r_2+r_3-\theta r'$ is a nonzero isotropic vector of
$\Lambda(V')$, and $(a_1a_2a_3b)^3$ acts on $\Lambda(V')$ by the
unitary transvection in $\xi$.
\end{lemma}

\begin{remark}
The given words in the Artin generators generate the centers of the
Artin groups, and in each case, $\xi$ spans the kernel of the
restriction of $\ip{\,}{}$ to the span of the roots.  So it isn't
surprising that the word acts by a transvection in a multiple of
$\xi$.  The point of the lemma is that this multiple is nonzero.
\end{remark}

\begin{proof}
We treat the $A_5$ case first.  
We remark that the $r_i$ are pairwise linearly independent by
lemma~\ref{lem-local-monodromy-noncyclic}, and by scaling them we may assume that their inner
products are as stated.  We need the sharper result that
$r_1,\dots,r_5$ are linearly independent.
Direct calculation using the given
inner products shows that $\xi$ is
isotropic and orthogonal to $r_1,\dots,r_5$.  One can show that the
locus of cubic threefolds with an $A_5$ singularity is irreducible, so
to prove $\xi\neq0$, it suffices to treat a single example.  If $T'$
has an $A_6$ singularity, then one can write down the inner product
matrix for its six roots and check that it is nondegenerate.
Therefore its six roots are linearly independent.  In particular, the
first five are, so $\xi\neq0$.  The fact that $(a_1\cdots a_5)^6$ acts
as the transvection in $\xi$ is a matrix calculation---one writes down
any linearly independent set of roots in $\C^{10,1}$ with these inner
products, and just multiplies the reflections together suitably.

The same idea works for the $D_4$ case.  To show $\xi\neq0$, one must
check (i) the irreducibility of the the set of threefolds with a $D_4$
singularity, (ii) that there is a threefold with a $D_5$ singularity
admitting a versal deformation in $P\forms$, and (iii) the inner
product matrix for 5 roots corresponding to the generators of
$\A(D_5)$ is nondegenerate.  Then one does a matrix calculation to
check the action of $(a_1a_2a_3b)^3$.
\end{proof}

\begin{theorem}
\label{thm-period-map-extends-to-semistable}
The period map $\widehat{P\forms}_s\to P\Gamma\backslash\ch^{10}$ extends
to a holomorphic map
$g:\widehat{P\forms}_{ss}\to\overline{P\Gamma\backslash\ch^{10}}$.  This
map sends $\widehat{P\forms}_{ss}-\widehat{P\forms}_s$ to the boundary
points of the Baily-Borel compactification.  
\end{theorem}

\begin{proof}
The extension of $g$ to $\widehat{P\forms}_{ss}-E$ follows section~8
of \cite{ACT}, but is  a little more complicated.  After explaining this,
we will extend $g$ to $\widehat{P\forms}_{ss}\cap E$ by modifying the
argument, in the same way that we modified the proof of theorem~\ref{thm-fox-completion-away-from-E}
to prove theorem~\ref{thm-fox-completion-over-E}.

We begin by supposing
$T\in\widehat{P\forms}_{ss}-(\widehat{P\forms}_s\cup E)$, with
defining function $F$.  We adopt the notation of lemma~\ref{lem-discriminant-away-from-E}, so that
$T$ has $m$ singularities $s_1,\dots,s_m$, and $U$ is a neighborhood
of $F$ in $\forms_{ss}$ with the properties stated there.  Let $\breve
U$ be the universal cover of $U$ with two-fold branching over
$U\cap\D$.  By Brieskorn's description \cite{brieskorn-versals} of versal deformations
of simple singularities, $\breve U$ is diffeomorphic to a neighborhood
of the origin in $\C^{34}$, such that the preimage $\breve\D$ of $\D$
is the union of the reflection hyperplanes for the Coxeter group
$W_1\times\dots\times W_m$, where $W_i$ is the Coxeter group of the
same type as the singularity $s_i$.  Let $\b$ be a loop lying in a
generic line through the origin in $\C^{35}$, and encircling the
origin once positively.  

We claim that the sixth power of $\rho(\b)$ is nontrivial and unipotent.
(Discussing $\rho(\b)$ requires choosing a basepoint in $\breve
U-\breve\D$ and a framing of the threefold represented by the
corresponding threefold $T'\in U-\D$.  These choices are immaterial.)  

If $T$ has only one singularity, of type $A_n$, then in terms of the
standard generators $a_1,\dots,a_n$ for $\pi_1(U-\D)\isomorphism
B_{n+1}$, $\b=(a_1\dots a_n)^{n+1}$.  The $a_i$ act by $\w$-reflections
in linearly independent roots $r_1,\dots,r_n$ of $\Lambda(V)$, with
$r_i\perp r_j$ except 
for $\ip{r_i}{r_{i\pm1}}=\pm\theta$.   
This lets one work out the
action of $\b$, by choosing such roots and multiplying matrices
together.  (Two choices of such roots are equivalent under $\U(9,1)$,
so the conjugacy class of $\rho(\b^6)$ in $\U(9,1)$ is independent of
the choice.)  Direct calculation shows that $\rho(\b)$ has order $3$, $2$,
$3$ or $6$ if $n=1,\dots,4$, and  if $n=5$ then it is a nontrivial
unipotent by lemma~\ref{lem-monodromy-for-A5-and-D4}. 
(No calculation is required to see that $\rho(\b)$ has finite order if
$n<5$, because $\rho(B_{n+1})$ is the finite group $R_n$.)  If $T$ has
only one singularity, of type $D_4$, with standard generators $a_1$,
$a_2$, $a_3$ and $b$, $b$ corresponding to the central node, then
$\b=(a_1a_2a_3b)^6$, again lemma~\ref{lem-monodromy-for-A5-and-D4} shows that $\rho(\b)$
is a nontrivial unipotent.

If $T$ has $m$ singularities, then $\b=\b_1\dots\b_m$, where each
$\b_i$ is as $\b$ above, one for each singularity.  Since the $\b_i$'s
commute, $\rho(\b^6)$ is a product of nontrivial commuting unipotent
isometries, one for each $A_5$ or $D_4$ singularity.  Since $T$ is not
stable, there is at least one $A_5$ or $D_4$ singularity.  If there is only one,
then this proves that $\rho(\b^6)$ is nontrivial and unipotent.  

If there are more than one, then the product is unipotent since it is
a product of commuting unipotents.  A little extra work is required to
show that it is nontrivial, i.e., that no cancellation occurs.
Suppose $s$ and $s'$ are two singularities of $T$, each of type $A_5$
or $D_4$.  ($T$ cannot have both an $A_5$ and a $D_4$ singularity, but
this isn't needed here.)  Let $\xi$ and $\xi'$ be the isotropic
vectors from lemma~\ref{lem-monodromy-for-A5-and-D4}.  They are orthogonal because they
correspond to distinct singularities.  Since $\Lambda(V')$ has
signature $(9,1)$, $\xi$ and $\xi'$ are proportional.  By using the
formula defining a unitary transvection, one can check that the
product of the transvections in $\xi$ and $\xi'$ is a transvection in
a nonzero multiple of $\xi$.  The same argument applies when there are
more than two singularities of type $A_5$ or $D_4$.  (There is only
one $PG$-orbit of such threefolds, which has three $D_4$
singularities.) 

Because a power of $\rho(\b)$ is nontrivial and unipotent, it
fixes a unique point of $\partial\ch^{10}$.  By unipotence, this fixed
point 
is represented by a null vector of
$\Lambda$, so there is an associated boundary point $\eta$ of
$\overline{\PGamma\backslash\ch^{10}}$.  We claim that for every
neighborhood $Z$ of $\eta$, there is a neighborhood $\breve Y$ of the
origin in $\breve U$, such that the composition
$$
\breve U-\breve\D\quad\to\quad
U-\D\quad\to\quad
\PGamma\backslash\ch^{10}\quad\to\quad
\overline{\PGamma\backslash\ch^{10}}
$$
carries $\breve Y-\breve\D$ into $Z$.  This uses a hyperbolic analysis
argument (see p.~708 of \cite{ACT}).  It follows that for every
neighborhood $Z$ of $\eta$, there is a neighborhood $Y$ of $F$ in
$\forms$, such that $g(Y-\D)\sset Z$.  Then, by Riemann extension,
$g:Y-\D\to\overline{P\Gamma\backslash\ch^{10}}$ extends
holomorphically to $Y$, carrying $F$ to $\eta$.

Now suppose $(T,\tau)\in\widehat{P\forms}_{ss}\cap E$, with $\tau$
having a point of multiplicity six, and let $\tilde U$ and $\tilde V$
be as in theorem~\ref{thm-discr-along-E-analytic-model}.  Then let $\breve V$ be the universal cover
of $\tilde V$ with 2-fold branching over $\tilde\D\sset\tilde V$.  We
write $\breve\D$ and $\breve E$ for the preimages of $\D$ and $E$ in
$\breve V$.  From here, the treatment is exactly the same as above, except
that $\breve E$ is present.  This affects nothing, because the
monodromy around $\breve E$ is trivial.  (The
monodromy around $E$ has order~$6$, and $\tilde U$ is the 6-fold cover
of $U$ branched over $U\cap E$.)
\end{proof}

\begin{theorem}
\label{thm-exactly-two-cusps}
There are exactly two
boundary points of $\overline{P\Gamma\backslash\ch^{10}}$.  The
threefolds having a $D_4$ singularity map to one, and those having an $A_5$ singularity map to the other.
The points $(T,\tau)$ of $E$, where $\tau$ has a point of
multiplicity~$6$, map to the latter boundary point.
\end{theorem}

\noindent
We will call the boundary points the $A_5$ and $D_4$ cusps of
$\overline{\PGamma\backslash\ch^{10}}$. 

\begin{proof}
We know that the preimage in $\widehat{P\forms}_{ss}$ of the boundary
consists of $\widehat{P\forms}_{ss}-\widehat{P\forms}_s$.  The GIT
equivalence classes in $\widehat{P\forms}_{ss}-\widehat{P\forms}_s$
are represented by the points \eqref{item-closed-orbits-of-semistable-threefolds} and \eqref{item-closed-orbits-of-semistable-12-tuples} of
theorem~\ref{thm-GIT-analysis}.  The threefolds $T_{A,B}$ form a 1-parameter family,
limiting to \eqref{item-closed-orbits-of-semistable-12-tuples}.  Therefore they all map to a single boundary
point.  We call the union of these GIT equivalence classes the
$A_5$ component of $\widehat{P\forms}_{ss}-\widehat{P\forms}_s$.  Only
one GIT equivalence class remains, which we call the $D_4$ component
of $\widehat{P\forms}_{ss}-\widehat{P\forms}_s$.  Therefore there are
at most two boundary points.

Also, the images of the $A_5$ and $D_4$ components in
$\widehat{P\forms}_{ss}//\SL(5,\C)$ are disjoint closed sets, so the
$A_5$ and $D_4$ components
have disjoint $PG$-invariant neighborhoods in
$\widehat{P\forms}_{ss}$.   If both
components mapped to a single boundary point, then the period map
could not be injective on 
$\moduli_0$, which it is by theorem~\ref{thm-main-theorem-smooth-case}.
\end{proof}

\section{Degeneration to a chordal cubic}
\label{sec-hodge-theory-chordal}

The aim of this section is to identify the limit Hodge structure for the
degeneration of cyclic quartic fourfolds 
associated to a generic degeneration to cubic threefolds to a chordal
cubic.  The following theorem is this section's contribution to the
proof of the main theorem of the paper, theorem~\ref{thm-MAIN-THEOREM}.

\begin{theorem}
\label{thm-chordal-maps-onto-divisor}
The period map $g:\widehat{P\forms}_{ss}\to \PGamma\backslash\ch^{10}$
carries the chordal locus onto a divisor.
\end{theorem}

This immediately implies that the chordal cubics are points of
indeterminacy for the rational map
$\forms_{ss}\dasharrow\overline{\PGamma\backslash\ch^{10}}$ obtained from the
period map $\forms_s\to\PGamma\backslash\ch^{10}$.

The theorem is a consequence of
theorem~\ref{thm-main-thm-chordal-limits} below, which most of the
section is devoted to proving.  This theorem describes the limit Hodge
structure in terms of a Hodge structure studied by Deligne and Mostow
\cite{Mostow},\cite{Deligne-Mostow}, associated to a
12-tuple of points in $P^1$.

\subsection{Statement of Results}

Let us establish notation and definitions.  The chordal cubic $T$ is
the secant variety of a rational normal curve $R$ in $P^4$.  $R$ is
the whole
singular locus of $T$.
Consider a pencil $\{T_t\}$ of cubic threefolds with $T_0=T$.  In
homogeneous coordinates, the family near $T_0$ is
$$
\set{\bigl([x_0{:}\dots{:}x_4],t\bigr)\in P^4\times\D}{F(x_0,\dots,x_4)+tG(x_0,\dots,x_4)=0}\;,
$$
where $F$ defines $T=T_0$, $G$ is some other member of the pencil, and
$\D$ is the unit disk in $\C$.  We assume that the pencil is generic
in the sense that $T_t$ is smooth for all sufficiently small nonzero
$t$.  By scaling $G$ we may suppose that $T_t$ is smooth for all
$t\in\D-\{0\}$.  We also make the genericity assumption that $G=0$
cuts out on $R$ a set $B$ of 12 distinct points, the ``infinitesimal
base locus''. 

Taking the threefold covers of $P^4$ over the $T_t$ gives a family
$\V$ of fourfolds, which is the restriction to $\D\sset\forms$ of the
family called $\V$ in the rest of the paper.  Explicitly,
\begin{align*}
\V=\bigl\{&\bigl([x_0{:}\dots{:}x_5],t\bigr)\in P^5\times\D\bigm|\\
&F(x_0,\dots,x_4)+tG(x_0,\dots,x_4)+x_5^3=0\bigr\}\;.
\end{align*}
All fibers of $\V$ are smooth except for $V_0$.

By lemma~\ref{lem-monodromy-around-E}, the monodromy of $\V|_{\D-\{0\}}$ on $H^4(V_t)$ has
finite index, so there is a well-defined limit Hodge structure.
To discuss this limiting Hodge structure, we define the
notion of an ``Eisenstein Hodge structure.''  This is an Eisenstein
module $H$ with a complex Hodge structure on the vector space
$H\otimes_{\E}\C$, where multiplication by $\omega$ preserves the
Hodge decomposition.  A morphism of such objects is a homomorphism of
Eisenstein modules which is compatible with the Hodge decomposition.

The Eisenstein Hodge structure of interest to us throughout the paper
is $\Lambda(V)$, for $V$ a smooth cyclic cubic 4-fold.  Recall that its
underlying abelian group is $H^4_0(V;\Z)$, and that the $\E$-module
structure is defined by taking $\w$ to act as $\sigma^*$, where
$\sigma$ is from \eqref{eq-def-of-sigma}.  Then $H^4_0(V;\R)$ is
identified with $\Lambda(V)\tensor_\E\C$ and is hence a complex vector
space.  The projection to $\sigma$'s $\w$-eigenspace identifies it
with $H^4_{\sigma=\w}(V;\C)$.  In light of this, the Hodge
decomposition 
\begin{equation}
\label{eq-Eisenstein-hodge-struc-on-V}
H^4_{\sigma=\w}(V;\C)=H^{3,1}_{\sigma=\w}(V)\oplus
H^{2,2}_{\sigma=\w}(V)
\end{equation}
gives an Eisenstein Hodge structure on
$\Lambda(V)$.

The limit mixed Hodge structure of a degeneration can be computed from
the mixed Hodge structure of the central fiber provided that the
degeneration is \emph{semistable}.  One computes the limit using the
Clemens-Schmid sequence \cite{Clemens-Schmid}.  A semistable model
$\hat\V$ for a degeneration $\V$ is a family over $\D$ disk with the
following properties: 

(a) There is a surjective morphism $\hat \V \to \V$ which makes the diagram
\[
\begin{CD}
\hat \V @>>> \V \\
@VVV @VVV \\
\Delta @>>> \Delta
\end{CD}
\]
commute, where the bottom arrow is $t= s^n$ for some $n$, with $s$ 
a parameter on the left $\D$ and $t$ a parameter on the right, and the
restriction of $\hat\V$ to $\D-\{0\}$ is the pullback of
$\V|_{\D-\{0\}}$.  We write $V$ and $\hat V$ for the central fibers of
$\V$ and $\hat\V$, and $V_t$ and $\hat V_s$ for other fibers.

(b) The total space of $\hat \V$ is smooth and the central fiber $s=0$
is a normal crossing divisor (NCD) with smooth components, each of
multiplicity one.

It is known 
that every one-parameter degeneration has
a semistable model, obtainable by artful combination of three moves:
base extension, blowing up, and normalization.  In our case we will
apply one base extension, replacing $t$ by $s^6$, followed by three
blowups along various ideal sheaves.

In our case the Clemens-Schmid sequence reads as follows:,
\[
   \cdots \to 
      H^4(\hat\V, \hat\V^*) 
   \to
      H^4(\hat V)
   \to \lim_{s \to 0} 
      H^4(\hat V_s) 
   \mapright{N}  
      H^4(\hat V_s)
      \to\cdots
\] 
where $N$ is the logarithm of the monodromy transformation.  The
asterisk indicates the restriction of $\hat\V$ to the punctured disk.
The terms of the sequence are abelian groups equipped with Hodge structures,
but in our calculations we will only need to work with rational coefficients.
Since the monodromy of $\V^*$ on $H^4$ of the fiber has order~$6$
(by lemma~\ref{lem-monodromy-around-E}), the 
base extension ensures that the
monodromy of the semistable model is trivial.  Therefore $N = 0$ and
the sequence reduces to
\begin{equation}
\label{eq-Clemens-Schmid-sequence}
   \cdots \to 
      H^4(\hat\V,\hat\V^*) 
   \to
      H^4(\hat V)
   \to \lim_{s \to 0} 
      H^4(\hat V_s) 
   \to 
     0.
\end{equation}
Since the monodromy is trivial, the limit is a (pure) Hodge structure
and is presented as a quotient of the mixed Hodge structure on
$H^4(\hat V)$.  

Because we are interested in the ``Eisenstein part'' of the sequence,
we make the following definition.
Let $H$ be a $\Z$-module on which an automorphism $\zeta$ of order
$n$ acts.  Let $p_k$ be the cyclotomic polynomial of degree $k$, where $k$
is a divisor of $n$.  Let $H_{(k)}$ be the kernel of $p_k(\zeta)$
acting on $H$.  Then there is a decomposition of $H\otimes\Q$ into a
direct sum of subspaces $H_{(k)}\otimes\Q$.  Moreover, $H$ is
commensurable with the direct sum of the $H_{(k)}$.  Since the
characteristic polynomial for the action of $\zeta$ is $p_k$ on
$H_{(k)}$, we call $H_{(k)}$ the $k$th characteristic submodule (or subspace)
of $H$, or alternatively, the \emph{$k$-part} of $H$.  

Passing to the 3-parts of the terms of
\eqref{eq-Clemens-Schmid-sequence}, we have the sequence
\begin{equation}
\label{eq-3-part-of-Clemens-Schmid}
   \cdots \to 
      H^4_{(3)}(\hat\V,\hat\V^*) 
   \to
      H^4_{(3)}(\hat V)
   \to \lim_{s \to 0} 
      H^4_{(3)}(\hat V_s) 
   \to 
     0,
\end{equation}
which need not be exact.  However, it is  exact when one takes $\Q$
as the coefficient group.

The Eisenstein Hodge structure of $\hat V$ is the focus of the rest of
this section; we describe it in
theorem~\ref{thm-main-thm-chordal-limits} below.  The description is
in terms of a curve $C$ determined by the inclusion $B\sset R$.
Identifying $R$ with $P^1$, $B$ is the zero locus of homogeneous
polynomial $f$ of degree~12.  The curve $C$ is a certain $6$-fold
cover of $P^1$ branched over $B$, namely
\begin{equation}
\label{eq-defn-of-C}
C=\set{[x{:}y{:}z]\in P(1,1,2)}{f(x,y)+z^6=0}\;,
\end{equation}
and it has an automorphism 
\begin{equation}
\label{eq-defn-of-zeta-the-order-6-deck-transformation}
\zeta\bigl([x{:}y{:}z]\bigr)=[x{:}y{:}-\w z]
\end{equation}  
of order~$6$.  The Griffiths
residue calculus \cite{tu} 
shows that $H^{1,0}(C)$ is the vector space
consisting of the residues of the rational differentials
\begin{equation}
\label{eq-rat-differential-on-P112}
\frac{a(x,y,z)\Omega}{f(x,y)+z^6}
\end{equation}
where $\Omega=xdydz-ydxdz+zdxdy$ and the total weight is zero, so $a$
is a polynomial of weight~$8$.  Now, $\zeta$ scales $z$ and $\Omega$
by $-\w$, and it follows that $H^{1,0}_{\zeta=-\wbar}(C)$ is
1-dimensional, spanned by \eqref{eq-rat-differential-on-P112} with
$a=1$.  One can also check that $H^{0,1}_{\zeta=-\wbar}(C)$ is
$9$-dimensional.  We define an Eisenstein module structure
$\Lambda_{10}(C)$ by taking the underlying group to be
$H^1_{(6)}(C;\Z)$ and taking $-\wbar$ to act as $\zeta^*$.  By the
same considerations as above, $\Lambda_{10}(C)\tensor_\E\C$ is
identified with $H^1_{\zeta=-\wbar}(C;\C)$, so the Hodge structure of
$C$ gives an Eisenstein Hodge structure on $\Lambda_{10}(C)$.  The
notation $\Lambda_{10}(C)$ is intended to emphasize the parallel
between the isomorphism in theorem~\ref{thm-main-thm-chordal-limits}
and the description of $\Lambda(V)$ in
theorem~\ref{thm-inner-product-matrix}.

\begin{theorem}
\label{thm-main-thm-chordal-limits}
There is an isogeny of Eisenstein Hodge structures
$$
\Lambda(\hat V)\to\E \oplus \Lambda_{10}(C)\;
$$
of weight $(-2,-1)$, 
where $\E$ indicates a $1$-dimensional Eisenstein lattice with 
Hodge type $(1,0)$.
\end{theorem}

\begin{remark}
The total weight of the map is $-3$.  It is odd for a map of Hodge
structures to have odd total weight.  Such odd weight morphisms were
first considered by van Geemen in his paper on half twists \cite{vangeeman}.
\end{remark}

\begin{proof}[Proof of theorem~\ref{thm-chordal-maps-onto-divisor}, given theorem~\ref{thm-main-thm-chordal-limits}]
A consequence of theorem~\ref{thm-main-thm-chordal-limits} is that $H^4_{(3)}(\hat V)$
has rank~22 as a $\Z$-module, which is the same as the rank of the
limit Hodge structure.  Therefore \eqref{eq-3-part-of-Clemens-Schmid} shows that the Eisenstein
Hodge structures of $\Lambda(\hat V)$ and $\lim_{s\to0} \Lambda(V_s)$
are isogenous.  
The period map for the family of Hodge structures $\Lambda_{10}(C)$
has rank~9, by a standard calculation with the Griffiths residue
calculus.  
Therefore, as the pencils in $P\forms$ through $T$
vary, their limiting Hodge structures sweep out a 9-dimensional
family.   It follows from theorem~\ref{thm-period-map-extends-to-semistable} that the period map
is proper, so the image of the chordal locus is closed, hence a divisor.

We remark that  Mostow \cite{Mostow} (see also
\cite{Deligne-Mostow})  proved that the Eisenstein Hodge structure
on $H^1_{(6)}(C)$ determines the point set $B$ in $P^1$ up to
projective transformation. This is a global Torelli theorem, rather than just the
local one provided by the residue calculation referred to above.
\end{proof}

\subsection{Overview of the Calculations}

This subsection outlines the calculations required for proving
theorem~\ref{thm-main-thm-chordal-limits}.  In order to discuss $\hat V$
we need to describe briefly our particular semistable model $\hat\V$.
In the next subsection, we define $\V_0$ as the degree $6$ base
extension of $\V$, $\V_1$ as a blowup of $\V_0$, $\V_2$ as a blowup of
$\V_1$ and $\V_3$ as a blowup of $\V_2$.  $\hat\V$ is $\V_3$.  We
write $E_1\sset\V_1$, $E_2\sset\V_2$ and $E_3\sset\V_3$ for the
exceptional divisors of the blowups, and indicate proper transforms of
$V$ and the $E_i$ in subsequent blowups by adding primes.  For
example, $V'$ is the proper transform of $V$ in $\V_1$, and the
central fiber of $\V_3$ is
\[
   \hat V = V''' \cup E_1'' \cup E_2' \cup E_3.
\]

Using the Mayer-Vietoris and Leray spectral sequences, we show in
lemma~\ref{lem-central-fiber-coho-chordal} that 
only $V'''$ and $E_1''$ contribute to $H^4_{(3)}(\hat V;\Q)$, and these
contributions depend only on $V$ and $E_1$.  That is,
\begin{equation}
   H^4_{(3)}(\hat V;\Q) \cong H^4_{(3)}(V;\Q) \oplus H^4_{(3)}(E_1;\Q).
\label{eq:cohofunion}
\end{equation}
A key point in the argument is that $H^*_{(3)}(A) = 0$ when $A$ is any
intersection of two or more components.  
In lemma~\ref{lem-1-dimensional-Hodge-summand-chordal-case} we show that the first term is the Eisenstein
Hodge structure $\E = \Z[\omega]$, where all elements are of type
$(2,2)$, and in lemma~\ref{lem-the-interesting-Hodge-summand-chordal-case} we show that the second term is isomorphic to $H^1_{(6)}(C)$ as an
Eisenstein Hodge structure.  
Theorem~\ref{thm-main-thm-chordal-limits} follows.

\subsection{Semistable reduction}
\label{subsec-semistable-reduction}

We now construct a semistable model $\hat\V$ for 
$\V$.  Once constructed, it will be called $\V_3$.  
As noted above, our approach is to apply one base extension followed by
three blowups, the first of which is centered in $R$. 

Let $\V_0$ denote the result of base extension
of $\V$, where $t = s^6$ is the parameter substitution.    Thus
$$
\V_0=\bigl\{\bigl([x_0{:}\dots{:}x_4,z],s\bigr)\in P^5\times\D
| F+ s^6 G+z^3=0\bigr\}\;,
$$
with $s=0$ defining the central fiber $\V_0(0)$.  This central fiber
is just the fourfold $V$ defined by $F+z^3=0$, which is singular along
$R$.    

The next step is to blow up $\V_0$ by blowing up $P^5\times\D$ along
an ideal sheaf $\J$ supported in $R$.  Let $\I$ denote the ideal sheaf
of $R$ in $P^4$. Extend it, keeping the same name, to the ideal sheaf
on $P^5\times\D$ that defines
$$
\bigl(\hbox{the cone on $R$ in $P^5$, with vertex $[0{:}0{:}0{:}0{:}0{:}1]$}\bigr)
\times\D\;.
$$
Let $\J$ be the ideal sheaf
\begin{equation}
\label{eq-defn-of-ideal-sheaf-J}
\J=\langle \I^2, \I zs, \I s^3, z^3, z^2s^4, zs^4, s^6\rangle\;.
\end{equation}
Finally, let  $\widehat{P^5\times\D}$ be the blowup of $P^5\times\D$
along $\J$, and let $\V_1$ be the proper transform of $\V_0$.  The
central fiber $\V_1(0)$ is the pullback of $s=0$, and consists of the proper
transform $V'$ of $V$ and the exceptional divisor $E_1$ of
$\V_1\to\V_0$.  

 The family $\V_1$ is nearly semistable; one could call it an
 ``orbifold semistable model''.  The precise meaning of this is that
 the inclusion of the central fiber into $\V_1$ is locally modeled on
 the inclusion of a normal crossing divisor into a smooth manifold,
 modulo a finite group.  We will now proceed to resolve the quotient
 singularities. 

Let $S$ be the singular locus of $\V_1$.  This turns out to lie in $V'
\cap E_1$ and be a $P^1$ bundle over $R$.  Let $\K$
be the ideal sheaf generated by (i) the regular
functions vanishing to order 2 along $S$, and
(ii) the regular functions vanishing to order 3 at a generic point of
$V'$.  Let $\V_2$ be the blowup of $\V_1$ along $\K$.  The central
fiber $\V_2(0)$ consists of the proper transforms $V''$ and $E_1'$ of
$V'$ and $E_1$, together with the new exceptional divisor $E_2$.

The family $\V_2$ is also an ``orbifold semistable model.''  It has
quotient singularities along a surface $\Sigma$ which projects
isomorphically to $S$ and is smooth away from $\Sigma$.  Define $\V_3$
to be the ordinary blowup of $\V_2$ along $\Sigma$.  This is
the desired semistable model $\hat\V$.  Its central fiber $\V_3(0)$ is
as always the pullback of $s=0$, and consists of the proper transforms
$V'''$, $E_1''$ and $E_2'$, and the exceptional divisor $E_3$.  Of
course, $\sigma$ acts on $\V_0$ by $z\to\w z$.  Since it preserves
$\J$, it acts on $\V_1$, and similarly for $\V_2$ and $\V_3$.

\begin{theorem}
\label{thm-is-semistable-model}
$\V_3\to\D$ is a semistable model for the family
$\V \to\D$.
\end{theorem}

The rest of this subsection is devoted to proving this theorem and
describing the various varieties in enough detail for the cohomology
calculations in subsection~\ref{subsec-cohomology-computations}.  We
will discuss weighted blowups and give two lemmas, and then construct
and study the three blowups.  We provide more computational detail than
we would if our blowups were just ordinary blowups.

We describe weighted blowups as follows.  Let $P(a) =
P(a_1,\ldots,a_n)$ be the weighted projective space with weights $a
=(a_1,\ldots,a_n)$.  The \emph{weighted blowup} of $\C^n$ with weights
$a$ is the closure of the graph of the rational map $f: \C^n \dasharrow
P(a)$, where $f(z_1,\ldots,z_n) = [z_1{:}\cdots{:}z_n]$ is the natural
projection.

One can write down orbifold coordinate charts for the weighted blowup
as in \cite[pp.~166--167]{Corti-Kollar-Smith}.  However, note that the definition of
weighted blowup on p.~166
of that reference contains a misprint.
The map $f$ used to define the graph should be the one we are using, not $f(z_1,\ldots,z_n) =
[z_1^{a_1}{:}\cdots{:}z_n^{a_n}]$.

Weighted blowups are a special case of blowups of $\C^n$ at an ideal
supported at the origin. Suppose given weights $(a_1,\ldots,a_n)$.
Let $d$ be the least common multiple of the weights.  Let
$h(z_1,\ldots,z_n)$ be the vector of monomials of weight $d$ in some
order, considered as a rational map from $\C^n$ to $P^{N-1}$, where
$N$ the number of monomials.  Let $\Gamma_h$ be the closure of the
graph of $h$.  By the naturality of blowups,
$\Gamma$ is the blowup
of $\C^n$ along the ideal generated by the monomials of weight $d$.
Let $v(z_1,\ldots,z_n)$ be the same vector of monomials regarded as a
map from $P(a)$ to $P^{N-1}$.  This ``Veronese map'' is well-defined,
and it embeds $P(a)$ in $P^{N-1}$.  Finally, let $f$ be the natural
quotient (rational) map from $\C^n$ to $P(a)$, and let $\Gamma_f$ be
the closure of its graph.  Then $v\circ f = h$.  Since $v$ is an
embedding, it induces an isomorphism $\Gamma_f \cong \Gamma_g$.  To
summarize:

\begin{lemma} 
Suppose given weights $(a_1,\ldots,a_n)$. Let $d$ be their least
common multiple, and let $\J$ be the ideal generated by the monomials
of weight $d$.  Then the blowup of $\C^n$ at the origin with weights
$(a_1,\ldots,a_n)$ is isomorphic to the blowup of $\C^n$ along $\J$.
\label{lemma:weightedblowup}
\qed
\end{lemma}

Let $X \subset \C^n$ be an analytic hypersurface with equation $f_d(z)
+ f_{d+1}(z) + \cdots$, where $f_k$ has weight $k$ and $d > 1$.  Then
the proper transform $\hat X$ of $X$ under the weighted blowup is
obtained as follows: delete the origin and replace it by the
hypersurface in $P(a)$ defined by $f_d(z) = 0$. That is, replace the
origin by the weighted projectivization of the weighted tangent cone
of $X$.  We call $\hat X$ the weighted blowup of $X$ at the point
corresponding to the origin.

Next we give a local equation for the chordal cubic at a point of
its singular locus.

\begin{lemma} Let $T$ be the chordal cubic, $R$ its rational normal curve,
and  $P$ a point of $R$.  Then there are  
local analytic coordinates $x,u,v,w$ on a neighborhood of $P$ in $P^4$,
in which $R$ is defined by $u = v = w = 0$ and $T$ by
\begin{equation}
\label{lemma:localeq}
  u^2 + v^2 + w^2 = 0.
\end{equation}
\end{lemma}

\begin{proof} 
A hyperplane $H$ transverse to $R$ at $P$ meets $T$ in a cubic surface
with a node at $P$.  Consequently there are analytic coordinates $u,
v, w$ on $H$ near $P$ such that the equation for $T \cap H$ is $u^2 + v^2 + w^2
= 0$.  Now, there exists a $1$-parameter subgroup $X$ of $\aut
T=\PGL(2,\C)$ that moves $P$ along $R$, and we take $x$ as a coordinate on
$X$.  The natural map $T\times X\to P^4$ gives the claimed local
analytic coordinates.
\end{proof}

\begin{remark}
If $H$ is a hyperplane in $P^4$ with a point of contact of order~$4$
with $R$, then there exist global algebraic coordinates on $P^4-H$
such that $T-H$ is given by \eqref{lemma:localeq}.  The proof is an unenlightening
sequence of coordinate transformations.
\end{remark}

{\bf First blowup.}  
One checks that $\V_0\sset P^5\times\D$ is smooth
away from $R$ and transverse to $\{s=0\}$ away from $R$.  To study
$\V_0$ near a point $P$ of $R$, introduce local analytic
coordinates $x,u,v,w,z,s$ around $P$ as follows.  Begin with affine
 coordinates $x,u,v,w$ around $P$ in $P^4$,  as in Lemma \ref{lemma:localeq}.
 Let  $z$ and $s$ be
the same $z$ and $s$ used above.  
 Then a neighborhood of $P$ in $\V_0$ is
\begin{align*}
\bigl\{(x,u,v,w,z,s)&{}\in(\hbox{some open set in $\C^6$})\bigm|\\
&u^2+v^2+w^2+s^6G(x,u,v,w)+z^3=0\bigr\}\;.
\end{align*}
If $P\in B$ then we may use the transversality of $R$ and $\{G=0\}$ to
change coordinates by $x\to x+(\hbox{a function of $u,v,w$})$, so that
$\{G=0\}$ is the same set as $\{x=0\}$.  That is, $G$ is $x$ times a
nonvanishing function.  Absorb a sixth root of this nonvanishing
function into $s$.  This yields a neighborhood of $P$ of the form
\begin{equation}
\label{eq-V0-in-local-analytic-coordinates}
\begin{split}
\bigl\{(x,u,v,w,z,s)&{}\in(\hbox{some open set in $\C^6$})\bigm|\\
&u^2+v^2+w^2+s^6x+z^3=0\bigr\}\;,
\end{split}
\end{equation}
where $\I=\langle u,v,w\rangle$ and $\sigma$ acts by $z\to\w z$.  If
$P\notin B$ then the same analysis leads to
\eqref{eq-V0-in-local-analytic-coordinates}, but with the $s^6x$ term
replaced by $s^6$.  We will not discuss this case, because it is
implicitly treated in the $P\in B$ case (by looking at points of $R$
near points of $B$).

We defined the ideal sheaf $\J$ in \eqref{eq-defn-of-ideal-sheaf-J},  defined
$\widehat{P^5\times\Delta}$ as the blowup of $P^5\times\Delta$
along $\J$, and  $\V_1$ as the proper transform of $\V_0$.  To
understand its geometry, we write down generators for $\J$ in our local
coordinates and recognize $\widehat{P^5\times\D}$ as a weighted blowup
along $R$.  Since $\I=\langle u,v,w\rangle$, $\J$ is generated by
the monomials of weight six, where we give the variables $u,v,w,z,s$
the weights $3,3,3,2,1$.  By lemma~\ref{lemma:weightedblowup}, $\V_1$ is the
weighted blowup of $\V_0$ along the $x$-axis.  By the discussion
following that lemma, above a neighborhood of $P$, the
exceptional divisor $E_1$ is
\begin{equation}
\label{eq-E1-in-C-times-P33321}
\begin{split}
\bigl\{&\bigl(x,[u{:}v{:}w{:}z{:}s]\bigr)\in\hbox{(an open set in $\C$)}\times
P(3,3,3,2,1)\bigm|\\
&u^2+v^2+w^2+z^3+xs^6=0\bigr\}\;.
\end{split}
\end{equation}
That is, over $R-B$, $E_1$ is a smooth fibration with fiber
isomorphic to 
the hypersurface in $P(3,3,3,2,1)$ defined by 
\begin{equation}
\label{eq-generic-fiber-of-E1-over-R}
u^2 + v^2 + w^2 + z^3 + s^6 = 0\;.  
\end{equation}
The special fibers lie above the points of $B \subset R$ and are
copies of the weighted homogeneous hypersurface with equation $u^2 +
v^2 + w^2 + z^3 = 0$.

The method we use for our detailed coordinate computations follows
\cite[pp.~166--167]{Corti-Kollar-Smith}.  The blowup is covered
by open sets $B_u$, $B_v$, $B_w$, $B_z$ and $B_s$, each the quotient
of an open set $A_u,\dots,A_s$ of $\C^6$ by a group of order $3,3,3,2$
or $1$.  We will treat $A_u\to B_u\sset\widehat{P^5\times\D}$ in some
detail; the treatment of $v$ and $w$ is the same, and we
will only briefly comment on $B_z$ and $B_s$.  One may take
coordinates $x,\udot,\vdot,\wdot,\zdot,\sdot$ on $A_u$, with
$B_u=A_u/\langle\eta\rangle$, where $\eta$ is the transformation
$\udot\to\w\udot$, $\zdot\to\w\zdot$, $\sdot\to\wbar\zdot$.  The map
$A_u\to P^5\times\D$ is given by $u=\udot^3$, $v=\vdot\udot^3$,
$w=\wdot\udot^3$, $z=\zdot\udot^2$ and $s=\sdot\udot$; these functions
of $\udot,\dots,\sdot$ are $\eta$-invariant, so they  define a
map $B_u\to P^5\times\D$.  The preimage in $A_u$ of the exceptional
divisor of $\widehat{P^5\times\D}\to P^5\times\D$ is $\{\udot=0\}$.  The
pullback to $A_u$ of the defining equation of $\V_0$ is
$\udot^6(1+\vdot^2+\wdot^2+\sdot^6x+\zdot^3)=0$, so the preimage in
$A_u$ of the proper transform $\V_1$ is the hypersurface $H$ given by
$1+\vdot^2+\wdot^2+\sdot^6x+\zdot^3=0$.  One checks that this is a
smooth hypersurface.  Also, the pullback to $A_u$ of the central fiber
$s=0$ is $\sdot\udot=0$, and one checks that $H$ meets $\{\udot=0\}$,
$\{\sdot=0\}$ and $\{\udot=\sdot=0\}$ transversely.  Therefore $s=0$
defines a NCD in $H$ with smooth components of multiplicity one.  The
one complication is that $\V_1\cap B_u$ is $H/\langle\eta\rangle$
rather than $H$ itself.  Since
$\eta$ has fixed points in $H$, $\V_1$ turns out to be singular.  We
will resolve these singularities after discussing $B_z$ and $B_s$.

The analysis of $\V_1\cap B_s$ is easy, since $B_s=A_s$.   One writes down
the equation for $\V_1\cap A_s$, and checks smoothness and
transversality to  $\{s=0\}$.  In fact, $A_s$ meets only one
component of the central fiber, $E_1$.    The analysis of $\V_1\cap B_z$
is only slightly more complicated.  We have $B_z=A_z/(\Z/2)$.  One
writes down the equation for the preimage of $\V_1$ in $A_z$, and checks
its smoothness and that $s=0$ defines a NCD with smooth components.
Then one observes that this hypersurface misses the fixed points of
$\Z/2$ in $A_z$, so that the same conclusions apply to $\V_1\cap B_z$.

Now we return to $\V_1\cap B_u$.  It is convenient to make the
coordinate change $\vdot\to(\vdot+\wdot)/2$,
$\wdot\to(\vdot-\wdot)/2i$ on $A_u$, so that $\vdot^2+\wdot^2$ is
replaced by $\vdot\wdot$.  Then the fixed-point set of $\eta$ in
$H\sset A_u$ is $\{(x,0,\vdot,\wdot,0,0)|1+\vdot\wdot=0\}$.  Therefore
every fixed point has $\vdot\neq0$, and we can use
$1+\vdot\wdot+\sdot^6x+\zdot^3=0$ to solve for $\wdot$ in terms of the
other variables.  The result is that there is an open set $W$ in
$\C^5$, with coordinates $x,\udot,\vdot,\zdot,\sdot$, mapping
isomorphically onto its image in $H$, with its image containing all
the fixed points of $\eta$ in $H$.  The induced action of $\eta$ on $W$ is
$(x,\udot,\vdot,\zdot,\sdot)\to(x,\w\udot,\vdot,\w\zdot,\wbar\zdot)$.
Therefore, a point on the singular locus of $\V_1$ has a neighborhood
in $\V_1$ locally modeled on $(W\sset\C^5)/\diag[1,\w,1,\w,\wbar]$.
Therefore the singular locus of $\V_1$ is a smooth surface $S$.  The
preimage in $W$ of the central fiber $s=0$ is $\sdot\udot=0$.  Note
that this is $\eta$-invariant, hence well-defined on
$W/\langle\eta\rangle$.  The preimage in $A_u$ of the exceptional
divisor $E_1$ is $\udot=0$ and that of $V'$ is $\sdot=0$.  Note that
neither $\udot$ nor $\sdot$ is $\eta$-invariant, so these equations
don't make sense on $W/\langle\eta\rangle$.  This reflects the fact 
that neither $E_1$ nor $V'$ is a Cartier divisor, but their sum is.
(It should not alarm the reader that $E_1$ is not Cartier, even though
blowing up an ideal always gives a Cartier divisor; it is $3E_1$
rather than $E_1$ which is Cartier.)
Finally, a lift to $W$ of $\sigma:\V_1\to\V_1$ is $\zdot\to\w\zdot$.
This shows that $\sigma$ acts trivially on $S$.

\emph{Geometry of the central fiber.}  The central fiber of $V_1$ is
the variety
\[
  \V_1(0) = V' \cup E_1,
\]
where the exceptional divisor has already been described as a
fibration over $R$.  The components have multiplicity one.  The first
component is a blowup of $V$ along $S$.  The exceptional divisor of
$V'\to V$ is $V'\cap E_1$, which in coordinates is given by
\eqref{eq-E1-in-C-times-P33321} with the extra condition $s=0$.
Therefore, $V'\cap E_1$ is a fiber bundle over $R$ with fiber
isomorphic to the hypersurface $u^2+v^2+w^2+z^3=0$ in $P(3,3,3,2)$.
One can check that this hypersurface is a copy of $P^2$ (by projecting
away from $[0{:}0{:}0{:}1]$).

{\bf Second blowup.}  Our aim is now to resolve the singularities of
$\V_1$.  Henceforth, we will work with the local description of $\V_1$
in terms of $W/\langle\eta\rangle$, rather than regarding $\V_1$ as a
hypersurface in a 6-dimensional space.  We will see below (lemma~\ref{lem-ideal-sheaf-K-gives-desired-blowup}) that
blowing up the ideal sheaf $\K$ to get $\V_2$ can be described in
terms of $W/\langle\eta\rangle$ as follows.  Recall that $W$ has
coordinates $x,\udot,\vdot,\zdot,\sdot$ with
$\eta=\diag[1,\w,1,\w,\wbar]$.  Define $\What$ to be the weighted
blowup of $W$ along the $x$-$\vdot$ plane, where $\udot$ and $\zdot$
have weight 1 and $\sdot$ has weight 2.  It is natural to use these
weights, because with respect to them, $\eta$ is weighted-homogeneous,
hence fixes the exceptional divisor pointwise.  Then the map
$$
\Bigl(\hbox{the preimage  of
  $W/\langle\eta\rangle\sset\V_1$ in $\V_2$}\Bigr)
\to W/\langle\eta\rangle
$$ turns out to be equivalent to $\What/\langle\eta\rangle\to
W/\langle\eta\rangle$.  This is the content of lemma~\ref{lem-ideal-sheaf-K-gives-desired-blowup}.  As before, one covers $\What$ by open
sets $B_\udot$, $B_\zdot$ and $B_\sdot$, which are quotients of open
sets $A_\udot$, $A_\zdot$ and $A_\sdot$ of $\C^5$ by cyclic groups of
orders $1$, $1$ and $2$.  By working in local coordinates, one can
check that $\eta$ acts on $B_\udot=A_\udot$ by multiplying a single
coordinate by $\w$.  Therefore $B_\udot/\langle\eta\rangle\sset\V_2$
is smooth.  One also checks that $s=0$ defines a NCD with smooth
components of multiplicity one.  (The multiplicity is
three along the preimage of $E_2$ in $B_\udot$, but in
$B_\udot/\langle\eta\rangle$ the multiplicity is only one.)  Exactly
the same considerations apply to
$B_\zdot/\langle\eta\rangle\sset\V_2$.

A similar analysis leads to the conclusion that
$B_\sdot/\langle\eta\rangle$ is isomorphic to an open set in $\C^5$,
modulo $\Z/2$, acting by negating three coordinates.  Therefore $\V_2$
is singular along a surface $\Sigma$ which maps isomorphically to $S$.
(It follows that $\sigma$ acts trivially on $\Sigma$.)  Furthermore,
at a singular point, $\V_2$ is locally modeled on
$$
\C^2\times\bigl(
\hbox{the cone in $\C^6$ on the Veronese surface $P^2\sset
P^5$}
\bigr)\;.
$$ 
One checks that away from $\Sigma$, the central fiber $\V_2(0)$ in
$B_\sdot/\langle\eta\rangle$ is a NCD with smooth components of
multiplicity one.

\emph{Geometry of the central fiber.}  The central fiber of $\V_2$ is
the variety
\[
  \V_2(0) = V''\cup E_1' \cup E_2.
\]
One can check the following:

(a) The new exceptional divisor $E_2$ is a $P(1,1,2)$-bundle over $S$.
Note that $P(1,1,2)$ is isomorphic to a cone in $P^3$ over on a smooth
plane conic.  The vertices of these cones comprise $\Sigma$.

(b) The map $V'' \to V'$ has exceptional divisor $V''\cap E_2$, which
is a $P^1$-bundle over $S$.  (Each fiber $P^1$ is a smooth section of
the quadric cone introduced in (a).  In particular, $V''\cap\Sigma=\emptyset$.)

(c) The map $E_1'\to E_1$ has exceptional divisor $E_1'\cap E_2$,
which is a $P^1$-bundle over $S$.  (Each fiber $P^1$ is a line through
the vertex of the cone in (a).)

(d) $V''\cap E_1'$ projects diffeomorphically to $V'\cap E_1$.  This
is obvious away from $S$, and is true over $S$ because the fiber of
$V''\cap E_1'$ over a point of $S$ is the intersection of the $P^1$'s
in (b) and (c), which is a point.

(e) The triple intersection $V''\cap E_1'\cap E_2$ is isomorphic to
$S$. 

{\bf Third blowup.}  The family $\V_2$ is singular along the surface
$\Sigma$, where it is locally modeled on $\C^2\times\bigl(\C^3/\{\pm
I\}\bigr)$. A single ordinary blowup along $\Sigma$, i.e., blow up
along the ideal sheaf of $\Sigma$, resolves the singularity, and it is
obvious that $E_3$ is a $P^2$-bundle over $\Sigma$.  One checks that
the central fiber is a NCD with smooth components of multiplicity one.
A convenient way to do the computations is to take the ordinary blowup
of $\C^5$ along $\C^2$, and then quotient by $\Z/2$.  (This gives the
same blowup.)

\emph{Geometry of the central fiber.}  The central fiber of $V_3$ is the variety
\[
  \V_2(0) = V'''\cup E_1''\cup E_2' \cup E_3.
\]
The components are of multiplicity one and smooth and transverse.  
We have already observed that $E_3$ is a $P^2$-bundle over
$\Sigma$.  Its intersections with $E_1''$ and $E_2'$ are $P^1$-bundles
over $\Sigma$; in fact the $P^1$'s are lines in the $P^2$'s.  The
triple intersection $E_3\cap E_2'\cap E_1''$ projects isomorphically
to $\Sigma$.  $E_3$ does not meet $V'''$, so $V'''=V''$ and the
intersections of $V'''$ with $E_1''$ and $E_2'$ are the same as
$V''$'s intersections with $E_1'$ and $E_2$, described earlier.

All that remains for the proof of
theorem~\ref{thm-is-semistable-model} is to explain why blowing up the
ideal sheaf $\K$ coincides with our description
$\What/\langle\eta\rangle\to W/\langle\eta\rangle$ of the second
blowup.  We formalize this as a lemma:

\begin{lemma}
\label{lem-ideal-sheaf-K-gives-desired-blowup}
Suppose $W$ is an open set in $\C^5$, with coordinates $x$, $\udot$,
$\vdot$, $\zdot$ and $\sdot$, and is invariant under
$\eta=\diag[1,\w,1,\w,\wbar]$.  Suppose $V'\sset W/\langle\eta\rangle$
is the image of $\{\sdot=0\}\sset W$, and let $\K$ be the ideal sheaf
on $W/\langle\eta\rangle$ generated by the regular functions vanishing
to order~$3$ at a generic point of $V'$ and the regular functions
vanishing to order~$2$ along the singular locus $S$ of
$W/\langle\eta\rangle$.  Let $\What$ be the weighted blowup of $W$
along $\udot=\zdot=\sdot=0$, with $\udot$ and $\zdot$ having
weight~$1$ and $\sdot$ having weight~$2$.  Then $W\to
W/\langle\eta\rangle$ induces an isomorphism from
$\What/\langle\eta\rangle$ to the blowup of $W/\langle\eta\rangle$
along $\K$.
\end{lemma}

\begin{proof}
We begin by observing that the nine monomials
$x$, $\vdot$, $\udot^3$, $\udot^2\zdot$, $\udot\zdot^2$, $\zdot^3$,
$\udot\sdot$, $\zdot\sdot$ and $\sdot^3$ generate the invariant ring
of $\eta$, so that evaluating them embeds $W/\langle\eta\rangle$ in
$\C^9$.  It is easy to see that $\K$ is the ideal generated by (i) the
quadratic monomials in $\udot^3,\udot^2\zdot,\dots,\sdot^3$, and (ii)
the linear function $\sdot^3$.  (We remark that the function $\sdot^3$
vanishes to order $3$ along $V'$, at a generic point of $V'$, and
generates the ideal of such functions.  However, even though $S$ lies
in $V'$, $\sdot^3$ only vanishes to order one at a point of $S$, because
$\sdot^3$ is one of the coordinate functions on $\C^9$.  This is
related to the fact that $V'$ is not a Cartier divisor, but $3V'$ is.)
Therefore, blowing up $\K$ amounts to defining $\widehat{\C^9}$ as the
weighted blowup of $\C^9$ along the $x$-$\vdot$ plane, with weights
$1,\dots,1,2$, and taking the proper transform therein of
$W/\langle\eta\rangle$.  One can cover $\What$ and $\widehat{\C^9}$ by
open sets, write down the rational map from $\What$ to
$\widehat{\C^9}$ explicitly, and check that it is regular and induces
an embedding $\What/\langle\eta\rangle\to\widehat{\C^9}$.
\end{proof}

\begin{remark}
The proof conceals the origin of the choice of $\K$.  We found it as
follows.  We knew we wanted to take the $(1,1,2)$ weighted blowup of
$W$, which is to say, blow up the ideal
$\langle\udot^2,\udot\zdot,\zdot^2,\sdot\rangle$.  The problem is that
its generators are not $\eta$-invariant, so they do not define
functions on $W/\langle\eta\rangle$.  The solution was to blow up the
cube of this ideal rather than the ideal itself.  We wrote down the
generators of the cube, and expressed them in terms of the generators
for the invariant ring of $\eta$.  Every generator of the ideal was a
quadratic monomial in the generating invariants, except for $\sdot^3$,
which was itself one of the generating invariants.  This suggested
that our weighted blowup of $\C^9$ would give the desired blowup
$\What/\langle\eta\rangle\to W/\langle\eta\rangle$.  Then we just checked the
construction.
\end{remark}

\subsection{Cohomology computations}
\label{subsec-cohomology-computations}

In this subsection we describe the $3$-part of the middle cohomology
of the central fiber $\hat V$ of the semistable model $\V_3$, in a
sequence of lemmas.

\begin{lemma} 
\label{lem-Leray-sequence-for-3-parts}
Let $M'$ and $M$ be algebraic varieties on which an automorphism
$\sigma$ of finite order acts.  Let $f: M' \to M$ be a morphism with
connected fibers which is equivariant with respect to this action.
Suppose $S \subset M$ is a subspace, let $D = f^{-1}(S)$, and assume
that (a) $f: M' - D \to M -S$ is an isomorphism, (b) $f: D \to S$ is a
fiber bundle, and (c) $\sigma$ acts trivially on $S$ and trivially on
the rational cohomology of the fiber.
Then the map of characteristic subspaces
\[
   H^*_{(k)}(M;\Q) \to H^*_{(k)}(M';\Q)
\]
is an isomorphism for each $k \neq 1$.
\label{eq:isomorphism-on-k-parts-of-cohomology}
\end{lemma}

\begin{proof}
Consider the Leray spectral sequence for the map $f: M' \to M$.  Its
abutment is the cohomology of $M'$.  Its $E_1$ terms are of the form
$E_1^{p,q} = H^p(M,R^qf_*\Q)$.  The automorphism $\sigma$ acts on $M$,
$M'$, and the spectral sequence.  Thus we can speak of the spectral
sequence for the $k$-part of the cohomology of $M'$.  Because $f$ is
an isomorphism over $M - S$ and has connected fibers over $S$, there
are two kinds of terms of the spectral sequence.  For the first, the
support of the coefficient sheaf is all of $M$.  These are the terms
$H^p\bigl(M,(f_*\Q)_{(k)}\bigr)$, which are isomorphic to
$H^p_{(k)}(M,\Q)$.  For the other terms the support of the coefficient
sheaf is $S$. These are the terms $H^p\bigl(S, (R^qf_*\Q)_{(k)}\bigr)$
with $q > 1$.  However, since the action on the cohomology of the
fibers is trivial, these spaces are zero.  Therefore the spectral
sequence degenerates and we have the isomorphism
\[
   H^*_{(k)}(M,\Q) \cong H^*_{(k)}(M,f_*\Q) \cong H^*_{(k)}(M',\Q).
\]
This completes the proof.
\end{proof}

\begin{lemma}
\label{lem-central-fiber-coho-chordal} 
The 3-part of the middle cohomology of the central fiber of $\V_3$, with $\Q$
coefficients,  is 
\[
   H^4_{(3)}(\hat V;\Q) \cong H^4_{(3)}(V;\Q) \oplus H^4_{(3)}(E_1;\Q).
\]
\end{lemma}

\begin{proof}
First, if $A$ is an intersection of two or more components of $\hat V$,
then $H^4_{(3)}(A;\Q)=0$.  One sees this by using
lemma~\ref{lem-Leray-sequence-for-3-parts} and the explicit
descriptions of the intersections, given in subsection~\ref{subsec-semistable-reduction}.  That
is, each of $E_3\cap E_2'$ and $E_3\cap E_1''$ is a $P^1$-bundle over
$\Sigma$, and $\sigma$ acts trivially on $\Sigma$.  Also, $E_3\cap
E_2'\cap E_1''$ projects isomorphically to $\Sigma$.  Finally, the
nonempty intersections involving $V'''$ are $V'''\cap E_1''$,
$V'''\cap E_2'$ and $V'''\cap E_1''\cap E_2'$, to which essentially
the same argument applies.  Then the Mayer-Vietoris spectral sequence
implies that with $\Q$ coefficients we have
$$
H^4_{(3)}(\hat V)=H^3_{(3)}(V''')
\oplus H^3_{(3)}(E_1'')
\oplus H^3_{(3)}(E_2')
\oplus H^3_{(3)}(E_3).
$$ We continue to use lemma~\ref{lem-Leray-sequence-for-3-parts}.
$E_3$ is a $P^2$-bundle over $\Sigma$, so the last term vanishes.
Also, $E_2'\to E_2$
induces an isomorphism on $H^4_{(3)}({}\cdot{};\Q)$, and similarly for
$E_1''\to E_1'\to E_1$ and $V'''=V''\to V'\to V$.  Finally, $E_2$ is a
$P(1,1,2)$-bundle over $S$, and $\sigma$ acts trivially on $S$, so
$H^4_{(3)}(E_2;\Q)=0$ and the lemma follows.
\end{proof}

\begin{lemma}  There is an isomorphism
\label{lem-1-dimensional-Hodge-summand-chordal-case}
\[
    H^4_{(3)}(V) \cong \E,
 \]
 where $\E$ is the rank one free Eisenstein module of type $(2,2)$.
\end{lemma}

\begin{proof}
By Morse theory, $P^5$ is obtained from a neighborhood of $V$ by
attaching cells of dimension five and larger.  Thus the cohomology of
$V$ is the same as that of $P^5$ in dimensions less than four.  It is
known that $V$ is a topological manifold (it follows from lemma~\ref{lemma:localeq}
that near $R$, $V$ is modeled on a sum of three squares and a cube;
apply the first example of \S9 of \cite{milnor}).  Therefore
Poincar\'e duality implies that $H^4(V)$ is torsion-free and that $V$'s
cohomology in all dimensions except four is that of $P^4$.

A chordal cubic $T$ has three natural strata: (i) $R$ itself,
 (ii) the union of the tangent
lines to $R$, minus their points of tangency, and (iii) the union of the
secant lines, minus the points in which they cut the $R$.  These
strata are a $P^1$, a $\C$ bundle over $P^1$, and a $\C^*$-bundle over
$P^1$.  Adding the Euler characteristics of the strata, we find that
$\chi(T) = 4$.  Therefore the Euler characteristic of $V$, a
three-sheeted cover of $P^4$ branched along $T$, is seven.  It follows
that $H^4(V)$ is free of rank three.  Since any $\sigma$-invariant
cohomology pulls back from $P^4$, we conclude that $H^4_{(3)}(V)$ has
rank two.  It is therefore a one-dimensional Eisenstein module,
necessarily of type $(2,2)$.
\end{proof}

{\bf Note.}  
Suppose that $\omega_t$ is a family of classes in
$H^{3,1}(V_t)_{\bar\omega}$.  Let $\delta_t$ be a family of homology
classes whose limit as $t$ approaches zero is a generator of the
Eisenstein factor in the above decomposition.  Then
\[
  \lim_{t \to 0} \int_{\delta_t} \omega_t = 0,
\]
since the $\E$-component of the limit of $\omega_t$ has type $(2,2)$.
Thus the linear equation
\[
    \int_\delta \omega = 0
\]
is the condition for a period vector to represent a chordal cubic.  In
other words, the limiting periods for degenerations to chordal cubics
near $V$ lie in the hyperplanes $\delta^\bot$.

Our final lemma treats the ``interesting'' part of the Hodge structure
of the central fiber, describing it in terms of that of the curve $C$
from \eqref{eq-defn-of-C}.

\begin{lemma} There is an isogeny 
\label{lem-the-interesting-Hodge-summand-chordal-case}
\[
    H^4_{(3)}(E_1) \to H^1_{(6)}(C)
\]
of Eisenstein Hodge structures, of weight $(-2,-1)$.
\end{lemma}

\begin{proof}
It suffices to exhibit an isomorphism 
\begin{equation}
\label{eq-bar}
H^4_{(3)}(E_1,\Q)\isomorphism H^1_{(6)}(C;\Q)
\end{equation}
that carries each component $H^{p,q}_{(3)}(E_1)$ of the Hodge
decomposition, $p+q=4$, to
$H^{p-2,q-1}_{(3)}(C)$.  The proof uses the following idea.  There are
maps $p: E_1 \to R$ and $p': C \to R$.  Thus we can compute the
3-parts of the cohomology of the varieties in question by using the
Leray spectral sequence.  The only relevant initial terms of the
spectral sequences are $H^1_{(3)}(R, R^3p_*\Q)$ and $H^1_{(6)}(R, p'_*\Q)$,
respectively.  The coefficient sheaves are isomorphic, and this leads
to the isomorphism \eqref{eq-bar}.

Consider $(p'_*\Z)_{(6)}$.  It is clear that its restriction to $R-B$
is a local system of abelian groups isomorphic to $\Z^2$.  This
becomes a local system of 1-dimensional $\E$-modules under our
definition of the action of $-\wbar\in\E$ as $\zeta^*$, where $\zeta$
is from \eqref{eq-defn-of-zeta-the-order-6-deck-transformation}.  If
$\gamma$ is a loop in $R-B$ encircling a single point of $B$ once
positively, then one can work out monodromy of $\gamma$ on the local
system by the method used for
lemma~\ref{lem-meridians-act-by-complex-reflections}.  The
degeneration around the point of $B$ is described by $z^6=t$, and we
will use the Sebastiani-Thom argument, as in the proof of
lemma~\ref{lem-meridians-act-by-complex-reflections}.  The result is
that the vanishing cycles are the 6-tuples summing to zero, indexed by
the 6th roots of unity, and the monodromy permuting these roots by
multiplication by $-\w$.  (Of course, the 3-part is a 2-dimensional
subspace.)  On the other hand, the action of $-\wbar\in\E$ is defined
to be that of $\zeta^*$.  Since the monodromy is the same as
$(\zeta^*)^{-1}$, it acts on the stalks by the $\E$-module
action of $(-\wbar)^{-1}=-\w$.  This
implies that the local system is supported away from $B$.  Passing to
rational coefficients, we obtain the local system $(p'_*\Q)_{(6)}$
of 1-dimensional vector spaces over $\E\tensor\Q$.

Now consider $(R^3p_*\Q)_{(3)}$.  Each fiber of $E_1\to R$ is a copy
of the hypersurface \eqref{eq-generic-fiber-of-E1-over-R} in
$P(3,3,3,2,1)$, whose middle cohomology can be worked out by the
Griffiths residue calculus \cite{tu}.  The result is that it is 2-dimensional
over $\Q$.  We defined the action of $\w\in\E$ to be that of $\sigma$,
which acts on $P(3,3,3,2,1)$ by $[u{:}v{:}w{:}z{:}s]\mapsto[u{:}v{:}w{:}\w z{:}s]$.
This makes $(R^3p_*\Q)_{(6)}$ into a local system of 1-dimensional
vector spaces over $\E\tensor\Q$.  Taking $\gamma$ as above, one can
work out its monodromy by the same method, with the result that it
acts by $-1\tensor-1\tensor-1\tensor\w$ in the notation used in the
proof of lemma~\ref{lem-meridians-act-by-complex-reflections}, while
$\sigma$ acts by $1\tensor1\tensor1\tensor\w$, which is to say, by the
$\E$-module action of $-\w$.  Therefore the monodromy is by $-\w$.
This shows that the local system is supported away from $B$, and that
the two local systems are isomorphic.

We claim next that the isomorphism \eqref{eq-bar} underlies an isomorphism
of complex Hodge structures.  A Hodge structure is defined via a
theory of differential forms and a harmonic theory.  These exist both
for cohomology with coefficients in the complex numbers $\C$ and with
coefficients in a local system, provided that the local system is
unitary.  This is the case for the eigensystems of the local systems
considered above.  Consider the local system $R^3p_*\Z$.  There is a
decomposition
\begin{equation}
    R^3p_*\Z\otimes\C = (R^3p_*\C)_{\omega} \oplus (R^3p_*\C)_{\bar\omega} 
\label{eq:rp3eigen}
\end{equation}
into eigenspaces of the deck transformation $\sigma$.  Each of the local systems on the right is unitary.
Moreover, each has a Hodge type inherited from the Hodge structure on the cohomology of the fiber.  To determine the relation between the types and the eigenvalues, consider the cohomology of a generic fiber $V$ with equation  $u^2 + v^2 + w^2 + z^3 + s^6 = 0$.  The differential form
\[
   \Phi = \frac{ \Omega(u,v,w,z)}{(u^2 + v^2 + w^2 + z^3 )^{11/6} }.
\]
is homogeneous of weight zero relative to the $\C^*$ action, and it generates $H^{2,1}(V)$.  It is clear that it is an eigenvector of $\sigma$ with eigenvalue $\omega$.  Thus 
\[
  (R^3p_*\C)_\omega =  (R^3p_*\C)_\omega^{2,1}.
\]
Consequently $H^1\bigl(R - B,  (R^3p_*\C)_\omega\bigr)$ carries a complex Hodge structure with types
$(3,1)$ and $(2,2)$, and there is an isomorphism of complex Hodge structures
\[
   H^4(E_1, \C)_\omega \cong H^1\bigl(R - B, (R^3p_*\C)_\omega^{2,1}\bigr) .
\]
Similar considerations yield an isomorphism
\[
  H^1_{(6),\omega}(C) \cong H^1_{(6),\omega}\bigl(R - B, (q_*\C)\bigr),
\]
where the decomposition is relative to the action of $\zeta^2$ and the local system can be considered as a local system of Hodge structures of type $(0,0)$.
Finally, the isomorphism
\[
   (R^3p_*\C)_\omega^{2,1} \cong  (q_*\C)_{(6),\omega}
\]
of unitary systems of complex Hodge structures gives an isomorphisms of complex Hodge structures
\[
  H^4_\omega(E_1, \C) \cong H^1_{(6),\omega}(C) .
\]
The proof is now complete.
\end{proof}

\begin{remark}
Another possible approach to the Hodge structure of $E_1$ is the following.
One can show that $E_1$ is birational to
the hypersurface $Z$ in $P(1,1,6,6,6,4)$ given by the equation
\[
   f(x,y) + u^2 +v^2 + w^2 + z^3 = 0.
\]
(This is plausible because the projection $(x,y,u,v,w,z) \mapsto
(x,y)$ blows up $Z$ along $\{x = y = 0\}$, yielding $\hat Z$.  And
$\hat Z$ fibers over $P^1$ with the same generic fiber and same special
fibers as in the fibration of $E_1$ over $R$.)  Presumably,
$E_1\dasharrow Z$ induces an isogeny of Hodge structures.  Standard
calculations with the Griffiths Jacobian calculus yields the Hodge
numbers $(0,1,9,0,0)$ for $H^4(Z) _\omega$, with the period map having
9-dimensional image.
\end{remark}

\section{Degeneration to a nodal cubic}
\label{sec-hodge-theory-nodal}

\noindent
The goal of this section is to identify the limit Hodge structure for
the degeneration of cyclic cubic fourfolds associated to a generic
nodal degeneration of cubic threefolds.  It is very similar to the
previous section, so we will be more brief.  The result of this
section that is used in the proof of the main theorem of the paper
(theorem~\ref{thm-MAIN-THEOREM}) is the following:

\begin{theorem}
\label{thm-discr-maps-onto-divisor}
The period map $g:\widehat{P\forms}_{ss}\to\PGamma\backslash\ch^{10}$ carries the
discriminant locus onto a divisor.
\end{theorem}

One should expect such a result, because at a node of a cubic
threefold $T$, the tangent cone is a cone over $P^1\times P^1$, and the
lines on $T$ through the node sweep out a $(3,3)$ curve in the
$P^1\times P^1$.  The generic genus~4 curve arises this way, providing
9 moduli.  Our
approach is to show that the interesting part of the limiting Hodge
structure is that of the K3 surface $K$ which is the 3-fold cover of
$P^1\times P^1$ branched over this curve.  This Hodge structure is the
same as the one studied by Kond\=o in his work on
moduli of genus four curves \cite{kondo:genus-4}.

Suppose $T$ is a generic nodal cubic threefold; we choose homogeneous
coordinates $x_0,\dots,x_4$ such that  the node is at
$[1{:}0{:}\dots{:}0]\in P^4$ and $T$ has defining equation
$$
F(x_0,\dots,x_4)=x_0(x_1x_2+x_3x_4)+f(x_1,\dots,x_4)=0\;,
$$ where $f$ is a homogeneous cubic.  Then $V\sset P^5$ is defined by
$F+x_5^3=0$, and we will project $V$ away from $P=[1{:}0{:}0{:}0{:}0{:}0]$ to
the $P^4$ with homogeneous coordinates $x_1,\dots,x_5$.  Every
reference to $P^4$ will be to this $P^4$.  We may suppose by
genericity that the intersection $K$ of $x_1x_2+x_3x_4=0$ and
$f(x_1,\dots,x_4)+x_5^3=0$ in $P^4$ is smooth.  The importance of this
surface is that projection away from the cusp of $V$ at
$P=(1,0,0,0,0,0)$ is a local diffeomorphism on $V-\{P\}$, except over
$K$, where it is a $\C$-bundle.  The notation reflects the fact that
$K$ is a K3 surface, since it is a smooth $(2,3)$-intersection in
$P^4$.  We will write $Q$ for the quadric $x_1x_2+x_3x_4=0$ in $P^4$.
The vertex of $Q$ is not in $K$, because $K$ is smooth.

By scaling the variables, we may suppose that $F+tx_0^3$ defines a
smooth threefold for all $t\in\D-\{0\}$.  This pencil of threefolds
is well-suited to our coordinates, but any pencil gives the same limiting Hodge
structure, provided that its generic member is smooth.  As in
section~\ref{sec-hodge-theory-chordal} we will write $\V$ for the associated family of cubic
fourfolds:
\begin{align*}
\V=\Bigl\{&\bigl([x_0{:}\dots{:}x_5],t\bigr)\in P^5\times\D\Bigm|\cr
&x_0(x_1x_2+x_3x_4)+f(x_1,\dots,x_4)+x_5^3+t\,x_0^3=0\Bigr\}\;.
\end{align*}
One checks that $\V$ is smooth away from the cusp $P$ of the central
fiber, and of course the central fiber itself is smooth away from $P$.
We let $\V_0$ be the degree six base extension, got by setting
$t=s^6$.  Since the monodromy of $\V$ on $H^4$ over $\D-\{0\}$ has
order~3, the monodromy of $\V_0$ is trivial.  Below, we will define
$\V_1$, $\V_2$ and $\V_3$ by repeated blowups.  The last blowup  $\V_3$ turns
out to be a semistable model for the degeneration.  The reason we take
six rather than three as the degree of the base extension is that this
choice makes the multiplicities of the components of the central fiber of
$\V_3$ all be $1$.

Regarding $x_1,\dots,x_5$ as affine coordinates on
$\C^5=\{x_0\neq0\}\sset P^5$, we have
\begin{align*}
\V_0\cap(\C^5\times\D)=\bigl\{&(x_1,\dots,x_5,s)\in\C^5\times\D
\bigm|\cr
&x_1x_2+x_3x_4+f(x_1,\dots,x_4)+x_5^3+s^6=0\bigr\}\;.
\end{align*}
By the Morse lemma, there is a neighborhood $W$ of the origin in
$\C^4$ with analytic coordinates $y_1,\dots,y_4$, such that the $x_i$
and $y_i$ agree to first order, and
$$
x_1x_2+x_3x_4+f(x_1,\dots,x_4)=y_1y_2+y_3y_4\;.
$$
To emphasize the analogy with the chordal case, we will write $z$ for
$x_5$.  Then
\begin{align*}
\V_0\cap(W\times\C\times\D)=\bigl\{&(y_1,\dots,y_4,z,s)\in
W\times\C\times\D\bigm|\cr
&y_1y_2+y_3y_4+z^3+s^6=0\bigr\}\;.
\end{align*}

The obvious thing to do is blow up the origin with weights
(3,\discretionary{}{}{}3,\discretionary{}{}{}3,\discretionary{}{}{}3,\discretionary{}{}{}2,\discretionary{}{}{}1).
To describe this intrinsically, we define $\I$ to be the 
ideal sheaf on $P^5\times\D$ defining the subvariety
$$
\bigl(\hbox{the line in $P^5$ joining $[1{:}0{:}0{:}0{:}0{:}0]$ and
  $[0{:}0{:}0{:}0{:}0{:}1]$}\bigr)
\times\D\;,
$$ and we define $\J$ as in \eqref{eq-defn-of-ideal-sheaf-J}.  We
define $\widehat{P^5\times\D}$ as the blowup of $P^5\times\D$ along
$\J$, and $\V_1$ as the proper transform of $\V_0$.  The exceptional
divisor in $\widehat{P^5\times\D}$ is a copy of
$P(3,\discretionary{}{}{}3,\discretionary{}{}{}3,\discretionary{}{}{}3,\discretionary{}{}{}2,\discretionary{}{}{}1)$,
and the exceptional divisor $E_1$ of $\V_1\to\V_0$ is the hypersurface
in $P(3,3,3,3,2,1)$ defined by
\begin{equation}
\label{eq-quasihomogeneous-eqn-for-E1}
y_1y_2+y_3y_4+z^3+s^6=0\;.
\end{equation}
As in section~\ref{sec-hodge-theory-chordal}, we will use primes to indicate proper transforms
of $V$ and the various exceptional divisors.  The exceptional
divisor of $V'\to V$ is the hypersurface in $P(3,3,3,3,2)$ defined by
\eqref{eq-quasihomogeneous-eqn-for-E1} and $s=0$.  This lets one see
that $\sigma$ acts trivially on $V'\cap E_1$, because
$\sigma=\diag[1,1,1,1,\w]$ acts on $P(3,3,3,3,2)$ in the same way as
the quasihomogeneous scaling by $\wbar$, which of course acts
trivially.  Finally, calculations strictly analogous to those of
section~\ref{sec-hodge-theory-chordal} show that $\V_1$ is smooth away from the surface $S$ in
$P(3,3,3,3)\isomorphism P^3$ defined by
\eqref{eq-quasihomogeneous-eqn-for-E1} and $s=z=0$.  $S$ is a
smooth quadric surface.  Furthermore, any point of $S$ admits a
neighborhood in $\V_1$ isomorphic to a neighborhood of the origin in
$\C^5$, modulo $\eta=\diag[\w,1,1,\w,\wbar]$.  This is exactly the
same as the local model near the singular set of section~\ref{sec-hodge-theory-chordal}'s
$\V_1$.

The second blowup $\V_2$ is the blowup of $\V_1$ along the ideal sheaf
$\K$, where $\K$ is defined exactly as in
section~\ref{sec-hodge-theory-chordal}: it is generated by the regular
functions which either vanish to order~2 
along $S$, or vanish to order~3 at a generic point of
$V'$.  Then the exceptional divisor $E_2$ is a $P(1,1,2)$-bundle over
$S$.  Now, $P(1,1,2)$ is isomorphic to a cone in $P^3$ over a smooth
plane conic, so each fiber has a singular point.  These turn out to be
singular in $\V_2$ as well, and constitute the entire singular locus
$\Sigma$ of $\V_2$, so $\Sigma$ is a copy of $S$.  Furthermore, near
$\Sigma$, $\V_2$ is locally modeled on $\C^2\times\bigl(\C^3/\{\pm
I\}\bigr)$, just as in section~\ref{sec-hodge-theory-chordal}.  It
turns out, also as in section~\ref{sec-hodge-theory-chordal}, that
$V''\cap\Sigma=\emptyset$.  Finally, we define $\V_3$ as the ordinary
blowup of $\V_2$ along $\Sigma$.  Then $E_3$ is a $P^2$-bundle over
$\Sigma$.  One can check that $\V_3$ is smooth, and that the central
fiber is a normal crossing divisor with smooth components $V'''$,
$E_1''$, $E_2'$ and $E_3$ of multiplicity one.

Now we will study the central fiber of $\V_3$, in order to determine
the limiting Hodge structure.  The first step is to rid ourselves of
most of the complication introduced by our blowups.  Then we will
study what remains, the Hodge structures of $V''$ and $E_1$.

\begin{lemma}
\label{lem-trivial-parts-of-(3)-cohomology-nodal-case}
$H_{(3)}^4\bigl(V'''\cup E_1''\cup E_2'\cup E_3;\Q\bigr)
=
H_{(3)}^4(V'';\Q)\oplus H_{(3)}^4(E_1;\Q)$.
\end{lemma}

\begin{proof}
This is analogous to lemma~\ref{lem-central-fiber-coho-chordal}.  
We use cohomology with rational coefficients throughout the proof.
For
the first step, the essential facts are the following. (i) The
restriction of the projection $\V_3\to\V_2$ to the central fiber is an
isomorphism to $V''\cup E_1'\cup E_2$, except over $\Sigma$. (ii) Over
$\Sigma$, it is a $P^2$-bundle.  (iii) $\sigma$ fixes $\Sigma$
pointwise.  Then lemma~\ref{eq:isomorphism-on-k-parts-of-cohomology}
implies
$$
H_{(3)}^4\bigl(V'''\cup E_1''\cup E_2'\cup E_3\bigr)
=
H_{(3)}^4\bigl(V''\cup E_1'\cup E_2\bigr)\;.
$$
For the second step, the essential facts are (iv) $E_2$ is a
$P(1,1,2)$-bundle over $S$, (v) $E_2\cap E_1'$ and $E_2\cap V''$ are
$P^1$-bundles over $S$, (vi) $E_2\cap E_1'\cap V''$ projects
isomorphically to $S$, and (vii) $\sigma$ fixes $S$ pointwise.  Then
Mayer-Vietoris implies
$$
H_{(3)}^4\bigl(V''\cup E_1'\cup E_2\bigr)\;
=
H_{(3)}^4\bigl(V''\cup E_1'\bigr)\;.
$$ For the third step, the essential fact is (viii) $\sigma$ acts
trivially on $V'\cap E_1$.  One can check that $V''\cap E_1'$ is the
proper transform of $V'\cap E_1$, so $\sigma$ acts trivially on
$V''\cap E_1'$.  Then Mayer-Vietoris implies
$$
H_{(3)}^4\bigl(V''\cup E_1'\bigr)
=
H_{(3)}^4(V'')\oplus H_{(3)}^4(E_1')\;.
$$
For the final step, the essential facts are (ix) $E_1'\cap E_2$ is a
$P^1$-bundle over $S$, so $E_1'\to E_1$ is a diffeomorphism except
over $S$, over which it is a $P^1$-bundle, and (x) $\sigma$ acts
trivially on $S$.  Then lemma~\ref{eq:isomorphism-on-k-parts-of-cohomology} implies
$H_{(3)}^4(E_1')=H_{(3)}^4(E_1)$.  This completes the proof.
\end{proof}

\begin{lemma}
\label{lem-Hodge-structure-of-E1-nodal-case}
$H_{(3)}^4(E_1;\Q)$ is $2$-dimensional,
of type $(2,2)$.
\end{lemma}

\begin{proof}
We regard $E_1$ as a hypersurface in $P(3,3,3,3,2,1)$ as in \eqref{eq-quasihomogeneous-eqn-for-E1}.
According to the Griffiths residue calculus \cite{tu},
the primitive cohomology of $E_1$ in dimension four
is obtained as follows.  Set $F = y_1y_2 + y_3y_4 + z^3 + s^6$.  Let 
\[
   \eta = y_1 \frac{\partial}{ \partial y_1} + \cdots +  y_4 \frac{\partial}{ \partial y_4} 
   +  z \frac{\partial}{ \partial z }  + s \frac{\partial}{ \partial s }
\]
be the Euler vector field, and let $\Omega$ be the contraction of
$\eta$ with $dy_1 \wedge \cdots \wedge dy_4 \wedge dz \wedge ds$.
Consider rational differential forms $A\Omega/F^k$ where the weight of
the numerator polynomial $A$ is chosen so that the rational
differential has weight zero.  Because $\Omega$ has weight fifteen,
the nonzero rational differentials with numerator polynomial of lowest
degree are those with $k = 3$ and $A = c_1 zs + c_2 s^3 + (\hbox{terms
in the $y_i$})$.  Let $res(A\Omega/F^k)$ be the Poincar\'e residue.
Such residues span the primitive cohomology in Hodge level $4-k+1$.
The Hodge level is greater than $4-k+1$ if and only if $A$ lies in the
Jacobian ideal of $F$; that is, the ideal generated by the partial
derivatives of $F$.  In the case at hand, we find that the primitive
part of $H^4$ is spanned by $res(sz\Omega/F^3)$ and
$res(s^3\Omega/F^3)$.  The automorphism $\sigma$ acts on these by
multiplication by $\omega^2$ and $\omega$, respectively.  Therefore
$H^4_{(3)}(E_1;\Q)$ is 2-dimensional.  It clearly has type $(2,2)$.
\end{proof}

\begin{lemma}
\label{lem-topology-of-V''-nodal-case}
$V''$ admits a morphism to $P^4$ which is an isomorphism except over
the K3 surface $K$ and the vertex $v$ of the quadric cone $Q$.  Over
$K$, $V''$ is a smooth $P^1$-bundle, and the preimage of $v$ is
$E_1'\cap V''$.  
\end{lemma}

\begin{proof}
Let $\phi$ be the rational map $V\dasharrow P^4$ given by projection
away from the cusp $P$.  Then $\phi$ induces rational maps
$\phi':V'\dasharrow P^4$ and $\phi'':V''\dasharrow P^4$.  One can
check by lengthy local coordinate calculations that $\phi''$ is a
morphism, not just a rational map, and has the properties claimed in
the lemma.  We just give a summary.

Consider a line $\ell$ in $P^5$ through $P$.  If $\ell$ lies in $V$
then $\phi(\ell-\{P\})$ lies in $K$.  We write $D$ for the union of
the lines on $V$ through $P$.
If $\ell$ makes contact of
order only~2 with $V$ at $P$, then it meets $V$ at exactly one further point,
and does so transversely.  If it makes contact of order exactly~3,
then $\ell$ meets $V$ only at $P$, and the corresponding point of
$P^4$ is in $Q$ but not in $\phi(V-\{P\})$.  We  conclude that $\phi$
identifies $V-D$ with $P^4-Q$, and realizes $D-\{P\}$ as a $\C$-bundle
over $K\sset Q$.   We write
$D''$ for the proper transform of $D$ in $V''$.

In this proof we only care about divisors in $V''$, not in $\V_2$, so
we will write $e_1$ for $V'\cap E_1$, $e_1'$ for $V''\cap E_1'$ and
$e_2$ for $V''\cap E_2$.  One can check that $e_1'$ is the proper
transform of $e_1$.

Local calculations show that $\phi'$ is regular except along $S$, and
carries $e_1-S$ to $v$.  Further calculations show that $\phi''$ is
regular on all of $V''$, and therefore $\phi''$ carries $e_1'$ to $v$.
The behavior of $\phi''$ on $e_2$ is easy to understand.  We know that
$e_2$ is a $P^1$-bundle over the smooth quadric surface $S$, and that
$Q$ is a cone over a smooth quadric surface.  Therefore it is not
surprising that (i) $\phi''$ carries each fiber $P^1$ of $e_2$
isomorphically onto a line in $Q$ through $v$, and (ii) the only
points of $e_2$ that $\phi''$ carries to $v$ are in $e_2\cap e_1'$.
It follows from (ii) that the fiber of $\phi''$ over $v$ is exactly
$e_1'$.  

It remains to show that $\phi''$ is a diffeomorphism over
$P^4-(K\cup\{v\})$ and a smooth $P^1$-bundle over $K$.  Because
$\phi''$ maps $V''$ onto $P^4$, the image of $e_2$ must contain
$Q-(K\cup\{v\})$ and hence be equal to $Q$.  Since $P^4$ is smooth and
$D''$ and $e_1'$
are the only divisors in $V''$ that $\phi''$ can crush to
lower-dimensional varieties, Zariski's main theorem implies that
$\phi''$ is a diffeomorphism except over the image $K\cup\{v\}$ of
$D''\cup e_1'$.  So all that remains is to show that $\phi''^{-1}(K)$
is a smooth $P^1$-bundle over $K$.  First, by (i) above, each point of
$K$ has only one preimage in $e_2$, so each point of $k$ has preimage
$\C\cup\{\hbox{point}\}$ in $V''$.    It follows that $D''$ is the
full preimage of $K$.  Also, $\phi''$ restricts to an isomorphism
$e_2-e_1'\to Q-\{v\}$.  This implies that
$\phi''^{-1}(K)\cap e_2$ is a copy of
$K$.

$D''$ is obviously smooth away from $e_2$.  To show it is smooth at a
point $d$ of $D''\cap e_2$, it suffices to exhibit a smooth curve in
$V''$ through $d$ that is transverse to $D''$.  Since $\phi''^{-1}(K)$
is a copy of $K$, we can just take a curve in $e_2$ transverse to this
copy of $K$.  Since $e_2-e_1'$ maps isomorphically to its image, the
image curve is transverse to $K$, so the original curve must be transverse
to $D''$.  To show that $D''$ is a fibration over $K$, we will show
that the rank of $\phi''|_{D''}$ is $2$ everywhere.  This clearly
holds at every point of $D''-e_2$.  If the rank were${}<2$ at a point
$d$ of $D''\cap e_2$, then the rank of $\phi'':V''\to P^4$ at $d$
would be${}<3$, which is impossible since $\phi''|_{e_2}$ is a local
diffeomorphism at $d$.
\end{proof}

\begin{lemma} 
\label{lem-isogeny-with-K3-nodal-case}
The map $\phi'':V''\to P^4$ of lemma~\ref{lem-topology-of-V''-nodal-case} induces an isomorphism
$$
H^4_{(3)}(V'';\Q)\isomorphism H^2_{(3)}(K;\Q).
$$
\end{lemma}

\begin{proof}
 From Lemma
$\ref{lem-topology-of-V''-nodal-case}$, one knows that

(i) if $x \in P^4 - (K \cup \{v\})$, then $\phi''^{-1}(x)$  is a point;

(ii) The part of $V''$ over $K$ is a $P^1$-bundle;

(iii) $\phi''^{-1}(v)=E_1'\cap V''$ is a copy of the hypersurface
$y_1y_2+y_3y_4+z^3$ in
$P(3,\discretionary{}{}{}3,\discretionary{}{}{}3,\discretionary{}{}{}3,\discretionary{}{}{}2)$;
by projecting away from $[0{:}\dots{:}0{:}1]$, one checks that this is a
copy of $P^3$.  We saw above that $\sigma$ acts trivially on it.

It follows that the only terms
of the Leray spectral sequence which can contribute to
$H^4_{(3)}(V'')$ are

(i$'$) $H^4_{(3)}(P^4, \phi''_*\Q)$, 

(ii$'$)  $H^2_{(3)}(K, R^2\phi''_*\Q)$, 

(iii$'$) $H^0_{(3)}(v,R^4\phi''_*\Q)$.  

\noindent
Only the middle term (ii$'$)  is nonzero, and it is $H^2_{(3)}(K;\Q)$.
The lemma follows.
\end{proof}

\begin{lemma}
There is an isogeny of Eisenstein Hodge structures
$$
H^4_{(3)}(\hat V)\to H^2_{(3)}(K)\oplus\E
$$
of weight $(-1,-1)$, where $\E$ indicates the $1$-dimensional
Eisenstein lattice with Hodge structure pure of type $(1,1)$.
\end{lemma}

\begin{proof}
This follows from the isomorphisms
\begin{align*}
H^4_{(3)}(\hat V;\Q)&\isomorphism 
H^4_{(3)}(V'';\Q)\oplus
H^4_{(3)}(E_1;\Q)\\
{}&\isomorphism
H^2_{(3)}(K;\Q)\oplus\Q^2
\end{align*}
of lemmas~\ref{lem-trivial-parts-of-(3)-cohomology-nodal-case},
\ref{lem-Hodge-structure-of-E1-nodal-case} and
\ref{lem-isogeny-with-K3-nodal-case}, together with the fact that the
isomorphism $H^4_{(3)}(V'')\to H^2_{(3)}(K;\Q)$ is a map of Hodge
structures of weight $(-1,-1)$.  The latter fact comes from the fact
that the restriction of $V''\to P^4$ to the preimage of $K$ is a
$P^1$-bundle over $K$.
\end{proof}

Now, $K$ lies in the quadric cone $Q$, and projection away from the
vertex realizes $K$ as a 3-fold cover of $P^1\times P^1$ branched over
a $(3,3)$ curve $C$.  This curve has genus~4, and the generic genus~4
curve arises this way.  Kond\=o \cite{kondo:genus-4} showed that $H^2_{(3)}(K;\Q)$
has dimension~20, and that the Hodge structure of $K$ contains enough
information to recover $C$, up to an automorphism of $P^1\times P^1$.
Therefore, as $C$ varies over the $(3,3)$ curves in $P^1\times P^1$,
the associated K3 surfaces provide a 9-dimensional family of Hodge
structures. 

\begin{proof}[Proof of theorem~\ref{thm-discr-maps-onto-divisor}:]
This is essentially the same as theorem~\ref{thm-chordal-maps-onto-divisor}.  The previous lemma
and the exactness of the 3-part of the Clemens-Schmid sequence
$$
\dots\to
H^4_{(3)}(\V_3,\V_3^*;\Q)
\to
H^4_{(3)}(\hat V;\Q)
\to
\lim_{s\to 0} H^4_{(3)}(\hat V_s;\Q)
\to 0
$$
imply that the limit Eisenstein Hodge structure is isogenous to that
of $K$ (plus a 1-dimensional summand).  Then Kond\=o's work shows that the
Eisenstein Hodge structures form a 9-dimensional family.
\end{proof}

\section{The main theorem}
\label{sec-main-theorem}

Recall that  $\moduli_0^f$ and $\moduli_0$ are the moduli spaces of framed and
unframed smooth cubic threefolds.  We also define $\moduli_{ss}$ as
the GIT moduli quotient $P\forms_{ss}//\SL(5,\C)$, and $\moduli_s$ as
the corresponding stable locus.  Since we needed to blow up $P\forms$
before extending the period map, we also define
\begin{align*}
\widehat\moduli_s^f&{}=P\framed_s/PG\\
\widehat\moduli_s&{}=\widehat{P\forms}_s/PG\\
\widehat\moduli_{ss}&{}=\widehat{P\forms}_{ss}//\SL(5,\C)\;.  
\end{align*}
These moduli spaces fit into the commutative diagram:
$$
\begin{CD}
\moduli_0^f @>>> \widehat\moduli_s^f\\
@VVV @VVV\\
\moduli_0 @>>> \widehat\moduli_s @>>> \widehat\moduli_{ss}
\end{CD}
$$

We saw in lemma~\ref{lem-G-acts-freely-and-per-map-has-full-rank} that
$\moduli_0^f$ is a complex manifold.  Since $PG$ acts
properly on $\widehat{P\forms}_s$, we see that 
$\moduli_0$, $\widehat\moduli_s^f$ and $\widehat\moduli_s$ are analytic spaces.
We will see below that $\moduli_s^f$ is smooth.
As a GIT quotient, $\widehat\moduli_{ss}$ is a compact algebraic variety.

The period maps $P\framed_s\to\ch^{10}$ and
$\widehat{P\forms}_{ss}\to\overline{P\Gamma\backslash\ch^{10}}$ are
$PG$-invariant, and hence induce maps on $\widehat\moduli_s^f$ and
$\widehat\moduli_{ss}$.  The main theorem of the paper, theorem~\ref{thm-MAIN-THEOREM}, says
that the first of these maps is an isomorphism and that the second is
almost an isomorphism.  For the statement of the theorem, it is
convenient to refer to the discriminant locus of $\widehat\moduli_s^f$, by
which we mean the image of the discriminant locus of $P\framed_s$,
and similarly for the chordal locus.  Let $\H_{\D}$ (resp. $\H_c$) be
the union of the discriminant (resp. chordal) hyperplanes in $\ch^{10}$
(defined just before theorem~\ref{thm-period-map-extends-to-stable}), and
let $\H=\H_{\D}\cup\H_c$.

\begin{theorem}
\label{thm-MAIN-THEOREM}
The period map $g:\widehat\moduli_s^f\to\ch^{10}$ is an isomorphism.  It
identifies the discriminant (resp. chordal) locus of $\widehat\moduli_s^f$
with $\H_{\D}$ (resp. $\H_c$), and $\moduli_0^f$ with $\ch^{10}-\H$.  The
induced map $\widehat\moduli_{ss}\to\overline{\PGamma\backslash\ch^{10}}$ is
an isomorphism except over the $A_5$ cusp, whose preimage is a
rational curve.  Finally, this map induces an isomorphism of
$\moduli_s$ with $P\Gamma\backslash(\ch^{10}-\H_c)$.
\end{theorem}

\begin{proof}[Proof of theorem~\ref{thm-MAIN-THEOREM}:\/]
Since $\widehat\moduli_{ss}$ is compact and the period map carries
$\widehat\moduli_{ss}-\widehat\moduli_s$ to the boundary of
$\overline{P\Gamma\backslash\ch^{10}}$,
$\widehat\moduli_s\to\PGamma\backslash\ch^{10}$ is proper.  It follows that
$\widehat\moduli_s^f\to\ch^{10}$ is proper.  The discriminant locus is carried
into $\H_{\D}$ by theorem~\ref{thm-period-map-extends-to-stable}, and has
image a divisor by theorem~\ref{thm-discr-maps-onto-divisor}.
Therefore its image is a union of discriminant hyperplanes.  By
lemma~\ref{lem-transitivity-on-nodal/chordal-hyperplanes} below, all
discriminant hyperplanes are $\PGamma$-equivalent, so the discriminant locus
has image exactly $\H_{\D}$.  The same argument, using theorem~\ref{thm-chordal-maps-onto-divisor},
shows that the chordal locus has image $\H_c$.

Now we claim that $\widehat\moduli_s^f\to\ch^{10}$ is a local isomorphism.
Since its restriction to the preimage of $\ch^{10}-\H$ is an
isomorphism, $g:\widehat\moduli_s^f\to\ch^{10}$ is a proper modification of
$\ch^{10}$.  (See \cite[pp.~214--215]{grauert-remmert}.)  Since
$\ch^{10}$ is smooth, $g$ is a local isomorphism at a point $x$ of
$\widehat\moduli_s^f$ unless it crushes some divisor passing through $x$ to a  variety
of lower dimension.  Since $\moduli_0^f$ maps isomorphically to its
image, and the discriminant and chordal loci are mapped onto divisors,
no divisors are crushed.  Therefore the period map is everywhere a
local isomorphism.  Since it is a generic isomorphism, it is an
isomorphism.  It follows that it identifies $\moduli_0^f$ with
$\ch^{10}-\H$.  

(One can avoid the machinery of proper modifications by applying
Zariski's main theorem to 
\begin{equation}
\label{eq-foobar-3}
P\Gamma'\backslash\widehat\moduli_s^f\to P\Gamma'\backslash\ch^{10}\;,
\end{equation}
where $P\Gamma'$ is a torsion-free finite index subgroup of
$P\Gamma$.  Because this is a birational isomorphism of algebraic
varieties, and the target space is smooth, one can apply a version
of Zariski's main theorem,  theorem~3.20 of \cite{mumford}, to deduce that
\eqref{eq-foobar-3} is a local isomorphism, hence an isomorphism.  That
$\widehat\moduli_s^f\to\ch^{10}$ is an isomorphism follows.)

To prove the claim about $\widehat\moduli_{ss}$, it suffices to examine the
cusps of $\overline{\PGamma\backslash\ch^{10}}$.  The $D_4$ cusp is
the image of just one point of $\widehat\moduli_{ss}$, so the period map is
finite there.  Since the Baily-Borel compactification is a normal
analytic space, it is an isomorphism there.  The preimage of the $A_5$
cusp consists of the (classes of the) threefolds $T_{A,B}$ from
theorem~\ref{thm-GIT-analysis}\eqref{item-closed-orbits-of-semistable-threefolds}, together with their limiting point, the (class
of the) points from theorem~\ref{thm-GIT-analysis}\eqref{item-closed-orbits-of-semistable-12-tuples}.  These form  a
rational curve because they are parameterized by $4A/B^2$.

For the last claim, we observe that $\moduli_s$ is $\widehat\moduli_s$
minus the chordal locus.  This follows from a comparison of the GIT
analyses in theorem~\ref{thm-GIT-analysis} and \cite{allcock-threefolds}.  This makes the last claim
obvious. 
\end{proof}

We used the following lemma in the proof of the theorem.

\begin{lemma}
\label{lem-transitivity-on-nodal/chordal-hyperplanes}
Any two discriminant (resp. chordal) hyperplanes in $\ch^{10}$ are
$P\Gamma$-equivalent. 
\end{lemma}

\begin{proof}
Suppose $r\in \Lambda$ is a chordal root.  Since $\ip{r}{r}=3$ and
$\ip{r}{\Lambda}=3\E$, $\langle r\rangle$ is a summand of $\Lambda$.
Its orthogonal complement therefore has the same determinant as
$\Lambda_{10}=E_8^\E\oplus E_8^\E\oplus
\bigl(\begin{smallmatrix}0&\theta\\\thetabar&0\end{smallmatrix}\bigr)$,
namely $3^5$.  Also, $\theta(r^\perp)^*=r^\perp$, since $r^\perp\sset
\Lambda$.  There is a 
unique $\E$-lattice $L$ of determinant $3^5$ and signature $(9,1)$
satisfying $L\sset\theta L^*$, by lemma~2.6 of \cite{basak}.  (The
proof given in \cite{basak} uses the uniqueness of unimodular $\E$-lattices
with given indefinite signature, which  is theorem~7.1 of
\cite{allcock-reflection-groups}.)  Therefore $r^\perp\isomorphism
\Lambda_{10}$.  If $s$ is another chordal root then the same argument
shows that $\Lambda=\langle s\rangle\oplus \Lambda_{10}$.  So there is
an isometry of $\Lambda$ carrying $r$ to $s$.  This proves
transitivity on chordal roots.

We will prove three claims below.  (i) Any two index~$3$ sublattices of
$\Lambda_{10}$  are isometric; we write $\Lambda_{10}'$
for such a lattice.  (ii) For any nodal root $r\in \Lambda$,
$r^\perp\isomorphism \Lambda_{10}'$. (iii) There are exactly two
enlargements of $r^\perp\oplus\langle r\rangle$ to a copy of $\Lambda$
in which $r$ is a nodal root, and these are exchanged by negating
$r^\perp$ while leaving $r$ fixed.  Given the claims, a standard
argument shows that $r$ can be carried to any other nodal root $s$ by
an element of $\Gamma$.  Namely, by (ii), there is an isometry
$r^\perp\oplus\langle r\rangle\to s^\perp\oplus\langle s\rangle$
carrying $r$ to $s$.  This isometry carries the enlargement $\Lambda$
of $r^\perp\oplus\langle r\rangle$ to an enlargement $M$ of
$s^\perp\oplus\langle s\rangle$.  By (iii), $M$ is either $\Lambda$
itself or else is carried to $\Lambda$ by negating $s^\perp$.  The
result is an isometry of $\Lambda$ carrying $r$ to $s$.

Now we prove (i).  We will even show that $\aut \Lambda_{10}$ acts
transitively on the index~$3$ sublattices of $\Lambda_{10}$.
The key ingredient is a symplectic form on
$\Lambda_{10}/\theta \Lambda_{10}$.  For any $\E$-lattice $L$
satisfying $L\sset\theta L^*$, the $\F_3$-vector space $L/\theta L$
admits an antisymmetric pairing, defined as follows: if $v,w\in L$
have images $\bar v,\bar w\in L/\theta L$, then $(\bar v,\bar
w)\in\F_3$ is the reduction of $\frac{1}{\theta}\ip{v}{w}$ modulo $\theta$.
It is easy to check that if $v\in L$ has norm~$3$ and has inner
product $\theta$ with some element of $L$, then $\bar v$ does not lie
in the kernel of $({\cdot},{\cdot})$, and the triflections in $v$ act
on $L/\theta L$ as the transvections in $\bar v$.  Since
$\Lambda_{10}=\theta \Lambda_{10}^*$, the pairing on
$\Lambda_{10}/\theta \Lambda_{10}$ is nondegenerate.  Index~$3$
sublattices of $\Lambda_{10}$ correspond to hyperplanes in the
$\F_3$-vector space, so to prove transitivity of $\aut \Lambda_{10}$
on such subspaces, it suffices to show that it acts as
$\Sp(10,\F_3)$.  This is easy, since every root of $\Lambda_{10}$
gives a transvection.

Before proving (ii), we prove the weaker claim that for any nodal root
$r\in \Lambda$, 
$$
r^\perp\oplus\langle r\rangle
\ \isomorphism\ 
\Lambda_{10}'\oplus(3)
\ \isomorphism\ 
(3)\oplus E_8^\E\oplus E_8^\E\oplus(-3)\oplus(3)\;.
$$
Because $\ip{r}{\Lambda}=\theta\E$, $r^\perp\oplus\langle r\rangle$
has index~$3$ in $\Lambda$ and contains $\theta \Lambda$.  Therefore
it is completely determined by its image $S$ in $W=\Lambda/\theta
\Lambda$, which is $\bar r^\perp$.  The kernel of the pairing is
1-dimensional, coming from the first summand of \eqref{eq-inner-product-matrix}, and $\bar r$
is not in this kernel.  Therefore $\bar r^\perp$ is the preimage under
$W\to W/\ker W$ of a hyperplane in $W/\ker W$.  As we saw above, $\aut
\Lambda_{10}\sset \Gamma$ acts on $W/\ker
W=\Lambda_{10}/\theta\Lambda_{10}$ as $\Sp(10,\F_3)$, so its acts
transitively on hyperplanes in $W/\ker W$.  It follows that the
isomorphism class of $r^\perp\oplus\langle r\rangle$ is independent of
$r$, so it can be described by working a single example.  We take
$r=(0,\dots,0,1,\w)$ in the coordinates of \eqref{eq-inner-product-matrix}, and write $a$ for
$(1,0,\dots,0)$ and $b$ for $(0,\dots,0,1,\wbar)$.  Then $r^\perp$ is
$(3)\oplus\E_8^\E\oplus E_8^\E\oplus(-3)$, with $a$ spanning the
first summand and $b$ spanning the last.

Now we prove (ii).  It suffices to prove that if $s$ is a norm~3
vector of $N=\Lambda_{10}'\oplus(3)$, such that $\langle s\rangle$ is
a summand, then $s^\perp\isomorphism \Lambda_{10}'$.  In order to do
this we need to refer to the $\F_3$-valued symmetric bilinear form on
$\theta N^*/N$, got by reducing inner products modulo $\theta$.
Because $\theta(E_8^\E)^*=E_8^\E$, $\theta N^*/N\isomorphism\F_3^3$,
with a basis consisting of $\bar a$, $\bar b$ and $\bar r$, which are
the reductions modulo $N$ of $a/\theta$, $b/\theta$ and $r/\theta$.
The norms of $\bar a$, $\bar b$ and $\bar r$ are $1$, $-1$ and $1$.
Since $\langle s\rangle$ is a summand of $N$, $\ip{s}{N}=3\E$, so
$s/\theta\in\theta N^*$.  We write $\bar s$ for the image of
$s/\theta$ in $\theta N^*/N$, and observe that $\bar
s^2=\ip{s/\theta}{s/\theta}=1$.  Enumerating the elements of $\theta
N^*/N$ of norm~$1$, we find that $\bar s$ is one of $\pm\bar a$,
$\pm\bar r$, or $\pm\bar a\pm\bar b\pm\bar r$.  In every case there is
a $1$-dimensional isotropic subspace $S$ of $\theta N^*/N$ orthogonal
to $\bar s$.  Two examples: if $\bar s=\bar a$ then we can take $S$ to
be the span of $\bar b+\bar r$, and if $\bar s=\bar a+\bar b+\bar r$
then we can take $S$ to be the span of $\bar a+\bar b$.  We define
$N^+$ to be the preimage of $S$ in $\theta N^*$.  This is spanned by
$N$ and a vector of norm divisible by~$3$, whose inner products with
elements of $N$ are divisible by $\theta$.  Therefore
$\theta(N^+)^*\sset N^+$.  Now, the determinant of $\Lambda$ is $3^6$,
so the index~$3$ lattice $N$ has determinant $3^7$, so $N$'s index~$3$
superlattice $N^+$ has determinant $3^6$.  Also, $\ip{s}{N^+}=3$, because
$S\perp\bar s$.  Therefore $\langle s\rangle$ is a summand of $N^+$,
so its complement in $N^+$ has determinant $3^5$ and is therefore a
copy of $\Lambda_{10}$.  Since the complement of $s$ in $N$ has index three
in the complement in $N^+$, it is $\Lambda_{10}'$ by (i), proving
(ii). 

Finally we prove (iii).  Using $N$, $\bar a$, $\bar b$ and $\bar r$ as
above, we must determine the enlargements $M$ of
$\Lambda_{10}'\oplus\langle r\rangle$ that are copies of $\Lambda$ in
which $r$ is a nodal root.  Such an $M$ must satisfy $M\sset\theta
M^*$ and contain a vector having inner product $\theta$ with $r$.
Such an $M$ is the preimage of an isotropic line in $\theta N^*/N$
which is not orthogonal to $\bar r$.  There are just two such
subspaces, the spans of $\bar r-\bar b$  and $\bar r+\bar b$.  The
first enlargement is $\Lambda$ itself.   The two enlargements are
exchanged by negating $r^\perp$.  This proves (iii).
\end{proof}

\section{The Monodromy Group and the Hyperplane Arrangement}
\label{sec-group-and-arrangement}

We have already discussed several facets of the action of $\Gamma$ on
various objects associated to $\Lambda$.
Theorem~\ref{thm-exactly-two-cusps} shows that $\Gamma$ acts with two
orbits on primitive null vectors, corresponding to the two boundary
points of $\overline{\PGamma\backslash\ch^{10}}$, and
lemma~\ref{lem-transitivity-on-nodal/chordal-hyperplanes} shows that
$\Gamma$ acts with two orbits on roots of $\Lambda$, corresponding to
discriminant and chordal hyperplanes.  In this section, we gather a few
more results of this flavor, and determine the image of the monodromy
representation.

\begin{theorem}
\label{thm-linear-monodromy-group}
For $F\in\forms_0$,
the  monodromy representation $\pi_1(\forms_0,F)\to\Gamma(V)$ has
image consisting of all isometries of $\Lambda(V)$
with determinant a cube root of unity.
\end{theorem}

\begin{proof}
We first show that the image of
$P\rho\bigl(\pi_1(\forms_0,F)\bigr)\sset\PGamma(V)$ is 
all of $\PGamma(V)$.  This amounts to the connectedness of
$\framed_0$, which is the same as the connectedness of $P\framed_0$,
which is the same as the connectedness of $P\framed_s$, which is the
same as the connectedness of $\widehat\moduli_s^f$, which by theorem~\ref{thm-MAIN-THEOREM} 
is the same as the
connectedness of $\ch^{10}$, which is obvious.

Now it suffices to show that $\rho\bigl(\pi_1(\forms_0,F)\bigr)$
contains $\wbar 
I$ and does not contain $-I$.  To see the latter, observe that
$\pi_1(\forms_0)$ is generated by meridians, which $\rho$ maps to
$\w$-reflections, so every element of $\pi_1(\forms_0)$ acts on
$\Lambda(V)$ by an isometry with 
determinant a cube root of~1.  (Another way to see
that $-I$ is not in the monodromy group is to use the ideas of the proof of
theorem~\ref{thm-main-theorem-smooth-case}: the negation map of
$H^4_0(V;\Z)$ does not extend to an isometry of $H^4(V;\Z)$ fixing
$\eta(V)$.)

Now we prove that $\wbar I$ lies in the monodromy group.
Consider the cubic threefold $T$ defined by 
$$
F=x_0x_3^2-x_0x_2x_4+x_1^2x_4+x_2^2x_3\;.
$$
By the method of \cite[\S2]{allcock-threefolds}, it has an $A_7$ singularity at
$[1{:}0{:}0{:}0{:}0]$, an $A_4$ singularity at $[0{:}0{:}0{:}0{:}1]$, and no other
singularities.   Using lemmas~\ref{lem-meridians-act-by-complex-reflections} and~\ref{lem-local-monodromy-noncyclic} 
shows that
$\rho(\pi_1(\forms_0,F))$ contains an image of the product $B_8\times B_5$ of the
8-strand and 5-strand braid groups, with standard generators
$a_1,\dots,a_7$ and $b_1,\dots,b_4$ mapping to the $\w$-reflections in
roots $r_1,\dots,r_7$ and $s_1,\dots,s_4$ of $\Lambda(V)$,
no two of which are proportional.
The method described before theorem~\ref{thm-finite-reflection-group-facts}  shows that the $r_i$
(resp. $s_i$) span a 
nondegenerate $\E$-lattice $R$ (resp. $S$) of dimension~7 (resp.~4).
Since $R$ and $S$ are orthogonal and nondegenerate, they meet only at
$0$, so $\Lambda(V)=R\oplus S$ up to finite index.
A generator $c_8$ (resp. $c_5$) for the center of $B_8$ (resp. $B_5$)
is $(a_1\dots a_7)^8$ (resp. $(b_1\dots b_4)^5$), so it acts on $R$
(resp. $S$) as a scalar of determinant $\w^{56}$ (resp. $\w^{20}$),
i.e., as the scalar $\wbar$ (resp. $\wbar$).  Therefore $c_5c_8$ acts
on $\Lambda(V)$ by $\wbar I$.
\end{proof}

\begin{theorem}
\label{thm-how-hyperplanes-meet}
(i) No two chordal hyperplanes meet in $\ch^{10}$.  
(ii) If a discriminant hyperplane and a chordal hyperplane meet in $\ch^{10}$, then they are
orthogonal. 
\end{theorem}

\begin{proof}
If $r$ and $s$ are a nodal and a chordal root, then
$s^2=r^2=3$ and $\ip{r}{s}$ is divisible by~$3$.  If $\ip{r}{s}=0$
then they are orthogonal, and otherwise the span of $r$ and $s$ has
inner product matrix which is not positive-definite.  In this case,
$r^\perp$ and $s^\perp$ do  not meet in $\ch^{10}$.  This proves
(ii). 

Now we prove (i).  Suppose $r$ and $s$ are non-proportional chordal
roots.  By the argument for (ii), if $r^\perp$ and $s^\perp$ meet,
then $r\perp s$.  So we must show that $r\perp s$ is
impossible.  By lemma~\ref{lem-transitivity-on-nodal/chordal-hyperplanes}, we may take $s$ to be $(1,0,\dots,0)$ in
the coordinates of \eqref{eq-inner-product-matrix}.  If $r\perp s$ then $r$ lies in
$\Lambda_{10}$.  Now, 
$r$ is primitive because its norm is~$3$, and by $\theta
\Lambda_{10}^*=\Lambda_{10}$, $r$ makes inner product $\theta$ with some element
of $\Lambda_{10}$.  Therefore $r$ cannot be a chordal root. 
\end{proof}

\begin{theorem}
\label{thm-nodal-and-chordal-loci-commensurable}
The discriminant and chordal loci of $\widehat\moduli_s$ have finite branched
covers which are isomorphic.
\end{theorem}

\begin{proof}
Consider the nodal root $r=(0,\dots,0,1,\w)$ and chordal root
$s=(1,0,\dots,0)$, in the notation of \eqref{eq-inner-product-matrix}.  The discriminant
(resp. chordal) locus of $\widehat\moduli_s$ is the image in
$P\Gamma\backslash\ch^{10}$ of $r^\perp$ (resp. $s^\perp$).  We have
$$
r^\perp\oplus\langle r\rangle
=\langle s\rangle\oplus E_8^\E\oplus E_8^\E\oplus(-3)\oplus\langle
r\rangle\;.
$$
Let $\Gamma'$ be the subgroup of $\Gamma$ preserving the sublattice
$r^\perp\oplus\langle r\rangle$, which has index $(3^{10}-1)/2=29524$ in
$\Gamma$.  The images of $r^\perp$ and $s^\perp$ in
$P\Gamma'\backslash\ch^{10}$ are isomorphic because
$r^\perp\oplus\langle r\rangle$ admits an isometry exchanging $r$ and $s$.
\end{proof}

\begin{remark}
The same argument proves theorem~3 in Kond\=o's work
\cite{kondo:genus-4} on genus 4 curves, and a similar result relating
nodal and hyperelliptic hyperplanes in his uniformization of the
moduli space of genus 3 curves by the 6-ball.
\end{remark}

\end{document}